\documentclass{article}
\usepackage{amssymb}
\usepackage{amsmath}
\usepackage{enumerate}
\usepackage{amscd}
\usepackage{graphics}

\usepackage{latexsym}
\usepackage{amssymb}
\usepackage{amsmath}
\usepackage{epsfig}
\usepackage[latin1]{inputenc}

\def\r{\mathbb{R}}
\def\n{\mathbb{N}}
\def\c{\mathbb{C}}

\def\s{\mathbb{S}}
\def\d{\mathbb{D}}
\def\l{\mathbb{R}_1}
\def\z{\mathbb{Z}}

\def\h{\mathbb{H}}
\def\u{\mathbb{U}}

\def\mb{\mathcal{M}}

\newenvironment{proof}{\trivlist
\item[\hskip\labelsep{\em Proof}\,:]}{\hfill{$\Box$}\endtrivlist}

\title{\huge On the uniqueness of the helicoid and Enneper's surface in the Lorentz-Minkowski space $\l^3$}
\author{\Large Isabel Fernandez and Francisco J. Lopez  \thanks{Research of both authors partially
supported by
MCYT-FEDER grant number MTM2004-00160.}}

\def\N{\mathcal N}

\newtheorem{lemma}{Lemma}[section]
\newtheorem{remark}{Remark}[section]
\newtheorem{theorem}{Theorem}[section]
\newtheorem{proposition}{Proposition}[section]\newtheorem{corollary}{Corollary}[section]

\newtheorem{definition}{Definition}[section]

\begin{document}
\maketitle

\begin{abstract} In this paper we deal with the uniqueness of the Lorentzian helicoid and Enneper's surface among properly embedded maximal surfaces with lightlike boundary of  mirror symmetry in the Lorentz-Minkowski space $\l^3.$
\end{abstract}

\section{Introduction}

\begin{figure}[htpb]
\begin{center}\includegraphics[width=.7\textwidth]{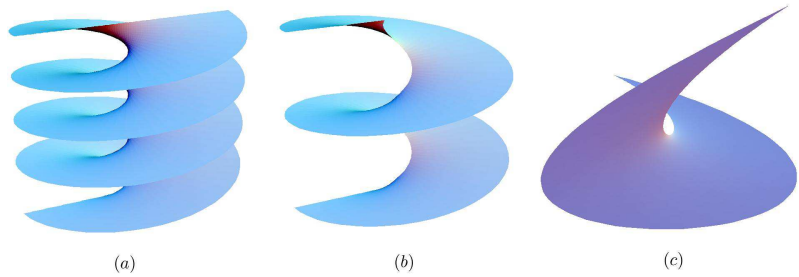}
\caption{$(a)$ the helicoid ${\cal H}_0;$ $(b)$ the Lorentzian
helicoid ${\cal H};$ $(c)$ Enneper's surface
$E_1$}\label{fig:helicoides} \end{center}
\end{figure}

The helicoid ${\cal H}_0:=\{(x,y,t)\in \r^3 \;:\; x\tan (t)=y\}$ was  first discovered by Jean Baptiste Meusnier in 1776. After the plane and the catenoid, is the third minimal surface  in Euclidean space $\r^3$ to be known. The helicoid is generated by spiraling a horizontal straight line along a vertical axis, and so, it is a ruled surface which is also foliated by helices (its name derives from this fact). As shown in Figure \ref{fig:helicoides}, it is shaped like the Archimedes' screw, but extends infinitely in all directions, see Figure \ref{fig:helicoides}, $(a).$

In analogy with minimal surfaces in $\r^3,$ a maximal surface in
3-dimensional Lorentz-Minkowski space
$\l^3=(\r^3,dx^2+dy^2-dt^2)$ is a surface which is spacelike (the
induced metric  is Riemannian) and whose mean curvature vanishes.
Maximal surfaces represent local maxima for the area functional, have conformal Gauss map and admit
a Weierstrass type representation (see equation
(\ref{eq:wei})).  Besides of their mathematical interest, they have a
significant importance in classical Relativity (see \cite{marsden-tipler}).

The relative complement of the rigid circular cylinder $C=\{(x,y,t) \,:\; x^2+y^2\leq 1\}$ in ${\cal H}_0$ is a spacelike surface when viewed in $\l^3,$ and consists of two congruent (in the Riemannian and Lorentzian sense) simply connected domains of ${\cal H}_0$ bounded  by a  lightlike  helix. The {\em Lorentzian helicoid} ${\cal H}$ is defined to be the closure  of the connected component of ${\cal H}_0-C$ containing $(1,0,1)$ in its boundary, see Figure \ref{fig:helicoides}, $(b).$

Amazingly, $\mbox{Int}({\cal H})$ is a maximal surface. As a matter of fact, O. Kobayashi \cite{kobayashi1} proved that the Lorentzian helicoid and spacelike planes  are the only maximal surfaces which are also minimal surfaces with respect to the Euclidean metric on the ambient 3-space.

The geometry of the Lorentzian helicoid is somehow concentrated in its boundary. If $X:\mb \equiv
\{z \in \c \,:\,\mbox{Im}(z) \geq 0\}  \to \l^3$ is a conformal embedding of ${\cal H},$  $\partial(\mb)$ is a  integral
curve of the {\em weighted} gradient $\frac{1}{(t\circ {\cal
N})^2}\nabla (t \circ X),$ where ${\cal N}$ is the Lorentzian
Gauss map  of $X$ and $\nabla$ is computed with respect to the
intrinsic metric induced by $\l^3$ (the factor
$\frac{1}{(t\circ {\cal N})^2}$ just controls the singularities
of $\nabla (t \circ X)$ along $\partial (\mb)$).  Equivalently, the immersion
folds back at $\partial(\mb),$ that is to say, $X$ extends
harmonically to the double of $\mb$ being invariant under the
mirror involution. Maximal surfaces with regular lightlike
boundary and satisfying this symmetry property are said to have
{\em lightlike boundary of mirror symmetry}, and will be written 
 $^*$maximal surfaces.

Another interesting example of properly embedded  $^*$maximal
surface with connected boundary in $\l^3$ is the so called Lorentzian Enneper surface
$E_1:=\{(x,y,t) \;:\; 32(y-t)^3 - 3(y+t) + 24(y-t)x = 0\},$ see
Figure \ref{fig:helicoides}, $(c).$  Unlike the Lorentzian
helicoid, the Lorentzian Enneper surface has finite {\em rotation number}, that is to say, 
the change of the tangent angle along the orthogonal projection over $\{t=0\}$ of its boundary is finite.

In this paper we will take interest in $^*$maximal surfaces with slightly controlled asymptotic behavior. To be more precise, a properly immersed maximal surface $\mb$ in $\l^3$ is said to be {\em asymptotically weakly spacelike}, or simply  a $\omega$-maximal surface,  if $\l^3-\mb$ contains  an {\em affinely spacelike} arc (i.e., a proper arc $\alpha \cong [0,1[$ in $\l^3$ lying in a closed spacelike wedge\footnote{A closed wedge in $\l^3$ is said to be spacelike if it is foliated by spacelike half planes with the same edge.} of $\l^3$). If in addition $\mb$ is $^*$maximal, we simply say that it is $\omega^*$-maximal. In the context of $^*$maximal surfaces with connected boundary, this apparently mild condition let us control the geometry of the homothetical blow-downs of $\mb.$ Moreover, it is automatically satisfied by surfaces with finite rotation number and by surfaces admitting gradient estimates (or more generally, being metrically complete) far from the boundary, like the helicoid. See subsection \ref{subsec:w*} for more details.

Recently, ${\cal H}_0$ has been characterized by W. H. Meeks III
and H. Rosenberg  \cite{meeks-ros} as the unique properly
embedded  non flat simply connected minimal surface in $\r^3.$
Likewise, J. Perez \cite{joa} has proved that half of the Enneper
minimal surface is  the only properly embedded non flat oriented
stable minimal surface bounded by a straight line and having
quadratic area growth. Somehow, this paper is devoted to obtain a
Lorentzian compilation of both Riemannian theorems.

We have proved the following:
\begin{quote}
{\bf Theorem I:} {\em The only properly embedded  $^*$maximal surface with connected boundary of finite rotation number in $\l^3$ is the Lorentzian Enneper surface.}

{\bf Theorem II:} {\em The only properly embedded  $\omega^*$-maximal  surface with connected boundary of infinite rotation number in $\l^3$ is the Lorentzian helicoid.}
\end{quote}

It remains open whether Theorem II can be extended to $^*$maximal surfaces.

The required theoretical background includes classical Calabi's
theorem \cite{calabi} (see also Cheng-Yau work \cite{c-y}) about
complete maximal surfaces, and some basic existence and
regularity properties of area maximizing surfaces in the Lorentz
Minkowski space  $\l^3,$ mainly proved by Bartnik and Simon in
\cite{bar-sim}.

In a first step, we obtain some regularity theorems and parabolicity criteria for maximal graphs,  and use these results to control the asymptotic behavior of maximal graphs over planar wedges. Among other things, we study under what geometrical conditions the homothetical blow-downs of a such graph do not converge to an angular region of the light cone. 
Taking advantage of this analysis and Calabis' theorem, we can derive an elementary Colding-Minicozzi theory  \cite{colding-mini} and prove that any homothetical blow-down of a properly embedded  $\omega^*$-maximal surface ${\cal S}$ with connected boundary is a plane viewed as a degenerated multigraph with a singular point. This means that any leaf of the blow-down sequence converges in the ${\cal C}^1$-topology outside the singularity to a plane $\Sigma_\infty$ (the blow-down plane) depending neither on the leaf nor the homothetical blow-down. 
Finally, and using some ideas by W.H. Meeks and Rosenberg  in \cite{meeks-ros}, we deduce that the Gauss map of ${\cal S}$ omits the normal direction of $\Sigma_\infty,$ and that any plane parallel to $\Sigma_\infty$ intersects ${\cal S}$ into a single arc. This reasoning strategy requires of a finiteness theorem for maximal graphs with planar boundary, whose proof has been deeply inspired by P. Li and J. Wang  work \cite{li-wang}. The natural dichotomy between  spacelike and lightlike blow-down plane leads to ${\cal S}={\cal H}$ and ${\cal S}=E_1,$ respectively.

The paper has been laid out as follows:

In Section \ref{sec:prelim} we introduce some terminology and  background material. A detailed description of the basic examples  (Helicoid, Enneper's surface and conjugate surfaces) is given in Section \ref{sec:ejemplos}. Section \ref{sec:parabo} is devoted to obtaining some parabolicity criteria for maximal surfaces. In Section \ref{sec:grafos} we deal with the geometry of maximal graphs, specially those over wedge-shaped regions. We also prove the Li-Wang type  finiteness theorem for maximal graphs. The deepest results are contained in Section \ref{sec:main}, which has been devoted to the global geometry of properly embedded $\omega^*$-maximal surfaces  with connected boundary. We construct the blow-down multigraph and prove the transversality of the surface and the  blow-down limit plane. Finally, in  Section \ref{sec:uniq} we prove the uniqueness theorems.

\begin{quote}
{\bf Acknowledgment:} {\small  The authors are deeply indebted to Prof. Rabah Souam for many useful conversations about maximal surfaces and the Lorentzian helicoid.}
\end{quote}

\section{Notations and Preliminaries} \label{sec:prelim}

As usual, $\overline{\c}=\c \cup \{\infty\},$ $\u=\{z \in \c \,:\, \mbox{Im}(z)>0\},$ $\d=\{z \in \c \,:\,|z|<1\}$ and $[-\infty,+\infty]=\r \cup \{-\infty,+\infty\}.$ We make the convention $x\pm\infty=\pm\infty,$ for all $x \in \r.$ If $I \subset \r$ is an interval, we call $|I|$ as its Euclidean length.

The Euclidean metric and norm in  $\r^n$  will be denoted  by
$\langle,\rangle_0$ and $\|\cdot \|_0,$ respectively, $n \geq 2.$
The origin in $\r^n$ will be written as $O.$  Given $W_1,$ $W_2
\subset \r^n$  we denote by  $\mbox{d}(W_1,W_2)$ and
$\mbox{d}_H(W_1,W_2):=\sup\{\sup\{\mbox{d}(w,W_j)\;:\;w \in W_1
\cup W_2\} \;:\; j=1,2\}$ the Euclidean and Hausdorff distance
between $W_1$ and $W_2,$ respectively.

A smooth divergent arc $\alpha(u):[0,+\infty[\to \r^n$ is defined to be {\em sublinear with direction} $v\in \s^n$ if $\lim_{u \to +\infty}\alpha'(u)=v,$ where $u$ is the arclength parameter of $\alpha$ with respect to $\langle,\rangle_0,$ $n \geq 2.$

We call $\l^3$ the three dimensional Lorentz-Minkowski space $(\r^3\equiv \r^2 \times \r,\langle , \rangle),$ where as usual $\langle(x_1,t_1) ,(x_2,t_2)\rangle=\langle x_1, x_2\rangle_0 -t_1 t_2,$ and write $\|(x,t)\|^2:=\|x\|_0^2-t^2.$ A vector ${ v} \in \r^3- \{ {(0,0,0)} \}$ is said to be spacelike,
timelike or lightlike if
$\|v\|^2>0,$  $\|v\|^2<0$ or $\|v\|^2=0$,
respectively. The vector
${(0,0,0)}$ is spacelike by definition.  A smooth curve in $\l^3$ is defined to be spacelike, timelike or lightlike if all its tangent vectors are spacelike, timelike or lightlike, respectively. A plane in $\l^3$ is spacelike, timelike or lightlike if the induced metric is Riemannian, non degenerate
 indefinite or degenerate, respectively. The spacelike plane $\{t=0\}$ will be denoted by $\Pi_0.$ We often use the identification $\Pi_0 \equiv \r^2$ given by $(x,0) \equiv x.$

A closed wedge $W$ of $\l^3$ is said to be spacelike if any half plane contained in $W$ is spacelike. Proper arcs lying in a spacelike wedge are said to be {\em affinely spacelike}. Given a spacelike plane $\Sigma \subset \l^3,$ $\pi_{_{\Sigma}}:\l^3 \to \Sigma$ will denote the Lorentzian orthogonal projection. If $\Sigma=\Pi_0$ we simply write $\pi$ instead of $\pi_{_{\Pi_0}},$  and in this case $$\pi((x,t)):=x, \quad (x,t) \in \l^3.$$

For any $p=(x_0,t_0) \in \l^3,$ we denote by ${\cal C}_p:=\{x \in \l^3: \|x-p\|^2=0\}$ the light cone with vertex at $p,$ and label ${\cal C}_p^+:={\cal C}_p \cap \{t \geq t_0\}$ and  ${\cal C}_p^-:={\cal C}_p \cap \{t \leq t_0\}.$ We also set $\mbox{Int}({\cal C}_p):=\{x \in \l^3: \|x-p\|^2< 0\}$ and $\mbox{Ext}({\cal C}_p):=\{x \in \l^3: \|x-p\|^2> 0\},$ and likewise define  $\mbox{Int}({\cal C}^+_p):=\mbox{Int}({\cal C}_p)\cap \{ t >t_0\}$ and  $\mbox{Int}({\cal C}^-_p):=\mbox{Int}({\cal C}_p)\cap \{t <t_0\}.$

A  spacelike arc $c \subset \l^3$ is said to be an {\em upward} (resp., {\em downward}) {\em lightlike ray} if,  up to removing a compact subarc, $\pi(c)$ is a closed half line and there exists $p \in \pi^{-1}(\pi(c))$  such that $\lim_{x \in c \to \infty} \mbox{d}(x,l_c)=0,$ where $l_c$ is the lightlike half line in ${\cal C}_p^+ \cap \pi^{-1}(\pi(c))$ (resp., in ${\cal C}_p^- \cap \pi^{-1}(\pi(c))$) with initial point $p.$

As usual, open connected subsets of  manifolds are called {\em domains} and their closures {\em regions}.  Throughout this paper we will deal with regions and domains of surfaces, namely $\Omega,$ with regular enough boundary. In the most cases that last means that $\partial (\Omega)$ is piecewise smooth. If $\overline{\Omega}$ lies in a Riemannian surface,  it suffices to require that $\partial (\Omega)$ is  ${\cal C}^0$ and locally Lipschitzian functions in $\mbox{Int}(\Omega)$ extend continuously to $\partial(\Omega).$ 

If $S$ is a manifold and $f:S \to \r$ is a function, the expression $\lim_{x \in S \to \infty} f(x)=[-\infty,+\infty]$ means that $\lim_{n \to \infty} f(x_n)=L$ for any divergent sequence $\{x_n\}_{n \in \n} \subset S.$

Let ${\cal R}^*:= \{(z,w) \in (\c-\{0\})\times {\c}\;:\; e^w=z\}$
denote the Riemann surface of $\log(z)$ endowed with the
Riemannian metric $|dz|^2.$ The function $w:{\cal R}^* \to {\c}$ is a
biholomorphism and $z:{\cal R}^* \to \c^*:={\c}-\{0\}$ is the
isometric  universal covering of the Euclidean once punctured
plane. The argument function is given by $\arg:{\cal R}^*\to \r,
\quad \arg=\mbox{Im}(w).$ For convenience, we add an extra point
$[0]$ to ${\cal R}^*,$ define $z([0])=0,$ and endow  ${\cal R}:={\cal R}^* \cup
\{0\}$ with the smallest topology containing the one of ${\cal R}^*$ and making   $z:{\cal R} \to \c$ continuous.

Let $W \subset {\cal R}$ be a proper subset  homeomorphic to $\overline{\d}-\beta,$ where $\beta$ is a non empty connected subset of $\partial(\d).$ $W$ is defined to be a (generalized) {\em  wedge} if $\partial(W)$ is smooth outside a compact subset $C \subset {\cal R},$ and for any proper Jordan arc $\alpha \cong [0,1[$ in ${\cal R}$ contained in $\partial(W)-C,$  either $\theta_\alpha:=\lim_{x \in \alpha \to \infty} \arg(x)=\pm \infty$ or $z(\alpha)$ is a
planar sublinear arc (hence   $\theta_\alpha \in \r$). If $\alpha_1,$ $\alpha_2$ are two such arcs in $\partial (W)-C,$   we set $\theta:=|\theta_{\alpha_2}-\theta_{\alpha_1}|\in [0,+\infty]$ the angle of $W.$ In case  $\partial(W)=\{[0]\}$ (i.e., $W={\cal R}$) or  $\partial(W)$ consists of a divergent Jordan arc with initial point $[0],$ $W$ is defined to have infinite angle. The wedge  $\arg^{-1}([-\theta,\theta])\cup \{[0]\}$  will be denoted by $W_\theta,$ $\theta \in [0,+\infty].$  When $z|_W:W \to z(W)$ is one to one,  $W$ and $z(W)\subset \c\equiv \r^2$ will be identified. Moreover, regions in $\r^2$ defined by translating  wedges of angle $<2\pi$ will be also named wedges. 

In the sequel,  $\mb$ will denote a differentiable surface, possibly with non empty regular enough 
boundary.

A continuous map $X:\mb\to\l^3$ is said to be {\em pseudo spacelike} (acrostically, PS) if  for any $p \in \mb$ there is an open
neighbourhood $U$ of $p$ in $\mb$ such that
$\|X(p_1)-X(p_2)\|\geq 0,$ for any $p_1,$ $p_2 \in U.$ If in
addition $\pi \circ X$ is a local embedding\footnote{By the
Domain Invariance Theorem, this simply means that $\pi \circ X$
is {\em locally injective}.},  then $X$ is defined to be a {\em
pseudo spacelike immersion}.  If $X:\mb\to\l^3$ is pseudo
spacelike immersion and  $\pi \circ X$ is one to one, then
$X(\mb)$ is said to be a pseudo spacelike graph over $\pi(X(\mb))
\subset \Pi_0.$

Let $X$ be a proper PS immersion of a symply connected surface $\mb$ into $\l^3.$ $X$ is said to be a PS {\em multigraph} of angle $\theta \in ]0,+\infty[$ if there is an open disc $D$ centered at $O$
such that $(\pi \circ X)^{-1}(\overline{D})$ is compact, $Y:\mb_0 \to {\cal R}$ is an embedding where $\mb_0=\mb-(\pi \circ X)^{-1}(\overline{D}),$  $z \circ Y=\pi \circ (X|_{\mb_0})$  (here we have identified $\Pi_0$ and the $z$-plane), and $W:=Y(\mb_0)$ is a wedge of angle $\theta.$

Let $Z$ be a simply connected topological space and consider a proper continuous map $X:Z \to \l^3.$ Suppose that $(\pi \circ X)^{-1}(O)$ is either empty or consists of a single point and call $\mb_0=\mb-(\pi \circ X)^{-1}(O).$ $X$ is said to be a PS {\em multigraph of infinite angle} if $X|_{\mb_0}$ is a PS immersion and there is an embedding $Y:\mb \to {\cal R}$ such that $z \circ Y=\pi \circ X|_{\mb}$ and  $W:=Y(\mb)$ is a wedge of infinite angle (obviously, $Y((X \circ \pi)^{-1}(O))=[0]$ provided that  $(X \circ \pi)^{-1}(O)$ is a point).

In any case,  $u=t \circ (X \circ Y^{-1}):\mb_0 \to \r$ is locally Lipschitzian with Lipschitz constant $1,$ i.e.,  $\nabla u$ is well defined in the weak sense and $\|\nabla u\|_0 \leq 1.$

Set  $G=\{(x,u(x))\;:\; x \in \Omega\}$ a PS graph over a  domain $\Omega \subset \r^2,$  and call $d_\Omega$  the inner metric in $\Omega$ induced by $\langle,\rangle.$ The PS condition again gives $|u(x)-u(y)| \leq d_\Omega(x,y)$ for all $x,$ $y \in \Omega.$ Thus,  if $\Omega$ is starshaped with center $x_0,$
\begin{equation} \label{eq:acausal}
G-\{(x_0,u(x_0))\} \subset \mbox{Ext}({\cal C}_{(x_0,u(x_0))}).
\end{equation}

 \begin{lemma} \label{lem:basico} Let $G$ be a PS graph over a convex region $R \subset \r^2.$  The following statements hold:

\begin{enumerate}[(a)]
\item If  $G$ contains a lightlike straight line $l,$ then $G$ lies in the lightlike plane $\Sigma_0$ containing $l.$
\item If $l_1=[p_1,p_2]$ and $l_2$ are lightlike segments in $G$ such that $]p_1,p_2[\cap l_2 \neq \emptyset$ then $l_1$ and $l_2$ lie in the same lightlike straight line.
\end{enumerate}
\end{lemma}
\begin{proof} To check $(a),$ use equation (\ref{eq:acausal}) and observe that $G \subset \cap_{x \in l} \overline{\mbox{Ext}({\cal C}_{x})}=\Sigma_0.$

To prove $(b),$ suppose up to a translation that $O\in ]p_1,p_2[ \cap l_2$ and consider the dilated graphs $G_n:=n \cdot G,$ $n \in \n.$ By Ascoli-Arzela Theorem, the sequence $\{G_n\}_{n \in \n}$ converges uniformly on compact subsets to a PS graph $G_\infty$ over a convex region containing the lightlike straight line  $l_0$ determined by $l_1.$ From $(a),$ $G_\infty$ lies in a lightlike plane, hence $l_2$ lies in $l_0$ too. \end{proof}

A smooth immersion $X:\mb \longrightarrow \l^3$ is said to be {\em spacelike } (and  $X(\mb)$ a spacelike surface in $\l^3$) if
the tangent plane at any  point is spacelike, that is to say, if the induced metric $ds^2:=X^*(\langle,\rangle)$ on $\mb$ is Riemannian. In this case,  the Gauss map ${\cal N}$ of $X$ is well defined and takes values in the Lorentzian sphere $\h^2:=\{x \in \l^3 \;;\; \langle x,x \rangle=-1\}.$
If we attach to $\mb$ the conformal structure induced by $ds^2,$ $\mb$ becomes a Riemann surface and $X$ a conformal spacelike immersion. It is easy to see that spacelike immersions are PS immersions.
%

Let $\mb$ be a Riemann surface, and let $S_X$ be a closed subset with empty interior (usually, a family of curves and points). A smooth map $X:\mb \to \l^3$ is said to be a {\em conformal spacelike immersion with singular set} $S_X$ (and $X(\mb)$ a spacelike surface with singular set $X(S_X)$) if $X^*(\langle,\rangle)=\lambda ds_0^2,$ where $X^*(\langle,\rangle)$ is the pull back metric of $\langle,\rangle,$ $ds_0^2$ is a conformal Riemannian metric on $\mb$ and $\lambda$ is a  function vanishing on $S_X$ and being positive on $\mb-S_X.$ A singular point $p \in S_X$ is said to be a {\em lightlike singularity} of $X$ if  $\lim_{q \in \mb-S_X \to p} {\cal N}(q)=\infty,$ where ${\cal N}:M-S_X\to\h^2$ is the Lorentzian Gauss map of $X|_{M-S_X}.$ If in addition $dX_p \neq 0,$ $p$ is said to be a {\em regular lightlike singularity}. See the papers \cite{um-ya,f-s-u-y} for a good setting about singularities.

A spacelike immersion $X:\mb \to \l^3$ is said to be {\em maximal} (and $X(\mb)$ a maximal surface) if its mean curvature vanishes.
A conformal maximal immersion $X:\mb \longrightarrow \l^3$  has harmonic coordinate functions and admits a Weierstrass type  representation $(g,\phi_3):$
\begin{equation} \label{eq:wei}
X=\mbox{Real} \int (\phi_1,\phi_2,i \phi_3),
\end{equation}
where $g$ is a meromorphic function (the meromorphic Gauss map) and $\phi_1=\frac{1}{2}(1/g-g) \phi_3,$ $\phi_2=\frac{i}{2}(1/g+g)\phi_3$ and $\phi_3$ are holomorphic 1-forms without common zeroes in $\mb.$
Recall that $g=\mbox{st} \circ {\cal N},$ where ${\cal N}:\mb \to \h^2$ is the Gauss map of $X$ and $\mbox{st}:\h^2 \to \overline{\c}$ is the   Lorentzian stereographic projection given by $\mbox{st}:{\h}^2 \to \d,$ $\mbox{st}((x_1,x_2,t))=(\frac{x_2}{t-1},\frac{x_1}{1-t}).$ We also set  $\mbox{st}_0:{\cal C}_0-\{(0,0,0)\} \to \{z \in \c \;:\; |z|=1\},$ $\mbox{st}_0((v_1,v_2,v_3)):=\frac{1}{v_3}(v_2,-v_1),$ and observe that $\mbox{st}_0(v)=\lim_{n \to \infty} \mbox{st}(w_n),$ provided that $\{w_n\}_{n \in \n} \subset \h^2$  and $\{\frac{w_n}{\|w_n\|_0}\}_{n \in \n} \to  \frac{v}{\|v\|_0}.$ 

For more details about the  Weierstrass representation of maximal surfaces see \cite{kobayashi}.

\begin{definition}
A map $X:\mb \to \l^3$ is defined to be a  conformal maximal immersion with singular set $S_X$ if $X$ is a conformal spacelike immersion with singular set $S_X$ and $X|_{\mb-S_X}$ is maximal. We also say that $X(\mb)$ is a maximal surface with  singular set $X(S_X).$
\end{definition}
The Weierstrass data $(g,\phi_3)$ of a  conformal maximal immersion $X:\mb \to \l^3$ with  singularities are well defined on $\mb,$ and since the intrinsic metric is given by $ds^2=\frac{1}{2} \left(1/{|g|}-|g|\right)^2 |\phi_3|^2,$ then $S_X$ coincides with the analytical set $|g|^{-1}(1) \cup |\phi_3|^{-1}(0).$ If in addition every singular point is regular and lightlike,  then  $\phi_3$ never vanishes and $S_X=|g|^{-1}(1).$
If $\mb_0$ is a connected component of $\mb-S_X,$ its Gaussian image ${\cal N}(\mb_0)$ lies in either $\h^2_+:=\h^2 \cap \{t >0\}$ or $\h^2_-:=\h^2 \cap \{t <0\}.$ We usually choose the orientation  in such a way that ${\cal N}(\mb_0) \subset \h^2_-,$ and so $g(\mb_0)\subset \d.$ In this case $g$ is said to be the holomorphic Gauss map of $X|_{\mb_0}.$

Let $\mb$ be a Riemann surface with analytical boundary. The mirror and double surface of $\mb$ with respect to $\partial(\mb)$ will be denoted by $\mb^*$ and $\hat{\mb}:=\mb \cup \mb^*,$  respectively. Recall that, up to natural identifications, $\partial(\mb)=\partial(\mb^*)=\mb \cap \mb^*=\{p \in \mb \;:\; J(p)=p\},$ where $J:\hat{\mb} \to \hat{\mb}$ is the antiholomorphic involution  mapping each point $p \in \mb$ into its mirror image $p^* \in \mb^*.$

\begin{definition} \label{def:simetria}
Let $\hat{X}:\hat{\mb} \to \l^3$ be a conformal maximal immersion with regular lightlike singularities. The map $X:=\hat{X}|_{\mb}$ is said to be a  conformal maximal immersion with lightlike boundary of mirror symmetry (or simply, a  conformal $^*$maximal immersion) if $S_{\hat{X}}=\partial(\mb)$ and $\hat{X} \circ J=\hat{X}.$ We also say that $X(\mb)=\hat{X}(\hat{\mb})$ is a $^*$maximal surface. 

In terms of the Weierstrass data $(g,\phi_3)$ of $\hat{X},$ $X$ is $^*$maximal if and only if:
\begin{equation} \label{eq:simetria}
\overline{g} \cdot  (g\circ J)={1},\quad J^*(\phi_3)=-\overline{\phi_3}.
\end{equation}
If $X$ is a proper embedding,  $\mb$ and $X(\mb)$ are identified via $X$  and $\mb\subset \l^3$ is said to be a properly embedded  $^*$maximal surface.
\end{definition}
With the previous notation, it is not hard to check that $X$ is a conformal $^*$maximal
immersion if and only if $S_{\hat{X}}=\partial(\mb)$ and
$\partial(\mb)$ consists of integral curves of $\frac{1}{\langle
{\cal N},w\rangle^2}\nabla \langle\hat{X},w\rangle,$ where $w$ is
any timelike vector and $\nabla$ is the gradient computed with
respect to the intrinsic metric\footnote{Despite the degeneracy
of $ds^2$ at $S_{\hat{X}},$ the weighted gradient
$\frac{1}{\langle {\cal N},w\rangle^2} \nabla
\langle\hat{X},w\rangle$ extends analytically to this set.}.
To see this, assume that up to a Lorentzian isometry $w=(0,0,1),$ and
write $\phi_3=-i f(z) dz.$ Then, suppose that
$S_{\hat{X}}=\partial(\mb)$ and take a conformal disc $(U,z=u+i
v)$ in $\hat{\mb}$ centered at $p \in \partial(\mb)$ and
satisfying $J(U)=U,$ $z \circ J=\overline{z}.$ A standard
computation gives $\lambda^2\nabla (t \circ
\hat{X})=\mbox{Re}(f)\frac{\partial \hat{X}}{\partial u}
-\mbox{Im}(f)\frac{\partial \hat{X}}{\partial v},$ where
$\lambda=\|\frac{\partial \hat{X}}{\partial u}\|=\|\frac{\partial
\hat{X}}{\partial v}\|=\frac{1}{2} \left(1/{|g|}-|g|\right) |f|.$
Therefore, $\partial(\mb) \cap U$ is an integral curve of
$\lambda^2\nabla (t \circ \hat{X})$ if and only if
$\mbox{Im}(f)=0$ (that is to say, $X$ is   $^*$maximal).
Taking into account that  $(t \circ {\cal N})^2 \lambda^2$ is
well defined and positive on $U,$  we are done.

We will need the following basic lemma.

\begin{lemma} \label{lem:prime}
Let $X:\mb \to \l^3$ be a conformal proper  $^*$maximal immersion with Weierstrass data $(g,\phi_3).$
Then, $dg$ and $\phi_3$ never vanish along $\partial(\mb)$ and $\pi \circ X:\mb \to \Pi_0$ is a local embedding.
\end{lemma}
\begin{proof} From (\ref{eq:simetria}), $g$ and $\phi_3$ extend by Schwarz reflection to the double surface $\hat{\mb}.$
Since $\partial(\mb)=|g|^{-1}(1)$ and $\partial(\mb)$ consists of a family of pairwise disjoint proper regular analytical curves in $\hat{\mb},$ the harmonic function $\log(|g|)$ has no singular points on $\partial(\mb)$ and  $dg(p) \neq 0$ for all $p \in \partial(\mb).$

On the other hand, $S_X$ consists of regular lightlike singularities, hence $dX \neq 0$ on $\partial(\mb)$ and  equation (\ref{eq:wei}) gives  $\phi_3(p) \neq 0$ for all $p \in \partial(\mb).$

Let us show that $\pi \circ X$ is a local embedding. Since
$X|_{\mb-\partial(\mb)}$ is spacelike then $(\pi \circ
X)|_{\mb-\partial(\mb)}$ is a local diffeomorphism, hence  we have to deal only with boundary points.  Fix $p \in
\partial(\mb),$ and up to a Lorentzian isometry, suppose $g(p)=1$
and $\hat{X}(p)=O.$ Then take a conformal disc $(D,z)$ in
$\hat{\mb}$ satisfying $z(D)=\d,$ $z(p)=0,$  $J(z)=\overline{z},$
$z(D \cap \mb)=\d^+:=\d \cap \{z \in \c \;:\; \mbox{Im}(z) \geq
0\}$   and $\phi_2(z)=dz.$ From equation (\ref{eq:wei}) and the facts $dg(p),$ $\phi_3(p) \neq 0,$ we get that
$\phi_1(z)= z h(z)\,{dz},$ where $h:D \to \c$ is holomorphic,
$h(p) \neq 0$ and $h\circ J=\overline{h}.$ In the sequel we
identify $D \equiv \d$ and call $D_\epsilon=\{z \in \c \;:\;
|z|<\epsilon\}.$ By the Domain Invariance Theorem, it suffices
to show that $(\pi \circ X)|_{D_\epsilon}$ is injective provided
that $\epsilon>0$ is small enough. Reason by contradiction and
take sequences $\{z_n\}_{n \in \n},$ $\{w_n\}_{n \in \n},$  in
$\d^+$ converging to $0$ satisfying $z_n \neq w_n,$
$\mbox{Re}(z_n)=\mbox{Re}(w_n)$ and $\mbox{Re}(\int_{z_n}^{w_n} z
h(z))=0.$ Therefore we can find $\xi_n$ in the vertical segment
$]z_n,w_n[$ such that $\mbox{Im}(\xi_n h(\xi_n))=0,$  $n \in \n,$ contradicting that
$\{z \in D_\epsilon\;:\; \mbox{Im} (z
h(z))=0\} \subset \r$ provided that $\epsilon$ is small enough. \end{proof}

The main global result about maximal surfaces was proved by Calabi \cite{calabi} (see also \cite{c-y} for further generalizations). It asserts the following:
\begin{theorem}[Calabi] \label{th:calabi}
Let  $X:\mb \to \l^3$ be a complete maximal immersion, where $\partial(\mb)=\emptyset.$ Then
$X(\mb)$ is a spacelike plane. The same  result holds if we
replace complete for proper.
\end{theorem}

\section{Basic examples} \label{sec:ejemplos}

The family of properly embedded $^*$maximal surfaces is very vast. We are goint to present only the most basic ones, already described by O. Kobayashi in \cite{kobayashi1}. 

Let $\hat{X}:\c \to \l^3$ be the conformal maximal immersion with regular lightlike singularities  associated to the Weierstrass data $g(z)=e^{i z},$ $\phi_3(z)=-idz.$ If $z=u+iv,$ equation (\ref{eq:wei}) gives $$X(u,v)=(\cosh (v) \cos (u) , \cosh (v) \sin(u), u).$$ Since $\hat{X}(\overline{z})=\hat{X}(z)$ and $X=\hat{X}|_{\overline{\u}}$ is a proper embedding, then  ${\cal H}:=X(\overline{\u})$ is a properly embedded   $^*$maximal surface  which has been named as the {\em Lorentzian helicoid}, see Figure \ref{fig:helicoides},$(b)$. 
The conjugate immersion of $\hat{X}$ is the universal converging of the  {\em Lorentzian catenoid}. The Lorentzian catenoid has Weierstrass data $\c-\{0\},$ $g(z)=z,$ $\phi_3(z)=\frac{i dz}{z},$  and it  is given by $$Y(m,s)=(\frac{1-m^2}{2m}\sin (s),\frac{m^2-1}{2m} \cos (s),\log (m)),$$ where $z=m e^{is}.$  In this case, $S_Y=\{|z|=1|\}$ consists of regular lightlike singularities, $Y(S_Y)$ is a single point, $Y(1/\overline{z})=-Y(z)$  and $C:=Y(\overline{\d}-\{0\})$ is an entire graph over $\r^2.$  Elementary characterizations of the Lorentzian catenoid can be found in \cite{kobayashi}, \cite{ecker}  and \cite{f-l-s}.\\

Consider now the data $\mb=\c,$ $g(z)=(z-i)/(z+i)$ and $\phi_3(z)=i(z^2+1)dz.$ Writing $z=m e^{i s},$ the corresponding maximal immersion $\hat{X}:\c \to \l^3$ is given by:  $$\hat{X}((m,s))=\left(-m^2 \cos (2 s), \frac{1}{3} (3m \cos(s) - m^3 \cos(3s)), -\frac{1}{3}m (3 \cos(s)+ m^2 \cos(3s))\right).$$
Since $\hat{X}({\overline{z}})=\hat{X}(z)$ and $X=\hat{X}|_{\overline{\u}}$ is a proper embedding, then  $E_1:=X(\overline{\u})$ is a properly embedded   $^*$maximal surface, that we call the {\em first Enneper's maximal surface}, see Figure \ref{fig:helicoides},$(c)$. $E_1$ contains a half line parallel to the $x_1$-axis and is invariant under the reflection about this line. 
The conjugate surface $E_1^*$ is called the {\em second Enneper's maximal surface}. Its Weierstrass data are $\mb=\c,\;g(z)=(z-i)/(z+i),\;\phi_3(z)=-(z^2+1)dz,$ and putting $z=m e^{i s},$ the immersion $X:\overline{\u} \to \l^3$ is given by $$X(m,s)=\left(m^2 \sin (2s),\frac{1}{3} (-3 m \sin(s) + m^3 \sin(3s)),\frac{1}{3}(3 m \sin(s) + m^3 \sin(3s))\right).$$  In this case  $S_X$ is the real axis, $X(S_X)$ is the origin and  $X$ is not proper. Indeed, $E_2=\overline{X(\overline{\u})}$ is an entire graph over $\r^2$ and $E_2-X(\overline{\u})$ is the open lightlike half line $ x_1=x_2-t=0,$ $t>0.$

\begin{figure}[htpb]
\begin{center}\includegraphics[width=.7\textwidth]{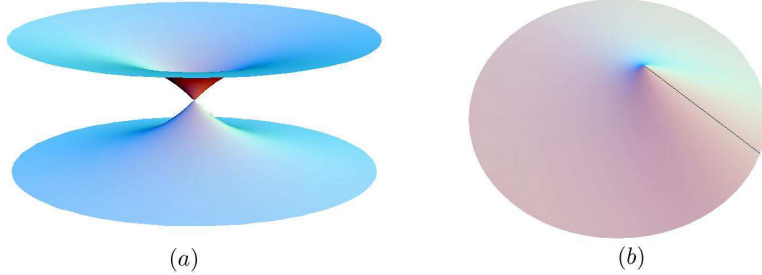}
\caption{$(a$) The Lorentzian catenoid; $(b)$ Enneper's graph
$E_2$ and the line $l$}\label{fig:conjugadas}\end{center}
\end{figure}

\begin{remark} \label{re:cono}
The maximal graphs $C$ and $E_1^*$ satisfy the implicit equations $x_1^2+x_2^2-\sinh^2 {(t)}=0$ and $3 (x_2-t)^2 - \frac{1}{4} (x_2-t)^4 + 3 x_1^2 + 6 (x_2-t) t = 0,$ respectively. Therefore,   any blow-up with center $O$ of these surfaces  converges in the ${\cal C}^0$-topology to half of the lightcone.
\end{remark}

Ecker \cite{ecker} proved that the Lorentzian catenoid is the unique entire maximal graph with one singular point. A similar result for $E_2$ can be found in Section \ref{sec:grafos} (Proposition \ref{pro:unienn}).

\section{Parabolicity of maximal surfaces in $\l^3$} \label{sec:parabo}

This section is devoted to proving some parabolicity criteria for properly immersed maximal surfaces in $\l^3.$ The required background can be found in
\cite{grigor},\cite{ahlfors},\cite{c-k-m-r}, \cite{lop-per} and \cite{joaquin}.

A non compact Riemann surface $\mb$ with non empty boundary is said to be {\em parabolic} if the only bounded harmonic function $f$ vanishing on $\partial (\mb)$ is the constant function $f=0,$ or equivalently, if there exists a proper positive  superharmonic function on $\mb.$ Otherwise, $\mb$ is said to be hyperbolic. If $\partial (\mb)=\emptyset,$ parabolicity means that positive superharmonic functions are constant.

For instance, $\overline{\u}$ is  parabolic, whereas $\d \cap \overline{\u}$ is hyperbolic.

Let $g:\overline{\u} \to \c$ be  continuous on $\overline{\u}$ and  harmonic on $\u.$  A divergent curve $\alpha \subset \overline{\u}$ is defined to be an {\em asymptotic curve} of $g$ if the limit $a:=\lim_{z \in \alpha \to \infty} g(z)\in \overline{\c}$ exists. In this case, $a$ is said to be an {\em  asymptotic value} of $g.$ The following theorem summarizes some well known classical results (see \cite{schiff}).

\begin{theorem} \label{th:sectorial}
Set $g: \overline{\u} \to \c$  continuous, holomorphic on $\u$ and omitting two finite complex values. Then:
\begin{enumerate}[(I)]
\item  $g$ has at most one asymptotic value, and in this case $g|_{\u}$ has angular limits at $\infty.$
\item If the boundary segments $[0,+\infty[$ and $]-\infty,0]$  are asymptotic curves of $g,$ then the limit $\lim_{z \to \infty} g(z)$ exists.
\end{enumerate}
\end{theorem}

Given a Riemann surface $\mb$  with non empty boundary  and $p \in \mb-\partial(\mb),$  we denote by $\mu _p$ the {\em harmonic measure} respect to the  $p.$
It is well known that $\mb$ is parabolic if and only if there exists $p_0 \in \mb-\partial(\mb)$ such that $\mu_{p_0}$ is full, i.e., $\mu_{p_0}(\partial (\mb))=1.$ In this case $\mu_p$ is full for any $p \in \mb-\partial(\mb),$ and bounded harmonic (superharmonic) functions $u$ on $\mb$ satisfy the mean  property
$$u(p)=(\geq)\int _{x\in \partial \mb} u(x)\, d \mu _p, \quad \mbox{ for any }p\in \mb.$$

Regions of parabolic Riemann surfaces  are parabolic, and if a Riemann surface is the union of two  parabolic regions with compact intersection then it is also parabolic.  The proof of the following theorem has been inspired by some ideas in \cite{c-k-m-r}.

\begin{theorem} \label{th:parabo}
Let $X:\mb \to \l^3$ be a conformal proper maximal immersion with singularities, where $\partial (\mb) \neq \emptyset,$ and suppose that there exists $\varepsilon>0$ and a compact subset $C\subset\mb$ such that
$$\langle X,X \rangle\geq \varepsilon \quad \mbox{on} \; \mb- C,$$

Then $\mb$ is parabolic.
\end{theorem}

\begin{proof}

Since parabolicity is not affected by adding compact subsets, we can suppose that
$\|X(p)\|^2\geq \varepsilon$ on $\mb.$

For any $n\in\n$  let $\mb_n:=\{p\in\mb\;:\; \langle X,X \rangle (p)\leq n\}.$
Let us see that $\mb_n$ is parabolic. Indeed, since $t \circ X$ is a proper positive harmonic function on $\mb_n^+:=(t \circ X)^{-1}([0,+\infty[)\cap \mb_n,$ $\mb_n^+$ is parabolic, and likewise $\mb_n^-:=(t \circ X)^{-1}(]-\infty,0])\cap \mb_n$ is  parabolic. As $\mb_n^+ \cap \mb_n^-$ is compact and $\mb_n=\mb_n^+ \cup \mb_n^-,$ then we are done.

Now define $h:\mb \to \r,$ $h(p)=\mbox{log}\langle X,X \rangle (p).$ Since $X$ is  maximal, a direct computation gives that $\Delta h=-4 \frac{<X,{\cal N}>^2}{<X,X>^2}\leq 0,$ where $\Delta$ is the intrinsic Laplacian and ${\cal N}$ is the Lorentzian Gauss map of $X.$ Therefore $h$  is superharmonic.

Without loss of generality, suppose there exists $p\in\mb$ with $h(p)> 0$ (otherwise $\mb=\mb_1$ and we have finished).
Up to rescaling assume that $h(p)=1.$
Since $h$ is a bounded superharmonic function on the parabolic surface $\mb_n,$  we have
$$1=h(p) \geq \int_{\partial(\mb_n) } h\, d\mu_p(n)=\int_{\partial(\mb)\cap\mb_n} h\, d\mu_p(n) + \int_{h^{-1}(\log (n))} \log (n)\, d\mu_p(n) $$
where $\mu_p^n$ denotes the harmonic measure in $\mb_n$ respect to $p.$ Since $0 \leq \int_{\partial(\mb)\cap\mb_n} d\mu_p^n \leq 1,$
$$1 \geq \log (\varepsilon)\int_{\partial(\mb)\cap\mb_n} d\mu_p^n +  \log (n) \int_{h^{-1}(\log (n))} d\mu_p^n \geq 
-|\log (\varepsilon)| +  \log (n) \int_{h^{-1}(\log (n))} d\mu_p^n.$$
Dividing by $\log (n)$ and taking the limit as $n$ goes to infinity, we get
$\lim_{n\to\infty} \int_{\{q\in\mb\;:\;h(q)=n\}} d\mu_p^n \leq 0,$ and so $\lim_{n\to\infty} \int_{h^{-1}(\log (n))} d\mu_p^n = 0.$

On the other hand, the parabolicity of  $\mb_n$ gives
$$1= \int_{\partial(\mb_n)} d\mu_p^n=\int_{\partial(\mb)\cap\mb_n} d\mu_p^n+ \int_{h^{-1}(\log (n))} d\mu_p^n$$
Taking the limit as $n \to \infty$ we get that $1=\int_{\partial(\mb)} d\mu_p,$
where $\mu_p$ is the harmonic measure in $\mb$ with respect to $p,$ concluding the proof. \end{proof}

\begin{corollary} \label{co:halfspace}
Let $X:\mb \to \l^3$ be a conformal proper maximal immersion with singularities, where $\partial (\mb) \neq \emptyset,$ and suppose that $X(\mb)$ lies in a spacelike half plane. Then $\mb$ is parabolic.
\end{corollary}
\begin{proof} Up to scaling and Lorentzian isometry, suppose $X(\mb) \subset \{t \geq 0\}.$

From Theorem \ref{th:parabo}, $\mb_n:=\{p\in\mb\;:\; (t\circ X)(p)\leq n\}$ is parabolic, $n \in \n.$ Defining now $h=t \circ X$ and reasoning as in the preceding proof we obtain the desired conclusion. \end{proof}
\section{Some results on maximal graphs} \label{sec:grafos}

The space of continuous functions $u$ on a domain $\Omega \subset \r^2$ with weak gradient satisfying $\|\nabla u\|_0 \leq 1$ will be denoted by  ${\cal C}^0_1(\Omega).$ We endow   ${\cal C}^0_1(\Omega)$ with the ${\cal C}^0$-topology of the uniform convergence on compact subsets of $\Omega.$ Likewise, and for any $k \in \n \cup \{\infty\},$  ${\cal C}_1^k(\Omega)$ will denote  the space of functions with continuous partial derivatives of order $< k+1,$ endowed with the ${\cal C}^k$-topology of the uniform convergence of $u$ and its partial derivatives of order $< k+1$ on compact subsets of $\Omega.$

A sequence of PS graphs $\{(x,u_n(x)) \;:\; x \in \Omega\}$, $n \in \n,$  is said to be convergent in the ${\cal C}^k$-topology to $\{(x,u(x) \;:\; x \in \Omega\}$  if $\{u_n\}_{n \in \n} \to u$  in the ${\cal C}^k$-topology, $k \in \n \cup \{\infty\}.$

\begin{remark} \label{re:convergencia}
If $\Omega$ is bounded, ${\cal C}^0_1(\Omega)$ lies in the Sobolev space ${\cal W}^{1,2}(\Omega)$ of $L^2$ functions with $L^2$ gradient, and the convergence in ${\cal C}^0_1(\Omega)$ implies the one in ${\cal W}^{1,2}(\Omega)$. By  Ascoli-Arzela theorem, any  sequence in ${\cal C}_1^0(\Omega)$ bounded at $x_0 \in \Omega$  contains a subsequence converging in both ${\cal C}_1^0(\Omega)$  and  ${\cal W}^{1,2}(\Omega).$
\end{remark}

Let  $u \in {\cal C}_1^\infty(\Omega).$ The associated graph $G=\{(x,u(x)) \;:\; x \in \Omega\}$ defines a maximal surface if and only if:
\begin{equation} \label{eq:max}
\|\nabla u \|_0<1\; \mbox{and}\; \mbox{div} \left({\nabla u}/{\sqrt{1-\|\nabla u\|_0^2}} \right)=0.
\end{equation}

The conjugate function $u^*$ is characterized by the identity
\begin{equation} \label{eq:conjugada}
\sqrt{1-\|\nabla u\|^2_0}\, du^*=\frac{\partial u}{\partial x_2} dx_1-\frac{\partial u}{\partial x_1} dx_2,
\end{equation}
 besides the initial condition. It is well defined if and only if $\frac{\partial u}{\partial x_2} dx_1-\frac{\partial u}{\partial x_1} dx_2$ is an exact 1-form (for instance, if $\Omega$ is simply connected), and satisfies the minimal surface equation  $$\mbox{div} \large( {\nabla u^*}/{\sqrt{1+\|\nabla u^*\|_0^2}} \large)=0.$$ Thus, $G^*=\{(x,u^*(x)) \;:\; x \in \Omega\}$ is a minimal surface in $\r^3.$
In terms of the Weierstrass representation, the conformal maximal and minimal immersions associated to $G$ and $G^*$ are given by
\begin{equation} \label{eq:weiconju}
X=\mbox{Real} \int (\phi_1,\phi_2,i \phi_3)\; \mbox{and}\; X^*=\mbox{Real} \int (\phi_1,\phi_2,\phi_3),
\end{equation}
 respectively.

The following theorem is the Lorentzian version of classical Plateau's problem.

\begin{theorem}[\cite{bar-sim}] \label{th:plateau}
Let $\gamma \subset \l^3$ be a Jordan curve bounding a PS embedded surface. Then there exists a PS  area maximizing disc ${\cal S}$ in $\l^3$ bounded by $\gamma.$ Furthermore, ${\cal S}$ is smooth (hence a maximal surface) except possibly on  piecewise linear lightlike arcs connecting points of $\gamma.$
\end{theorem}
This result  applies to curves $\gamma$ whose projection $\pi(\gamma)$ is a Jordan curve and $|t(p)-t(q)|\leq d_{\Omega}(\pi(p),\pi(q))$ for any $p,$ $q \in \gamma,$ where $\Omega$ is the domain bounded by $\pi(\gamma)$ and $d_\Omega$ is the inner distance in $\Omega$ (see \cite{kly-mik}). Furthermore,  ${\cal S}$ is smooth provided that $|t(p)-t(q)|< d_{\Omega}(\pi(p),\pi(q))$ for any $p,$ $q \in \gamma.$\\

Let $\Omega\subset \r^2$ be a bounded domain and consider a sequence $\{u_n\}_{n \in \n} \subset {\cal C}_1^\infty (\Omega)$  of functions satisfying equation (\ref{eq:max}). Assume that  $\{u_n\}_{n \in \n} \to u $ in the ${\cal C}^0$-topology, where $u \in {\cal C}_1^0(\Omega)$ (see Remark
\ref{re:convergencia}). 

Given $(x,y) \in \partial(\Omega)^2,$ the segment $]x,y[$ is said to be singular if $]x,y[ \subset {\Omega}$ and $|u(x)-u(y)|=\|x-y\|_0.$\footnote{Recall that $u$ is locally Lipschitzian, hence it extends continously to $\overline{\Omega}.$} We write
$A=\{(x,y) \in \partial(\Omega)^2 \,:\, ]x,y[ \;\mbox{is singular}\}$ and  set ${\cal A}:=\cup_{(x,y) \in A} ]x,y[ \subset \Omega$ the (closed) {\em singular set} of $u.$ Next theorem summarizes some known results mainly proved by Bartnik and Simon in \cite{bar-sim}:

\begin{theorem} \label{th:converge} The function $u$ defines an area maximizing PS graph. Moreover,
\begin{enumerate}[(A)]
\item $u|_{\Omega-{\cal A}}$ satisfies (\ref{eq:max}) and $\{u_n|_{\Omega-{\cal A}}\}_{n \in \n} \to u|_{\Omega-{\cal A}}$ in the ${\cal C}^\infty$-topology,
\item $x_0 \in {\cal A}$ if and only if there exists $\{x_n\}_{n \in \n} \to x_0$ such that $\{\|\nabla u_n (x_n)\|_0\}_{n \in \n} \to 1,$ and
\item for any $(x,y) \in A,$ $\{(z,u(z)) \;:\; z \in ]x,y[\}$ is a lightlike segment.
\end{enumerate}
\end{theorem}

\begin{remark} \label{re:strip}
From Lemma \ref{lem:basico}, $(b),$  different points of $A$ determine disjoint singular segments of $u,$ hence ${\cal A}$ is a closed subset of $\Omega$ foliated by singular segments of $u.$ If $\Omega$ is unbounded and $\{\Omega_n\}_{n \in \n}$ is an exhaustion  of $\Omega$ by bounded domains, we set  ${\cal A}:=\cup_{n \in \n} {\cal A}_n,$ where ${\cal A}_n$ is the singular set of $u|_{\Omega_n},$ $n \in \n.$ Since ${\cal A}_n \subset {\cal A}_{n+1}$ for any  $n,$  ${\cal A}$ is foliated too by inextensible segments in $\Omega,$ but in this case some of them could have infinite length.

Since  maximal surfaces are locally graphical, the notions of singular set and singular segment can be straightforwardly extended to the limit of a sequence of maximal surfaces in $\r^3_1.$
\end{remark}

\begin{remark} The singular set ${\cal A}$ of the area maximizing PS graph associated to $u$ and the singular set $S_X$ of a conformal maximal immersion $X:\mb \to \l^3$ have different nature and should not be confused.
The set $S_X$ lies in the conformal support of $X$ and has vanishing measure, whereas ${\cal A}$ is contained in the domain of $u$ possibly with non zero Lebesgue measure. 

For instance, the blow up the catenoid is the lightcone. In this case ${\cal A}=\r^2$ and $S_X$ has no sense because the limit is not a conformal maximal surface with singularities.
\end{remark}

Given $(x,y) \in A,$ we call $\Sigma_{(x,y)}$ as the unique lightlike plane containing $\{(z,u(z)) \;:\; z \in ]x,y[\},$ and
set $\sigma_{(x,y)}\in \r^2 $ the  unitary vector for which $\Sigma_{(x,y)}=\{(z,\langle z,\sigma_{(x,y)}\rangle_0) \;:\; z \in \r^2\}.$ Since $u \in {\cal W}^{1,2}(\Omega),$ $\nabla u$ is well defined {\em almost everywhere} on $\Omega$ (that is to say, on a subset $\Omega_0\subset \Omega$ having the same Lebesgue measure as $\Omega$). Furthermore, item $(A)$ in  Theorem \ref{th:converge} implies that $\Omega-{\cal A} \subset \Omega_0,$ whereas item $(C)$ and the PS property give that $\nabla u(z)=\sigma_{(x,y)}$ provided that $(x,y) \in A$ and $z \in ]x,y[\cap \Omega_0.$ 

Therefore, it is natural to define ${\cal D} u:\Omega \to \r^2$ by  ${\cal D} u|_{\Omega_0}=\nabla u$ and ${\cal D} u(z)=\sigma_{(x,y)}$ if $z \in ]x,y[,$ $(x,y) \in A.$ Obviously ${\cal A}=\|{\cal D} u\|_0^{-1}(1).$

\begin{proposition} \label{pro:converge}
$u\in {\cal C}_1^1(\Omega),$ ${\cal D} u=\nabla u$ and $\{ u_n\}_{n \in \n} \to u$ in the ${\cal C}^1$-topology.
\end{proposition}
\begin{proof}
Take a sequence  $\{x_n\}_{n \in \n} \subset \Omega$ converging to $x_0 \in \Omega$ and such that the limit  $\sigma:=\lim_{n \to \infty} \nabla u_n(x_n)$ exists.

\begin{quote}
{\bf Claim 1:} {\em If $\|\sigma\|_0<1$  then $x_0 \in \Omega-{\cal A}$ and $\sigma={\cal D} u(x_0)=\nabla u (x_0).$}
\end{quote}
\begin{proof} Let $D_0 \subset \Omega$ be a  closed disc of positive radius centered at $x_0.$ Take $\epsilon \in ]\|\sigma\|_0,1[$ and without loss of generality suppose $\|\nabla u_n (x_n)\|_0<\epsilon,$ for all $n \in \n.$  Label $u_n^*$ as the conjugate function of $u_n|_{D_0}$ satisfying $u_n^*(x_n)=0$ (well defined because $D_0$ is simply connected,  see equation (\ref{eq:conjugada})), and denote by  $S_n:=\{(x,u_n^*(x)) \;:\; x \in D_0\}$ the associated minimal graph, $n \in \n.$ Standard curvature estimates for minimal graphs give that  $|K_n|\leq C_1$ on $D_0$ for any $n \in \n,$ where $K_n$ is the Gaussian curvature of $S_n$ and $C_1$ is a constant depending only on $\mbox{d}(D_0,\partial(\Omega))>0.$  From our hypothesis, $\|\nabla u_n^*\|_0(x_n)<\frac{\epsilon}{\sqrt{1-\epsilon^2}},$ and taking $\delta>\frac{\epsilon}{\sqrt{1-\epsilon^2}},$ the Uniform Graph Lemma for minimal surfaces \cite{joa-ros} implies the existence of a smaller disc $D \subset D_0$ centered at $x_0$ such that $\|\nabla u_n^*\|_0<\delta$  on $D,$ for any $n \in \n.$ Thus,   $\|\nabla u_n\|_0<\frac{\delta}{\sqrt{1+\delta^2}}<1$ on $D$ for all $k \in \n.$ Barnik-Simon results in \cite{bar-sim} give that $\{{u_{n_k}}|_D\}_{k \in \n}\to u|_D$  in the ${\cal C}^\infty-$topology  and $u|_D$ satisfies the maximal surface equation, (that is to say, $D \subset \Omega-{\cal A}$). In particular, $\sigma =\nabla u(x_0)={\cal D} u(x_0)$ and we are done. \end{proof}

\begin{quote}
{\bf Claim 2:} {\em If  $\|\sigma\|_0=1$  then $x_0 \in {\cal A}$ and $\sigma={\cal D} u(x_0).$}
\end{quote}

\begin{proof} It is clear that $x_0 \in {\cal A}$ (see Theorem \ref{th:converge}, $(B)$).
Consider $\{\mu_n\}_{n \in \n} \to 0,$ $\mu_n>0,$ and define $\Omega_n:=\frac{1}{\mu_n} \cdot (\Omega-x_n)$ and $v_n:\Omega_n \to \r, \;v_n(y):= \frac{1}{\mu_n} \left(u_n(\mu_n y+x_n)-u_n(x_n)\right),$ $n \in \n.$

Let us show that up to subsequences,  $\{v_n\}_{n \in \n} \to v$   in the ${\cal C}^0$-topology, where $v:\r^2\to \r,$ $v(y):= (y,\langle \sigma,y \rangle_0).$ Since $v_n$  lie in ${\cal C}_1^0(\Omega_n)$ and vanish at the origin, Remark \ref{re:convergencia} yields that, up to subsequences,  $\{v_n\}_{n \in \n}\to v_0$  in the ${\cal C}^0$-topology, where $v_0 \in {\cal C}_1^0(\r^2).$ We have to show that $v=v_0.$

Call $G$ as the entire graph defined by $v_0,$ and for any bounded  domain $\Omega' \subset \r^2$ label ${\cal A}_{\Omega'}$ as the singular set  of $v_0|_{\Omega'}.$   If ${\cal A}_{\Omega'}=\emptyset$ for any $\Omega',$ Theorem \ref{th:converge} and Calabi's theorem  would imply that $\{v_n\}_{n \in \n} \to v_0$ in the ${\cal C}^\infty$-topology and $v_0$ is a linear map defining a spacelike plane, contradicting that $\sigma=\lim_{n \to \infty} \nabla v_n(0)$ is a unitary vector.  Therefore ${\cal A}_{\r^2} \neq \emptyset$ and  $G$ contains  a lightlike straight line. From Lemma \ref{lem:basico}, $G$ must be a lightlike plane and so $\r^2$ is foliated by singular straight lines of $v_0.$ As a consequence, Claim 1 implies that $\{\|\nabla v_n\|_0\}_{n \in \n} \to 1$ in the ${\cal C}^0$-topology over $\r^2.$

In the sequel we will assume that $\Omega$ is simply connected (otherwise, replace $\Omega$ for a small enough disc centered at $x_0$).
Let $v_n^*:\Omega_n \to \r$ denote the conjugate function of $v_n$ with initial condition $v_n^*(0)=0,$  and label $S_n$ as its associated minimal graph.
Let $\Pi_n$ denote the tangent plane of $S_n$ at $0,$ i.e., the plane passing through $0$ and orthogonal in the Euclidean sense to the vector $\sqrt{1-\|\nabla v_n \|^2_0} \left(-\nabla v_n^*,1\right)(0)=\left(\nabla u_n ^\bot,\sqrt{1-\|\nabla u_n \|^2_0}\right)(x_n),$ where  $\nabla u_n^\bot (x_n)=(-\frac{\partial u_n}{\partial y},\frac{\partial u_n}{\partial x})(x_n).$ The limit plane $\Sigma=\lim_{n \to \infty} \Pi_n$ is orthogonal to $(\sigma^\bot,0),$ where $\sigma^\bot=(-w_2,w_1)$ provided that $\sigma=(w_1,w_2).$

By standard curvature estimates and the Uniform Graph Lemma for minimal graphs \cite{joa-ros}, we can find a graph $S_n'\subset S_n,$ $n \in \n,$ such that   $\{S_n'\}_{n \in \n}\to \Sigma$ in the ${\cal C}^\infty$-topology as graphs over $\Sigma.$

Let $\hat{\gamma}_n(s): [-L'_n,L_n] \to \Pi_0$  be the arc-length parameterized inextensible arc  in $S_n \cap \Pi_0$ satisfying $\hat{\gamma}_n(0)=0.$ Write $\hat{\gamma}_n=(\gamma_n,0),$  and note that $\{[L'_n,L_n]\}_{n \in \n}\to \r$ and  $\{\gamma_n(s)\}_{n \in \n} \to \gamma_0$  in the ${\cal C}^\infty-$topology over $\r,$ where $ \gamma_0:\r \to \r^2$ is given by $\gamma_0(s)= s\sigma.$

For each $n \in \n,$ call $G_n$ the maximal graph determined by $v_n$ and  set $\alpha_n:[-L'_n,L_n] \to G_n,$ $\alpha_n(s):=(\gamma_n(s),v_n(s)),$ where $v_n(s):=v_n(\gamma_n(s)).$  It is not hard to see that $\{\alpha_n'(s)\}_{n \in \n} \to (\sigma,1)$ in the ${\cal C}^\infty-$topology. Indeed, since $\alpha_n'(s)=(\gamma_n'(s),\langle \gamma_n'(s),\nabla v_n(\gamma_n(s))\rangle_0),$ it suffices to check that $\{\nabla v_n(\gamma_n(s))\}_{n \in \n}\to \sigma.$ Taking into account that $v_n^*(\gamma_n(s))=0$ for any $s,$ we have that $\langle \nabla v_n^*(\gamma_n(s)), \gamma_n'(s) \rangle_0 =0,$ and so, $\nabla v_n (\gamma_n(s))=\lambda_n(s) \gamma_n'(s).$ As $\{\|\nabla v_n\|_0\}_{n \in \n} \to 1$ in the ${\cal C}^0$-topology,  then $|\lambda_n(s)| \to 1$ uniformly on compact subsets of $\r.$ But  $\{\nabla v_n (0)\}_{n \in \n} \to \sigma$ implies $\lambda_n(0)=1,$ hence $\{\lambda_n(s)\}_{n \in \n} \to 1$ and $\{\nabla v_n(\gamma_n(s))\}_{n \in \n} \to \sigma.$ 

As a consequence, $\{\alpha_n\}_{n \in \n}$ converges to the lightlike straight line $\alpha_0:\r \to  G,$ $\alpha_0(s)=(s \sigma,s),$ hence $G$ is the lightlike plane containing $\alpha_0$ and $v=v_0,$ proving our assertion.

To finish the claim,  take a closed disc $D \subset \Omega$ of positive radius centered at $x_0,$ and without loss of generality, suppose  $\{x_n, \; n \in \n\} \subset D.$
Label $\mu_n:=\max\{|u_n(x)-u(x)|,\; x \in D\},$  and define $v_n,$ $w_n:\Omega_n \to \r$ by $v_n(y)= \frac{1}{\mu_n} (u_n(\mu_n y+x_n)-u_n(x_n)),$ $w_n(y)= \frac{1}{\mu_n} (u(\mu_n y+x_n)-u(x_n)),$ $n \in \n.$ We know that  $\{v_n\}_{n \in \n}\to v,$  and by Remark \ref{re:convergencia}  $\{w_n\}_{n \in \n}\to w$ in the ${\cal C}^0$-topology, where $w\in {\cal C}_1^0(\r^2).$

If $]x,y[\subset {\cal A}$ be the inextensible singular segment of $u$ containing $x_0,$ then the PS graph $G':=\{(y,w(y))\;:\; y \in \r^2\}$ contains the  straight line passing through $O$ and parallel to the lightlike vector $(y-x,u(y)-u(x)).$ From Lemma \ref{lem:basico}, $G'$ is the lightlike plane parallel to this vector, and so $w(y)=(y,\langle y,{\cal D} u(x_0) \rangle_0),$ for any $y \in \r^2.$

Setting $D_n:= \frac{1}{\mu_n} (D-x_n),$ the graphs $G_n:=\{(y,v_n(y)) \,:\, y \in D_n \}$ and $G'_n:=\{(y,w_n(y)) \,:\, y \in D_n \}$ satisfy  $\mbox{d}_H(G_n,G'_n )\leq 2,$ $n \in \n,$ hence $\mbox{d}_H(G,G' )\leq 2.$ This implies that $G$ and $G'$ must be parallel and so $\sigma={\cal D} u(x_0),$ which proves the claim. \end{proof}
Claims 1 and 2 imply that $\{\|\nabla u_n-{\cal D} u\|_0\}_{n \in \n} \to 0$ in the ${\cal C}^0$-topology over $\r^2.$ Let us show that ${\cal D} u$ is continuous on $\Omega.$ From Theorem \ref{th:converge}, ${\cal D} u$ is  continuous on  $\Omega-{\cal A},$ and Lemma \ref{lem:basico}, $(ii)$ and  Remark \ref{re:strip}) show that $\sigma_{(x,y)}$ depends continuously on $(x,y)\in A,$ hence ${\cal D} u$ is continuous on ${\cal A}$ too.
Therefore, it suffices to prove that $\lim_{k \to \infty} {\cal D} u (y_k)={\cal D} u(x_0),$ provided that $\{y_k\}_{k \in \n} \subset \Omega-{\cal A}$  and $\lim_{k \to +\infty} y_k=x_0 \in \partial({\cal A}).$ To see this, use Theorem \ref{th:converge}, (A) to find a divergent sequence $\{n_k\}_{k \in \n}$ in $\n$ such that $\|\nabla u_{n_k} (y_k)-{\cal D} u (y_k)\|_0<1/k,$ for  any $k \in \n.$ From Claims 1 and 2,  $\lim_{k  \to +\infty} \nabla u_{n_k} (y_k)= {\cal D} u(x_0),$ and so $\lim_{k \to \infty} {\cal D} u (y_k)={\cal D} u(x_0).$

Finally, fix $x_0 \in \Omega$ and define $du=\langle {\cal D} u, (dx,dy)\rangle_0.$ For any $x \in \Omega$ and any smooth curve $\alpha\subset \Omega$ connecting $x_0$ and $x$ one has $u(x)=\lim_{n \to \infty} u_n(x)=\lim_{n \to \infty} (u_n(x_0)+ \int_{\alpha} d u_n)= u(x_0)+ \int_{\alpha} d u.$ Since $du$ is continuous then $u \in {\cal C}_1^1(\Omega),$ concluding the proof. \end{proof}

\subsection{Asymptotic behavior of maximal multigraphs of finite angle}

Set $G=\{(z(x),u(x)) \;:\; x \in W\}$ a PS multigraph over a wedge $W \subset {\cal R}$ of finite angle. Let $\mbox{d}_W$ be the intrinsic distance in $W$ induced by $|dz|^2,$ and fix $x_0 \in W.$  Since $\|\nabla u\|_0\leq 1$ and $\lim_{x \in W \to \infty} \frac{\mbox{d}_W(x_0,x)}{\|z(x)\|_0}=1,$ then $\limsup_{x \in W \to \infty} \frac{|u(x)|}{\|z(x)\|_0}\leq 1$ and  $\tau^+(r):=1-\min\{\frac{u(x)}{r}\,:\; \|z(x)\|_0=r\}\in [0,2]$ for any $r >0.$

Define $\tau^+(G):=\limsup_{r \to +\infty} \tau^+(r),$ $\tau^-(G):=\tau^+(-G)$ and $\tau(G)=\min\{\tau^+ (G),\tau^-(G)\}.$

Likewise, $\tau_0^+(G):=\liminf_{r \to +\infty} \tau^+(r),$ $\tau_0^-(G):=\tau_0^+(-G),$ and $\tau_0(G)=\min \{\tau_0^+ (G),\tau_0^-(G)\}.$

These numbers give different measures of the asymptotic closeness between $G$ and the light cone. Obviously, $\tau^+(G)\geq \tau_0^+(G)$ and  $\tau^-(G)\geq \tau_0^-(G).$

For $\theta\in]0,+\infty[,$ call
$$\Xi(\theta):=\inf\{\tau(G) \;:\; G\;\mbox{is a maximal multigraph over a wedge of angle} \; \theta\}$$
$$\Xi_0(\theta):=\inf\{\tau_0(G) \;:\; G\;\mbox{is a maximal multigraph over a wedge of angle} \; \theta\}$$ and notice that $\Xi(\theta) \geq \Xi_0(\theta).$

The monotonicity formulae $\tau(G') \leq \tau(G)$ and  $\tau_0(G') \leq \tau_0(G),$ provided that $G' \subset G,$ hold straightforwardly. As a consequence,  $\Xi(\theta') \leq \Xi(\theta)$ and $\Xi_0(\theta') \leq \Xi_0(\theta)$ provided that $\theta' \leq \theta.$

\begin{lemma} \label{lem:universal}
$\Xi(\theta)>0$ for any $\theta \in]0,+\infty[.$
\end{lemma}
\begin{proof} Since $\tau^+(G)=\tau^-(-G),$ we have $$\Xi(\theta)=\inf\{\tau^+(G) \;:\; G\;\mbox{is a maximal multigraph over a wedge of angle} \; \theta\}\geq 0.$$

On the other hand, any multigraph of angle $\theta$ contains, up to a translation, a graph over the wedge $W_{\theta'}$ for any $\theta'<\min\{\frac{\theta}{2},\pi\}$ (see Section \ref{sec:prelim}). 
By the above  monotonicity argument, if suffices to prove that  $$\inf\{\tau^+(G) \;:\; G\;\mbox{is a graph over} \; W_\theta\}>0$$ for any $\theta \in ]0,\pi[.$ Reason by contradiction, and assume that there exists $\theta \in ]0,\pi[$ and sequence of maximal graphs $\{G_n\}_{n \in \n}$ over $W_{\theta}$ such that
$\lim_{n \to \infty}\tau^+ (G_n)=0.$ Write $G_n=\{(x,u_n(x))\;:\; x\in W_{\theta}\},$ and without loss of generality suppose $u_n((0,0))=0,$ $n \in \n.$ From equation (\ref{eq:acausal}) and up to scaling depending on $n,$ we can also assume that
\begin{equation} \label{eq:oju}
u_n(x)/\|x\|_0\in [1-\tau^+(G_n)-\frac{1}{n},1], \; \mbox{for all}\; x \in W_\theta\cap \{\|x\|_0\geq 1\} \; \mbox{and}\;n \in \n.
\end{equation}
Define $v:W_{\theta}\to \r,$ $v(x)=\|x\|_0,$ and let us see that $$\lim_{n \to \infty} \sup \left\{\|\nabla u_n -\nabla v\|_0 \;:\; x \in W_{\theta'}\cap \{\|x\|_0 \geq 1\}\right\}=0,$$ for any $\theta' \in ]0,\theta[.$
Indeed, reason by contradiction and suppose there is a sequence $\{x_n\}_{n \in \n}$ in $W_{\theta'}\cap \{\|x\|_0 \geq 1\}$ such that, and up to subsequences, $\|\nabla u_n(x_n) -\nabla v(x_n)\|_0 \geq \epsilon >0.$ Call  $v_n(y):=\frac{1}{\|x_n\|_0} u_n (\|x_n\|_0 \,y),$ for each $n \in \n.$
The fact that $\{\tau^+(G_n)\}_{n \in \n} \to 0,$ equation (\ref{eq:oju}) and Proposition \ref{pro:converge} imply that $\{v_n\}_{n \in \n} \to v$ in the ${\cal C}^1$-topology on $\mbox{Int}(W_{\theta}),$ contradicting that $\|\nabla v_n(\frac{x_n}{\|x_n\|_0}) -\nabla v(\frac{x_n}{\|x_n\|_0})\|_0 \geq \epsilon >0$ for all $n \in \n$ and proving our assertion.

Let $g_n$ be the holomorphic Gauss map of of $G_n.$ Writing $\nabla u_n=\frac{\partial u_n}{\partial x_1}+i \frac{\partial u_n}{\partial x_2},$ one has that $g_n= \frac{-i}{1+\sqrt{1-\|\nabla u_n\|_0^2}} \nabla u_n.$ Rewriting the above limit in polar coordinates we infer that $\lim_{n \to \infty} \sup \{|g_n(s e^{i \xi})+i e^{i \xi}|\;:\; (\xi,s)\in [-{\theta'},{\theta'}]\times [1,+\infty[\}=0,$ for any  $\theta'\in ]0,\theta[.$

Therefore, fixing $\theta'\in ]0,\theta[$ and $\epsilon \in ]0,\min(\frac{1}{2},\frac{\theta'}{2})[,$ we can find $n_0 \in \n$ large enough in such a way that $|\mbox{Im}(\log(g_{n_0})(s e^{i \xi}))-\xi+\pi/2|<\epsilon$
 and $|g_{n_0}(s e^{i \xi})|>1-\epsilon>\frac{1}{2},$  for every $s \geq 1$ and $\xi \in [-{\theta'},{\theta'}].$ An intermediate value argument gives that $C_\delta:=\{x \in W_{\theta'}\cap \{\|x\|_0\geq 1\}\;:\;\mbox{Im}(\log(g_{n_0})(x))=\delta\}$ is  non compact  and  $C_\delta \cap \partial \left(W_{\theta'} \cap \{\|x\|_0\geq 1\}\right) \subset \{\|x\|_0=1\},$ for any $\delta \in ]-{\theta'}-\pi/2+2\epsilon,{\theta'}-\pi/2-2\epsilon[.$

Choose $\delta$ in such a way that $dg_{n_0}$ never vanishes along $C_\delta,$ and take a divergent regular arc  $\alpha_\delta\subset C_\delta.$ As  $\log |g_{n_0}|$  is  strictly monotone on $\alpha_\delta$, then $\lim_{x \in \alpha_\delta \to \infty}g_{n_0}(x)=r_\delta e^{i\delta},$ $1-\epsilon\leq r_\delta\leq 1.$  In other words, $r_\delta e^{i \delta}$ is an asymptotic value of  $g_{n_0}$ at the unique end of $G_{n_0}.$  This argument works for infinitely many $\delta' \in ]-\theta'-{\pi}/{2}+2\epsilon,\theta'-{\pi}/{2} -2\epsilon[$ different from $\delta,$ and so $g_{n_0}$ has infinitely many asymptotic values. This contradicts the parabolicity of $G_{n_0}$ (see for instance Corollary \ref{co:halfspace}) and Theorem \ref{th:sectorial}, and proves the lemma. \end{proof}

For any $\delta \in ]0,\frac{\pi}{4}[,$ set  $U_\delta$ the region $\{(x_1,x_2,t) \in \l^3 \,:\, \arg{((x_1,t))} \in [\frac{\pi}{2},\pi+\delta]\}.$ 

\begin{lemma} \label{lem:universal1}  Consider a wedge $W \subset {\cal R}$ and a region $\Omega \subset W$ satisfying that  $W_{\frac{\pi}{2}}-D\subset \Omega\subset W,$ where $D$ is an open disc centered at the origin. Let $G=\{(x,u(x)) \,:\, x \in W\}$ be  a maximal multigraph over $W,$ and call $G_0=\{(x,u(x)) \,:\, x \in \Omega\}.$ Assume that $\tau_0^+(G) <1-\tan(\delta)$ and $\partial(G_0)-\pi^{-1}(D) \subset U_\delta$ for some $\delta \in ]0,\frac{\pi}{4}[.$
 
Then $\liminf_{y \in +\infty} \frac{u((y,0))}{y} \geq 1-\tau_0^+(G).$ In particular,  $\tau_0^+(G)=0$ implies that $\{(y,u((y,0)))\,:\, y \in \r\}$ is an upward lightlike ray.
\end{lemma}
\begin{proof} Take $\epsilon \in ]\tan(\delta),1-\tau_0^+(G)[,$ and let $H_\epsilon\subset \l^3$ denote the smallest half space with boundary plane parallel to $\{(x_1,x_2,t)\,:\,  t -\epsilon x_1=0\}$ containing $\partial(G_0) \cup U_\delta.$  

Since $\tau_0^+(G)< 1-\epsilon,$ we can find a divergent sequence $\{r_k\}_{k \in \n}$ in $[1,+\infty[$ such that $G \cap \pi^{-1}(\{x \in \Pi_0 \;:\; \|x\|_0=r_ k\}) \subset H_\epsilon$ for any $k \in \n.$ As $\partial(G_0)\cup  U_\delta \subset H_\epsilon$ then $\partial(G_k)\subset H_\epsilon,$ where $G_k=G_0 \cap \pi^{-1}(\{x \in \Pi_0 \;:\; \|x\|_0\leq r_ k\}.$ The convex hull property gives that $G_k \subset H_\epsilon,$ for any $k>0,$ and therefore $\liminf_{y \in +\infty} \frac{u((y,0))}{y} \geq \epsilon.$ Since $\epsilon$ is an arbitrary real number in $]\tan(\delta),1-\tau_0^+(G)[$ we are done. \end{proof}

\begin{corollary} \label{co:universal}
Let $W$ be a wedge of angle $\geq 4 \pi$ and write $\partial(W)=\alpha_1\cup \alpha_2,$ where $\alpha_j\cong [0,1[,$ $j=1,2.$ Call $\theta_j(W)=\lim_{x \in \alpha_j} \arg(x),$ $j=1,2,$ and $\theta_0(W)=\frac{\theta_1(W)+\theta_2(W)}{2}.$ Consider $G:=\{(x,u(x))\,:\, x \in W\}$ a properly embedded maximal multigraph over $W,$ and assume that there is a region $\Omega \subset W$ and an open disc $D \subset \c\equiv \Pi_0$ such that $\arg^{-1}([\theta_0(W)-\frac{\pi}{2},\theta_0(W)+\frac{\pi}{2}])-z^{-1}(D)\subset \Omega$ and  $\partial(G_0)-\pi^{-1}(D) \subset U_\delta,$ where  $G_0=\{(x,u(x)) \,:\, x \in \Omega\}$ and   $\delta \in ]0,\frac{\pi}{4}[.$

Then there exists a positive constant $\Xi$ not depending on $G$  such that $\tau_0(G) \geq \Xi.$
%
%
%
\end{corollary}
\begin{proof} Since  $\partial(G_0)-\pi^{-1}(D) \subset U_\delta$ then $\limsup_{r \to +\infty} \min\{-\frac{u(x)}{r}\,:\,x\in \Omega, \, \|z(x)\|_0=r\}\leq \tan (\delta),$ hence $\tau_0^-(G) \geq 1-\tan(\delta).$ Thus, it suffices to check that $\tau_0^+(G) \geq \Xi$ for a suitable constant $\Xi>0.$

Reason by contradiction and take a sequence  $\{(W_n,G_n,\Omega_n,D_n)\}_{n \in \n}$ of wedges, embedded multigraphs, domains and discs satisfying the above hypothesis and such that $\lim_{n \to \infty} \tau^+_0(G_n)=0.$ Up to a rotations about the $t$-axis and  reparameterizations  we will suppose that $\theta_0(W_n)=0$  for any $n \in \n.$ Label $l_n$ as the proper arc $\{(y,0,u_n(y))\,:\, y \in [1,+\infty[\}$ in $G_n.$ Using Lemma \ref{lem:universal1} we get that$\liminf_{y \in +\infty} \frac{u_n((y,0))}{y} \geq 1-\tau_0^+(G_n)\geq 0.$ Therefore,   up to removing from $W_n$ a suitable compact subset and choosing a larger $D_n,$ we can suppose that $l_n \subset \{t \geq 0\}$ for any $n \in \n.$ 

Up to scaling, we will also assume that $D_n=\d$ for any $n \in \n,$ and call  $W=W_{2 \pi}-\d.$ Moreover, we replace $G_n$ for $\{(x,u_n(x)) \,:\, x \in W=W_{2 \pi}-\d\}$  keeping the same name for the new multigraph, for any $n \in \n.$

Since $G_n$ is embedded, $\partial(G_n)$ contains an unique proper arc $l_n'$ lying above $l_n$ (and so contained in $\{t \geq 0\}$) and such that $\pi(l_n')=\pi(l_n).$ Write $W':=\{p \in {\cal R}\,:\, \|z(p)\|_0\geq 1\; \mbox{and}\; \arg(p)\in [0,2\pi]\}$ and  $F_n=\{(x,u_n(x)) \;:\; x \in W'\},$ and up to the reflection about the origin, suppose that $\partial(F_n)$ consists of $l_n \cup l_n'$ and a suitable subarc of $\partial(G)\cap \pi^{-1}(\partial(D)).$

Take $c_n \in ]\tau_0^+(F_n),2 \tau_0^+(F_n)[$ and call $V_n:=\{(x,t) \in \l^3 \;:\; |t|>(1-c_n)\|x\|_0\}$ and  $V_n^+=V_n \cap \{ t \geq 0 \},$ $n \in \n.$  For each $\xi \in I:=[\pi/2,3 \pi/2]$  call $H_{n}(\xi)$ as the closed half space being tangent to $\partial(V^+_{{n}})$ at  $L_{n}(\xi):=\{(s e^{i \xi},(1-c_{n}) s) \;:\;s \in \r\}$ and containing $V_{n}^+.$ Consider an increasing divergent sequence $\{r_k\}_{k \in \n}$ in $]1,+\infty[$  such that  $F_{n} \cap \pi^{-1}(\{x \in \Pi_0 \;:\; \|x\|_0=r_ k\}) \subset V_{n}^+.$ Set $F'_{n}:=F_{n}-V^+_{{n}}$ and let $C$ be an arbitrary connected component of $F'_{n} \cap \pi^{-1}(\{x \in \Pi_0\,:\, \|x\|_0\geq r_1\}).$ Obviously, $C$ is compact and $\partial(C) \subset \partial(V^+_n) \cup 
\{(y,0,t)\,:\, y,\,t \geq 0\}\subset H_{n}(\xi),$ for any $\xi \in I.$ 
The convex hull property implies that $C \subset H_{n}(\xi),$ and by a standard application of the maximum principle,  $\partial(H_{n}(\xi)) \cap C=\emptyset$ for any $\xi \in I.$ Since this is valid for any connected component $C$ of $F'_{n}\cap \pi^{-1}(\{x \in \Pi_0\,:\, \|x\|_0\geq r_1\})$ and $\xi \in I,$ we infer that $F_{n}'':=F_{n} \cap \pi^{-1}(\{x \in \Pi_0 \;:\;\|x\|_0\geq r_1,\; \arg(x) \in I\}) \subset V^+_{n}$ and $\tau^+(F_{n}'')\leq c_{n}.$ However $\lim_{n \to \infty} \tau_0^+(F_n)=0$ implies that  $\{c_n\}_{n \in \n} \to 0,$ hence $\lim_{n \to \infty} \tau^+(F''_n)=0$ too. We infer that $\Xi(\pi)=0,$ contradicting Lemma \ref{lem:universal} and proving the corollary. \end{proof}

The existence of lightlike rays in a maximal surface imposes some restrictions on its geometry. We start with the following lemma.
\begin{lemma} \label{lem:ray}
Let $W \subset {\cal R}$ be a  wedge of angle $\theta \in ]0,+\infty[,$  write $\partial(W)=L_1 \cup L_2,$ where $L_1$ and $L_2$ are divergent arcs with the same initial point, and call $\theta_j=\lim_{x \in L_j} \arg(x),$ $j=1,2.$ Let $c_0 \subset W$ be an arc such that $z(c_0)$ is a half line and $\lim_{x \in z(c_0)\to \infty} \arg(x)=\xi \in ]\theta_1,\theta_2[.$

If $X:W \to \l^3$ is a maximal multigraph and $c:=X(c_0)$ is a lightlike ray then $c$ is sublinear with direction $v_{\xi}\in \{\frac{1}{\sqrt{2}}(e^{i \xi},1),\frac{1}{\sqrt{2}}(e^{i \xi},-1)\}.$
\end{lemma}
\begin{proof}  Any blow-down of $c$ with center $O$ is a lightlike half line in ${\cal C}_0$ (that is to say, if $\{\mu_n\}_{n \in \n} \to 0,$ $\mu_n>0,$ then $\{\mu_n c\}_{n \in \n} \to l,$ $l \subset {\cal C}_0$). It suffices to consider the blow-down sequence $\{\mu_n X\}_{n \in \n}$ of maximal multigraphs and  take into account Proposition \ref{pro:converge}. \end{proof}

\begin{proposition}\label{pro:wedge}
Let  $N \subset \l^3$ be a properly embedded maximal multigraph of finite angle $\theta>0,$ and assume that $\partial(N)$ can be split into two  proper sublinear  arcs $l_1$ and $l_2 $  with lightlike directions $v_1$ and $v_2,$ respectively.

Then $v_1=\pm v_2,$ and   $\lim_{x \in N \to \infty}
g(x)=\mbox{st}_0(w),$  where $g$ is the holomorphic Gauss map of
$N$ and  $w=-\frac{1}{t(v_1)}v_1.$  In particular, the underlying
complex structure of $N$ is parabolic.

\end{proposition}
\begin{proof}   Up to removing a suitable compact  subset suppose that $O \notin \pi(N),$ and as usual call $\arg:N \to \r$ a branch of the argument of $\pi|_N.$  Write $\theta_j=\lim_{x \in l_j \to \infty}\arg(x),$ $j=1,2,$ and suppose without loss of generality that $\theta_1<\theta_2.$ Fix a compact arc $l_0 \subset \partial(N).$
From the definition of multigraph, it is not hard to construct a foliation ${\cal F}(\xi,u):[\theta_1,\theta_2]\times [0,+\infty[ \to N$ satisfying:
\begin{enumerate}[(i)]
\item $l_\xi:={\cal F}(\xi,\cdot)$ is a proper arc with initial point in $l_0,$ for any $\xi \in ]\theta_1,\theta_2[,$ and $l_{\theta_j}=l_j$ up to a compact subset, $j=1,2.$
\item For any $\epsilon>0,$ there is a closed disc  $D(\epsilon)\subset \Pi_0$ such that $\pi(l_\xi)-D(\epsilon)$ is a half line, for any $\xi \in ]\theta_1+\epsilon,\theta_2-\epsilon[.$
\item $u$ is the Euclidean arclength parameter of $\pi(l_\xi),$ $\xi \in [\theta_1,\theta_2],$ and $$\lim_{u \to \infty} \max\{|\arg({\cal F}(\xi,u))-\xi|\;:\; \xi \in [\theta_1,\theta_2]\}=0.$$
\end{enumerate}

Let $F \subset [\theta_1,\theta_2]$ be the closure of $F_0:=\{\xi \in [\theta_1,\theta_2] \;:\; l_\xi\; \mbox{is a lightlike ray}\}.$ Since blow-downs of lightlike rays are lightlike half straight lines and $F_0$ is dense in $F,$  any blow-down of $N_F:={\cal F}(F \times [0,+\infty[)$ with center the origin is a closed countable collection of angular regions\footnote{$W \subset {\cal C}_0$ is said to be an angular region if either $W=\pi^{-1}(W_{\theta}) \cap {\cal C}_0^+$ or $W=\pi^{-1}(W_{\theta}) \cap {\cal C}_0^-,$ where $\theta \in [0,2 \pi].$} in ${\cal C}_0$ (some of them could be lightlike rays).  This argument and   Proposition \ref{pro:converge} (see also Lemma \ref{lem:ray}) show that $l_\xi$ is a sublinear arc with lightlike direction $v_{\xi}\in \{\frac{1}{\sqrt{2}}(e^{i \xi},1),\frac{1}{\sqrt{2}}(e^{i \xi},-1)\},$  for every $\xi \in F.$

Let us see that $F$ is a compact totally disconnected set. Reason by contradiction and suppose there exists a closed interval $J \subset F$ of length $|J|>0.$ Then, any blow-down of $N_J:={\cal F}(J \times [0,+\infty[)$ with center $O$ is an angular region of ${\cal C}_0$ of positive angle, and thus $\tau(N_j)=0.$ This contradicts Lemma \ref{lem:universal} and proves our assertion.

\begin{quote}
{\bf Claim 1:} {\em If  $\{p_n\}_{n \in \n} \subset N$ is  divergent  and  $\lim_{n \to \infty} \arg(p_n) =\xi \in F,$  then $\lim_{n\to \infty} g(p_n)=\mbox{st}_0(w_\xi),$ where $w_\xi=-\frac{1}{t(v_\xi)} v_\xi.$}
\end{quote}
\begin{proof} Call $\mu_n:=\mbox{d}(p_n,l_{\xi})$ and  take $q_n \in l_{\xi}$ satisfying $\|p_n-q_n\|_0=\mu_n,$ $n \in \n.$ Set $G:= \arg^{-1}([\xi-\delta,\xi+\delta]),$ $\delta \in ]0,\pi[,$ and put $G_n:=\frac{1}{\lambda_n}(G-q_n),$ where $\lambda_n:=\max\{\mu_n,1\},$ for any $n \in \n.$

Since $\{\frac{q_n}{\lambda_n}\}_{n \in \n}$ is divergent,
$\{G_n\}_{n \in \n}$ converges in the ${\cal C}^0$-topology to
either an entire graph over $\Pi_0$ containing a lightlike straight line parallel to $v_\xi$ (if $\xi \notin
\{\theta_1,\theta_2\}$) or a graph over a closed lightlike half plane $H
\subset \Pi_0$  with boundary parallel to $v_{\theta_j},$ $j \in \{1,2\}$ (if $\xi \in \{\theta_1,\theta_2\}$). Anyway,  Lemma \ref{lem:basico}
gives that  $G_\infty$ is  either a lightlike plane or a
lightlike half plane bounded by a lightlike line. The
claim follows from Proposition \ref{pro:converge}. \end{proof}

The closure  of a connected component of $[\theta_1,\theta_2]-F$ is defined to be a {\em good} component of $[\theta_1,\theta_2].$ As above, if $I$ is a good component we set $N_I:={\cal F}(I \times [0,+\infty[).$

\begin{quote}
{\bf Claim 2:} {\em If $I=[\xi_1,\xi_2]$ is a good component of $[\theta_1,\theta_2],$  then the limit $w_I:=\lim_{x \in N_I \to \infty} g(x)$ exists. In particular, $w_I=w_{\xi_1}=w_{\xi_2}.$}
\end{quote}

\begin{proof} Define ${\cal H}_1=\{x \in \l^3\;:\; \|x\|\leq 1\},$ and let us show that $l_\xi \cap {\cal H}_1$ is compact for every $\xi \in ]\xi_1,\xi_2[.$
It suffices to check that any divergent subarc $l_\xi'\subset l_\xi$ satisfying that $\pi(l_\xi') \subset \{s e^{i \xi} \;:\; s \geq 0\}$ intersects ${\cal H}_1$ in a compact set. Assume that $l_\xi' \cap {\cal H}_1\neq \emptyset$ (otherwise we are done), and note that $l_\xi'$ can not lie in ${\cal H}_1$ because otherwise $l_\xi'$ would be a lightlike ray. Therefore,  $l_\xi' \cap \partial({\cal H}_1)\neq \emptyset.$ Since $l_\xi'$ has slope $<1$ and the hyperbola $\pi^{-1}(\{s e^{i \xi} \;:\; s \geq 0\}) \cap {\cal H}_1$ is timelike, then $l_\xi'$ and $\partial({\cal H}_1)$ meet only once, proving our assertion.

Thus we can find an smooth proper arc in $c \subset N_I-{\cal H}_1,$ $c \cong ]0,1[,$ satisfying that $\arg|_c$ is monotone and $\arg(c)=]\xi_1,\xi_2[.$ Set $N_I'\subset N_I$ the simply connected region bounded by $c$ and disjoint from $\partial(N_I),$ and note that $N_I'$ has parabolic underlying conformal structure (use that $N_I' \cap {\cal H}_1=\emptyset$ and Theorem \ref{th:parabo}).
On the other hand, splitting $c$ into two divergent arcs $c_1$ and $c_2$ with the same initial point and using Claim 1, we have that, up to relabeling,   $\lim_{x \in c_j \to \infty} g(x)= w_{\xi_j},$ $j=1,$ $2.$ By Theorem \ref{th:sectorial}, $w_{\xi_1}=w_{\xi_2}$ and $\lim_{x \in N_I' \to \infty} g(x)=w_{\xi_1}.$

To finish, take a compact arc $c_0 \subset N_I$  connecting  $\partial(N_I)$ and $\partial(N_I')$ and splitting $N_I-\mbox{Int}(N_I')$ into two regions $N_I^j,$ $j=1,2$. The spacelike property guarantees that  $N_I^j$ is contained in a spacelike half space, hence  it is  parabolic by Corollary \ref{co:halfspace}. As above $\lim_{x \in N_I^j \to \infty} g(x)=w_{\xi_1},$ $j=1,2,$ completing the proof. \end{proof}

Claim 2 and a connection argument give that $w:=w_{I}$ does not depend on the good component $I$ of  $[\theta_1,\theta_2]$ and  $\lim_{x \in N \to \infty} g(x)=w.$ Since $|g|<1$ on $N,$ then $h:=-\log |g-w| +\log 2$ is a positive proper harmonic function on $N,$ proving that $N$ is parabolic and concluding the proof. \end{proof}

\begin{corollary} \label{co:wedge}
Let $N$ be as in Proposition \ref{pro:wedge}, but allowing that $\partial(N)$ contains  lightlike subarcs.
Then $v_1=\pm v_2$ and   $\lim_{x \in N-\partial(N) \to \infty} g(x)=\mbox{st}_0(w),$  where $w=-\frac{1}{t(v_1)}v_1.$

Moreover, if in addition $N$ is a  $^*$maximal surface then its underlying complex structure is parabolic.
\end{corollary}
\begin{proof}
Consider a properly embedded maximal multigraph $N'$ in $\l^3$  contained in $N-\partial(N)$ and with the same angle of $N.$ If $\partial(N')$   is close enough to $\partial(N),$ $\partial(N')$ can be also split into  two proper {\em spacelike}  sublinear subarcs $l'_1$ and $l'_2$  with  directions $v_1$ and $v_2,$ respectively, and Proposition \ref{pro:wedge}  applies to $N'.$ Since $N'$ is any arbitrary region of $N$ satisfying these conditions, the first part of the  corollary easily holds.
For the second part,  take a conformal parameterization $X:\mb \to \l^3$ of $N$ and note that $g$ has well defined limit at the end of $\mb.$ Reasoning like in Proposition \ref{pro:wedge} $\mb$ is parabolic. \end{proof}

\begin{remark} \label{re:weakly}
Any maximal multigraph $N$ satisfying the hypothesis of Proposition \ref{pro:wedge} or Corollary \ref{co:wedge} is asymptotically weakly spacelike. Indeed, simply observe that any spacelike straight line not orthogonal to $v_1$ intersects $N$ into a compact set.

\end{remark}

Now we can prove the following uniqueness result.
\begin{proposition}\label{pro:unienn}
Let $G$ be an entire PS graph which is a maximal surface except on a closed lightlike half line $l \subset G.$
Then $G$ is congruent in the Lorentzian sense to  Enneper's graph $\overline{E}_2.$
\end{proposition}
\begin{proof} Up to a Lorentzian isometry, put $l=\{(x_1,x_2,t) \in \l^3\;:\; x_1=x_2-t=0, \;t\geq 0\}.$

Set $G_0=G-l$ and $l_0={l}-\{O\}.$ From Riemann's uniformization theorem, $G_0$ is conformally equivalent to either $\c$ or $\u.$ Since $g$ is not constant and $|g|<1$ on $G_0,$ then necessarily $G_0\equiv\u.$ Label $X:\u \to \l^3$ as the associated conformal parameterization of $G_0.$

\begin{quote}
{\bf Claim 1:} {\em If $\{z_n\}_{n \in \n} \subset \u$ and $\{X(z_n)\}_{n \in \n}\to p_0 \in l_0 \cup\{\infty\}$ then  $\lim_{n \to \infty} g(z_n)=1.$}
\end{quote}
\begin{proof} Label $\lambda_n=\mbox{d}(X(z_n),l),$ and let us see that $\{g(z_n)\}_{n \in \n} \to 1$ provided that  $\lim_{n \to \infty} X(z_n)/\lambda_n =\infty.$ Indeed, the sequence $\{G_n:=\frac{1}{\lambda_n} \cdot (G-X(z_n))\}_{n \in \n}$ converges in the ${\cal C}^0$-topology to an entire PS graph $G_\infty$ containing a lightlike straight line parallel to $\{x_1=x_2-t=0\},$ hence $ G_\infty=\{x_2-t=0\}$ from Lemma \ref{lem:basico}.   The assertion follows from Proposition \ref{pro:converge}.

Applying Proposition \ref{pro:wedge} to $G - \pi^{-1}(\{|x_1|<\delta,\; x_2>-\delta\}$ for any $\delta>0,$ the claim holds. \end{proof}

Fatou's theorem guarantees that $g:\u \to \d$ has well defined angular limits a. e. on $\partial(\u)\equiv \r,$ and since $g$ is not constant,   Privalov's theorem  gives that these limits are different from $1$ a. e. on $\partial (\u).$ By  Claim 1 and a connectedness argument, we infer that any two sequences $\{z_n\}_{n \in \n},$ $\{z'_n\}_{n \in \n}$  satisfying $\lim_{n \to \infty} X(z_n),$ $\lim_{n \to \infty} X(z'_n) \in l_0 \cup\{\infty\}$ converge to the same point $z_0 \in\r\cup \{\infty\}$ (up to a conformal transformation we will suppose $z_0=\infty$). Therefore, $\lim_{z \to r} X(z)=O$ for all $r \in  \r,$ and from equation (\ref{eq:wei}) we get that $|g|=1$ on $\r\cup \{\infty\}.$   By Schwarz reflection, $X$ and $g$ extend to $\c$ and $\overline{\c},$ respectively,  and $dg \neq 0$ on $\partial(U)\cup \{\infty\}.$ The extended map $X:\c \to  \l^3$ is a conformal maximal immersion with lightlike singular set $\r$ and $X(\overline{\u})=G_0 \cup \{O\}.$

Set $u:=\left((t-x_2) \circ X\right)|_{\overline{\u}}$ and label $u^*$ as its harmonic conjugate.
\begin{quote}
{\bf Claim 2:} {\em  The holomorphic map $h:\overline{\u}\to \c,$ $h:=u+i u^*,$  is injective and $h(\overline{\u})=\{z \in \c \;:\; \mbox{Re}(z) \geq 0\}.$}
\end{quote}
\begin{proof}
From equation (\ref{eq:acausal}),  $G \subset \cap_{x \in l} \overline{\mbox{Ext}({\cal C}_{x})}\subset \{t-x_2 \geq 0\}.$ Then, the maximum principle gives that $G_0 \subset  \{t-x_2 > 0\},$ that is to say, $u^{-1}(0)=\r.$ Furthermore, as $\overline{\u}$ is parabolic and $u$ is not constant (see Section \ref{sec:parabo}), then  $u$ is non negative and unbounded.

Basic theory of harmonic functions says that $u^{-1}(a)$ consists
of a proper family of analytical curves meeting at equal angles
at singular points of $u,$ $a \geq 0.$ Let us show that
$u^{-1}(a)$ consists of a unique regular analytical arc, for any
$a \geq 0.$ Indeed, otherwise we can found a region $\Omega
\subset\overline{\u}$ such that $0\leq u|_{\Omega} \leq a$ and
$u|_{\partial(\Omega)}=a,$ contradicting the parabolicity of
$\Omega.$

Since $u^*|_{u^{-1}(a)}$ is one to one for any $a\geq 0,$ then $h$ is injective. Furthermore, $h(\overline{\u})$ is parabolic simply connected {\em open} subset of $\{z \in \c \;:\; \mbox{Re}(z) \geq 0\},$ and so $h(\overline{\u})=\{z \in \c \;:\; \mbox{Re}(z) \geq 0\},$ which proves the claim. \end{proof}

Up to a conformal transformation, we have $h(z)=i B z,$ $B \in \r-\{0\},$ $B<0,$ and since $dh=i\phi_3-\phi_2=-i \frac{(g-1)^2}{2 g} \phi_3,$ then $\phi_3=-B \frac{2 g}{(g-1)^2} dz.$ As $G$ has a unique topological end, then $g^{-1}(1)=\infty.$ Moreover, $dg\neq 0$ along $\r \cup \{\infty\}$ gives that $g|_{\r \cup \infty}$ is one to one, and so $g(z)=(z-i r)/(z+i r),$ where  $r \in \r$ and  $|(1-r)/(1+r)|<1.$ Up to conformal reparameterizations, Lorentzian isometries and homotheties, these are the Weierstrass data of $E_2.$ \end{proof}

\subsection{Finiteness of maximal graphs with planar boundary}
Let $\Omega$ be a region in $\r^2.$ A non flat maximal graph $G=\{(x,u(x)) \;:\; x \in \Omega\}$ is said to be supported on $\Omega$ if $u|_{\Omega-\partial(\Omega)}$ satisfies equation (\ref{eq:max}) and  $u=0$ on $\partial (\Omega)$ (in particular, $\Omega$ can not be compact).

Assume that $G=\{(x,u(x)) \;:\; x \in \Omega\}$ is supported on $\Omega$ and denote by $G(R)$ (resp., $\Omega(R)$) the intersection $G \cap (D(R)\times \r)$ (resp., $\Omega \cap D(R)$), where $D(R)=\{x \in \r^2 \;:\; \|x\|_0 \leq R\},$ $R>0.$
Let $A(G(R))$ denote the area of $G(R)$ computed with the Riemannian metric induced by $\langle ,\rangle$ on $G.$ The spacelike condition $\|\nabla u \|_0<1$ gives the following trivial estimate:
\begin{equation}\label{eq:li-wang1}
A(G(R)) = \int_{\Omega(R)} \sqrt{1-\|\nabla u\|_0^2} \;da\leq A_0(\Omega (R)) \leq \pi R^2
\end{equation}
where $da$ is the Euclidean area element in $\r^2$ and $A_0(\Omega(R))$ is the Euclidean area of $\Omega(R)$ in $\r^2.$

The following theorem has been inspired by Li-Wang work \cite{li-wang}.
\begin{theorem} \label{th:li-wang2}
Let $\{G_i\}_{i=1}^k$ be a set of $k$ maximal simply connected graphs in $\l^3$ defined by the functions $\{u_i\}_{i=1}^k$
with disjoint supports
 $\{\Omega_i\}_{i=1}^k$ in $\r^2.$ Let us assume that $\|\nabla u_i \|_0\leq 1-\varepsilon,$ for any $i=1,\ldots,k,$ where $\varepsilon\in ]0,1[.$

Then $k\leq \frac{8}{\varepsilon(2-\varepsilon)}.$
\end{theorem}
\begin{proof}
Without loss of generality, suppose $O \notin \cup_{i=1}^k \Omega_i.$ Since $O \notin \Omega_i,$ $u_i|_{\partial(\Omega_i)}=0$  and $\|\nabla u_i\|_0\leq 1-\varepsilon,$ we get from (\ref{eq:acausal}) that $|u_i(x)|\leq (1-\varepsilon)\|x\|_0 $ on $\Omega_i,$ $i=1,\ldots,k.$

Fix $R_0>0$ such that $G_i(R_0) \neq \emptyset$ for $i=1,\ldots,k.$
As $|u_i|\leq (1-\varepsilon) R$ on $G_i(R)$ and $A(G_i(R))\leq \pi R^2$ for any $R\geq R_0$ and $i\in \{1,\ldots,k\},$ then for any $m \geq 1$ we obtain
$$\prod_{j=0}^m \prod _{i=1}^k \frac{\int_{G_i (2^j R_0)}|u_i| dx}{\int_{G_i (2^{j+1} R_0)}|u_i| dx}=\prod _{i=1}^k \frac{\int_{G_i (R_0)}|u_i| dx}{\int_{G_i (2^{m+1} R_0)}|u_i| dx}\geq \alpha ((1-\varepsilon) 2^{3(m+1)}\pi R_0^3)^{-k},$$ where $dx$ is the intrinsic area element and $\alpha=\prod _{i=1}^k \int_{G_i (R_0)}|u_i| dx.$

Hence, there exists $0\leq t\leq m$ such that $\prod _{i=1}^k \frac{\int_{G_i (2^t R_0)}|u_i| dx}{\int_{G_i (2^{t+1} R_0)}|u_i| dx}\geq 2^{-3k}\alpha^{1/(m+1)} ((1-\varepsilon) \pi R_0^3)^{-k/(m+1)},$ and by the arithmetic means we infer that
\begin{equation} \label{eq:prime}
 \sum _{i=1}^k \frac{\int_{G_i (2^t R_0)}|u_i| dx}{\int_{G_i (2^{t+1} R_0)}|u_i| dx}\geq k  ( \prod _{i=1}^k \frac{\int_{G_i (2^t R_0)}|u_i| dx}{\int_{G_i (2^{t+1} R_0)}|u_i| dx} )^{1/k} \geq 2^{-3} k \alpha^{\frac{1}{(m+1)k}} ((1-\varepsilon)\pi R_0^3)^{-1/(m+1)}
 \end{equation}
On the other hand, labeling $M_i= \max\{\frac{|u_i(x)|}{\int_{G_i (2^{t+1} R_0)}|u_i| dx}\;:\; x \in  G_i (2^t R_0) \},$ $i=1,\ldots,k,$  we have
$$\sum _{i=1}^k \frac{\int_{G_i (2^t R_0)}|u_i| dx}{\int_{G_i (2^{t+1} R_0)}|u_i| dx} \leq \max\{M_i\;:\; i=1,\ldots,k \} \large( \sum_{i=1}^k A(G_i(2^t R_0)) \large).$$ Using that  $u_i$'s are disjointly supported and equation (\ref{eq:li-wang1}),  we have $\sum_{i=1}^k A(G_i(2^t R_0)) \leq \sum_{i=1}^k A_0(\Omega_i(2^t R_0)) = A_0(\cup_{i=1}^k \Omega_i(2^t R_0)) \leq \pi(2^t R_0)^2,$ and so
\begin{equation} \label{eq:segun}
 \sum _{i=1}^k \frac{\int_{G_i (2^t R_0)}|u_i| dx}{\int_{G_i (2^{t+1} R_0)}|u_i| dx} \leq \frac{|u_{i_0}(x_0)|}{\int_{G_{i_0} (2^{t+1} R_0)}|u_{i_0}| dx} \pi (2^t R_0)^2,
\end{equation}
for some $1 \leq i_0 \leq k$ and $(x_0,u_{i_0}(x_0)) \in G_{i_0}(2^t R_0).$ 

As $\Omega_{i_0}$ is simply connected, then the conjugate minimal graph $G_{i_0}^*=\{(x,u_{i_0}^*(x)), \;:\; x \in \Omega_{i_0}\}$ of  $G_{i_0}$ is well defined (see equation (\ref{eq:conjugada})). Since $u_{i_0}$ is harmonic  on $G_{i_0}^*$ and vanishes on $\partial (G_{i_0}^*),$ the mean value property for subharmonic functions on minimal surfaces gives
$$|u_{i_0}(x_0)|\leq \frac{1}{\pi (2^t R_0)^2} \int_{G_{i_0}^* \cap B(p_0,2^t R_0)} |u_{i_0}(y)|d^*y,$$ where $p_0=(x_0,u_{i_0}^*(x_0)),$ $B(p_0,2^t R_0)$ is the Euclidean ball of radius $2^t R_0$ centered at $p_0,$ and $d^* y$ is the Euclidean area element on $G_{i_0}^*.$

From  equation (\ref{eq:conjugada}), $d^*y=\sqrt{1+\|\nabla u_{i_0}^*\|^2}da= \frac{1}{\sqrt{1-\|\nabla u_{i_0}\|_0^2}} da \leq \frac{1}{\varepsilon(2-\varepsilon)} dx.$ Since $G_{i_0}^* \cap B(p_0,2^t R_0) \subset D(2^{t+1} R_0)\times \r,$ we deduce that
$$|u_{i_0}(x_0)|\leq  \frac{1}{\varepsilon(2-\varepsilon)} \frac{1}{\pi (2^t R_0)^2} \int_{G_{i_0}(2^{t+1} R_0)} |u_{i_0}(y)|dx,$$ and from
(\ref{eq:segun}),
$\sum _{i=1}^k \frac{\int_{G_i (2^t R_0)}|u_i| dx}{\int_{G_i (2^{t+1} R_0)}|u_i| dx} \leq  \frac{1}{\varepsilon(2-\varepsilon)}.$ Equation (\ref{eq:prime}) gives
$$ 2^{-3}k\alpha^{\frac{1}{(m+1)k}} ((1-\varepsilon)\pi R_0^3)^{-1/(m+1)} \leq \frac{1}{\varepsilon(2-\varepsilon)},$$ and taking the limit as
 $m \to \infty,$ we obtain $k  \leq   \frac{8}{\varepsilon(2-\varepsilon)},$ concluding the proof. \end{proof}


\section{Maximal surfaces with connected lightlike boundary of mirror symmetry} \label{sec:main}

Let us go over some basic definitions and properties, thereby fixing some  notations and conventions.

Throughout this section, $\mb$ be a Riemann surface whose boundary  consists of a non compact analytical arc $\Gamma,$ and $X:\mb \to \l^3$  a conformal proper  $^*$maximal embedding.  As usual we identify $\mb\equiv X(\mb)\subset \l^3$ and $\Gamma \equiv X(\Gamma)\subset \l^3.$

From Lemma \ref{lem:prime}, $\pi|_\mb:\mb \to\Pi_0$ is a local embedding. Moreover, since $\Gamma$ is a regular lightlike arc then its orthogonal projection over the $t$-axis is one to one and we can write   $\Gamma \equiv \{\Gamma(s):=(\gamma(s),s) \in \l^3 \,:\,s \in \r\},$ where $s:=t|_\Gamma$ is the Euclidean arclength of $\gamma \subset \r^2 \equiv \Pi_0.$  Up to a translation, we assume that $\Gamma$ passes through $q_0:=(1,0,0)$ (and so $\Gamma(0)=(1,0,0)$ and  $\gamma(0)=(1,0)$).
Since any spacelike plane and $\Gamma$ meet transversally at a unique point, then  
\begin{equation} \label{eq:normal}
\Gamma-\{q\} \subset \mbox{Int}({\cal C}_q), \quad \mbox{for any} \; q \in \Gamma.
\end{equation}

Let $g:\mb \to \overline{\d}$ denote the holomorphic Gauss map of $\mb,$ and recall that $|g|(p)=1$ if and only if $p \in \Gamma.$ Thus the argument function $\mbox{Im}(\log(g))$ has a well defined one to one branch, namely $\theta,$ along $\Gamma.$ 
Labeling $\theta^-:=\inf (\theta(\Gamma))$ and  $\theta^+:=\sup(\theta(\Gamma)),$ the function $\theta(s):\r \to ]\theta^-,\theta^+[,$ $\theta(s):= \theta(\Gamma(s))$ is a diffeomorphism  and provides a global parameter along $\Gamma.$ Up to a ambient isometry, we will assume that $\theta'(s)>0$ and $\theta(0)=0.$


From equation (\ref{eq:wei}), $\gamma'(s)=ig(s)=i e^{i \theta(s)},$ and so in complex notation  $\gamma(s)=1+i \int_{0}^{s} e^{i \theta(x)}  dx.$ Up to a symmetry with respect to a timelike plane, we can assume that $\theta'(s)>0,$ hence $\lim_{s \to+\infty} \theta(s)=\theta^+$ and $\lim_{s \to -\infty} \theta(s)=\theta^-.$ 
By definition, the {\em rotation number} $\theta_\mb$ of $\mb$ is the change of the tangent angle along $\gamma.$ Obviously, $\theta_\mb=\theta^+-\theta^- \in [0,+\infty].$

As a consequence of Lemma \ref{lem:prime}, $\pi:\mb \to \Pi_0$ satisfies the following path-lifting property:\\

{\em Given an interval $I,$  $a_0 \in I,$ a curve $\beta (a):I \to \Pi_0$  and $p_0 \in \pi^{-1}(\beta(a_0))\cap \mb,$ there exists a unique inextendible curve $\widetilde{\beta}:J \to \mb$ in $\mb$ such that $a_0 \in J \subset I,$ $\widetilde{\beta}(a_0)=p_0$ and $\pi \circ \widetilde{\beta}=\beta|_{J}.$} \\

The curve $\widetilde{\beta}$ is said to be the  {\em lifting of $\beta$} to $\mb$ with initial condition $\widetilde{\beta}(a_0)=p_0.$ Note that either $J=I$ or at least one of the endpoints of $\widetilde{\beta}$ lies in $\Gamma.$


For each $s \in \r,$ take  an open simply connected neighbourhood
$V_s$ of $\Gamma(s)$ in $\mb$ such that  $\Gamma \cap V_s$ is
connected and $\pi|_{V_s}:V_s \to \pi(V_s)$ is one to one. Then
label  $n(s) \in \Pi_0$ as the unit normal to $\gamma$ at
$\gamma(s)$ interior to $\pi(V_s).$ Obviously, $n(s)$ does not
depend on the chosen neighbourhood $V_s.$

Set $\alpha_s:[0,+\infty[ \to \Pi_0,$ $\alpha_s(a)=\gamma(s)+a n(s),$ and consider the lifting  $\widetilde{\alpha}_s:J \to \mb$ of  $\alpha_s$ to $\mb$ with initial condition $\widetilde{\alpha}_s(0)=\Gamma(s).$   The property $\Gamma-\{\Gamma(s)\}\subset \mbox{Int}({\cal C}_{\Gamma(s)})$ and equation (\ref{eq:acausal}) give   $\widetilde{\alpha}_s\cap \Gamma=\Gamma(s),$ hence from the properness of $\widetilde{\alpha}_s$ we have $J=[0,+\infty[.$

On the other hand, $\gamma'(s)=ig(s):=i g(\gamma(s))$ implies that $n(s)=\pm g(s)$ (this ambiguity will be solved in the next Lemma). Taking into account that $g'(s)\neq 0$ for all $s \in \r,$ we deduce that $\gamma$ is {\em locally convex}, and so   $\gamma \cap \pi(V_{s})$ lies at one side of the tangent line  $r_s$ of $\gamma$ at $\gamma(s).$ Let $P_s^+$ and $P_s^-$ denote the two closed half  planes in $\Pi_0$ bounded by $r_s,$ and up to relabeling suppose $\gamma \cap \pi(V_{s}) \subset P_s^-.$

\begin{figure}[htpb]
\begin{center}\includegraphics[width=.4\textwidth]{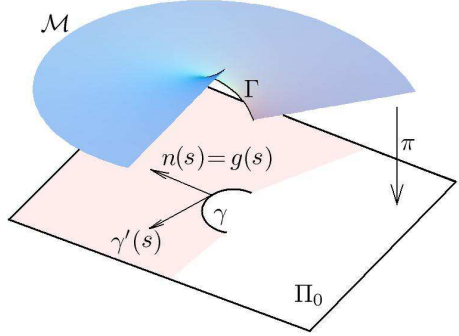}
\caption{$\Gamma,$ $\gamma,$ $\gamma'(s),$ $g(s)$ and
$n(s).$}\label{fig:primera}\end{center}
\end{figure}

\begin{lemma} \label{lem:param}
The normal vector $n(s)$ points to $P_s^+$ and $n(s)=g(s)=-i \gamma'(s).$

As a consequence, $F_0:\r \times [0,+\infty[ \to \mb,$ $F_0(s,a):=\widetilde{\alpha}_s (a)$ is a diffeomorphism.
\end{lemma}
\begin{proof} Reason by contradiction and suppose there exists $s_0 \in \r$ such that $n(s_0)$ points to $P_{s_0}^-.$ By a connection argument,  $n(s)$ points to $P_{s}^-$ for any $s\in \r.$ Take $V_{s_0}$ as above, and observe that  without loss of generality we can suppose that $\pi(V_{s_0})$ is  convex and contained in $P_{s}^-,$ for every $s$ such that $\Gamma(s) \in V_{s_0}.$ Take a segment $\zeta\subset \pi(V_{s_0})$ connecting two points  $p,$ $q \in \pi(\Gamma).$  Call $\widetilde{\zeta}\subset V_{{s_0}}$ its corresponding lifting, and observe that $\widetilde{\zeta}$ connects two points $\widetilde{p},$ $\widetilde{q} \in \Gamma \cap V_{{s_0}}.$ The spacelike property gives $|t(q)-t(p)|<\|p-q\|_0,$ which contradicts  equation (\ref{eq:normal}) and proves that $n(s)$ points to $\Pi_s^+$ for any $s.$ As a consequence of the convexity of $\gamma,$ $n'(s)= \kappa (s) \gamma'(s),$ $\kappa(s)>0.$ Taking into account that $n(s)=\pm g(s),$ $g'(s)=\theta'(s) \gamma'(s)$ and $\theta'(s)>0,$ we get that $n(s)= g(s).$

To finish, note that $F_0$ is a local diffeomorphism (take into account Lemma \ref{lem:prime}). Hence, it suffices to check that $F_0$ is proper. Take a divergent sequence $\{(s_n,a_n)\}_{n\in \n} \subset \r \times [0,+\infty[,$ and write $p_n:=F_0(s_n,a_n),$ $n \in \n.$ 

If $\{s_n\}_{n \in \n}$ is bounded and $\{a_n\}_{n \in \n}$ is divergent, the properness of $\mb$ implies that  $\{\pi(p_n)\}_{n \in \n}$ and $\{p_n\}_{n \in \n}$ diverge. 

Assume that $\{s_n\}_{n \in \n}$ diverges and, reasoning by contradiction, suppose that $\{p_n\}_{n \in \n} \to p_0\in \l^3.$  The properness of $\mb$ gives that $p_0 \in \mb.$ Furthermore, since $\Gamma$ is a lightlike  curve (see also equation (\ref{eq:normal})) it is not hard to check that $p_0 \notin \Gamma.$  Let $V\subset \mb-\Gamma$ be a  neighbourhood of $p_0$ whose projection $\Pi(V)$ is a closed disc, and without loss of generality suppose $p_n \in V$ for any $n \in \n$ (recall that $\mb$ is properly embedded). For any $q \in \pi(V)$ and  $\eta \in \s^1$ set $\beta_{q,\eta}:\r \to \Pi_0,$ $\beta_{q,\eta}(a)=q+a \eta,$ and let $\widetilde{\beta}_{q,\eta}$ denote the lifting of $\beta_{q,\eta}$ with initial condition $\widetilde{\beta}_{q,\eta}(0)=q.$ Since no spacelike arc projecting onto a segment can connect two points of $\Gamma,$ we deduce that $\widetilde{\beta}_{q,\eta} \cap \Gamma$ consists of at most one point. The first part of the lemma gives that   $A:=\Gamma\cap \left(\bigcup_{(q,\eta) \in \pi(V) \times \s^1} \widetilde{\beta}_{q,\eta}\right)$ is a connected closed subset of $\Gamma$ and  $\theta(A)$ is a relatively closed interval  in $]\theta^-,\theta^+[$  of finite length.

 On the other hand, by the unique lifting property there are $q_n \in \pi(V)$ and $\eta_n \in \s^1$ such that $\widetilde{\beta}_{q_n,\eta_n}=\widetilde{\alpha}_{s_n}$ for any $n \in \n,$ and therefore  $\{\Gamma(s_n)\}_{n \in \n}\subset A.$ Since $\{\Gamma(s_n)\}_{n \in \n}$ is divergent, $A$ can not be compact and either  $\{\Gamma(s_n)\}_{n \in \n}\to \theta^+<+\infty$ or  $\{\Gamma(s_n)\}_{n \in \n}\to \theta^->-\infty.$ Suppose that $\{\Gamma(s_n)\}_{n \in \n}\to \theta^+<+\infty$ (the other case is similar) and note that $\lim_{n \to \infty} \Gamma'(s_n)=(e^{i \theta^+},1)$ and $\lim_{n \to \infty} n(s_n)=e^{i \theta^+}.$ Hence, the sequence of curves $\{\widetilde{\alpha}_{s_n} \,:\, n \in \n\}$ is uniformly divergent (i.e., for any compact $C \subset \r^3$ there is $n_0 \in \n$ such that $\widetilde{\alpha}_{s_n} \cap C=\emptyset$ for any $n\geq n_0$), which contradicts that $\{p_n\}_{n \in \n}\to p_0$ and concludes  the proof. \end{proof}

\begin{definition} \label{def:argu}
The submersion ${\Theta_0}:\mb \to ]\theta^-,\theta^+[,$  ${\Theta_0} (F_0(s,a)):=\theta(s),$ is said to be the argument function of $\mb.$ For any subset $I \subset \Theta_0(\mb)=]\theta^-,\theta^+[,$ we call  $\mb^I:={\Theta_0}^{-1}(I) \subset \mb,$ $\Gamma^I:=\mb^I \cap \Gamma=\theta^{-1}(I)$ and $\gamma^I:=\pi(\Gamma^I).$
\end{definition}
Note that $\partial (\mb^I)= \Gamma^I \cup \widetilde{\alpha}_{s_1} \cup \widetilde{\alpha}_{s_2},$ $ \Gamma^I \cup \widetilde{\alpha}_{s_2}$ or $ \Gamma^I \cup \widetilde{\alpha}_{s_1},$ provided that $I=]\theta(s_1),\theta(s_2)[,$ $I=]\theta^-,\theta(s_2)[$ or $I=]\theta(s_1),\theta^+[,$ respectively, where $s_1,$ $s_2 \in \r.$

An open interval $I \subset \Theta_0(\mb)$ is said to be {\em good} if  $|I|\leq\pi.$ A good interval $I$ is said to be a  {\em tail interval} if one of its endpoints lies in $\{\theta^-,\theta^+\}.$

Let $I=]\theta(s_1),\theta(s_2)[\subset \Theta_0(\mb)$ be a good interval, where $s_1,$ $s_2 \in \r.$ Since the change of the tangent angle along $\gamma^I$ is  $\leq \pi,$ then $\gamma^I$ is an embedded  arc. From  Lemma \ref{lem:param}, the arcs ${\alpha}_{s_1},$ $\gamma^I$ and ${\alpha}_{s_2}$ are laid end to end and form an embedded divergent arc. Furthermore, the family $\{\alpha_s(]0,+\infty[),\;s \in ]s_1,s_2[\}$ foliates  the domain $\Omega^I\subset \Pi_0$ bounded by  $\pi(\partial (\mb^I))$ and with  interior normal $n$ along $\gamma^I.$ Thus  $\overline{\Omega^I}$ is a wedge of angle $\theta(s_2)-\theta(s_1)$ and  $\overline{\mb}^I$ is a maximal graph over $\overline{\Omega^I}.$

The case when $I$ is a  tail interval admits a similar discussion. First define 
$$\gamma'(-\infty):=\lim_{s \to -\infty} \gamma'(s)=-ie^{i \theta^-} \;\;\mbox{and}\;\; \gamma'(+\infty):=\lim_{s \to +\infty} \gamma'(s)=ie^{i \theta^+} $$
provided that $\theta^->-\infty$ and $\theta^+<+\infty,$ respectively.
 If $I=]\theta^-,\theta(s_2)]$ (resp., $I=]\theta(s_1),\theta^+[$),  then $\gamma^I$ is a sublinear arc with direction $-\gamma'(-\infty)$  (resp., $\gamma'(+\infty)$),  and $\overline{\Omega^I}$ is the wedge of angle $\theta(s_2)-\theta^-+\frac{\pi}{2}$ (resp.,  $\theta^+-\theta(s_1)+\frac{\pi}{2}$) bounded by  $\gamma^I \cup \alpha_{s_2}$ (resp., $\alpha_{s_1} \cup \gamma^I$). If $I=\Theta_0(\mb)$ is a good interval, $\mb^I=\mb$ and $\overline{\Omega^I}=\pi(\mb)$ is a wedge of angle $\theta^+-\theta^-+\pi.$

These facts have been summarized in the following lemma:
\begin{lemma} \label{lem:sheet}
If $I \subset \Theta_0(\mb)$ is a  good interval then $\Omega^I:=\pi(\mb^I-\Gamma^I)$ is a planar domain bounded by the Jordan arc 
$\pi(\partial(\mb^I)).$ Moreover, $\overline{\Omega^I}$ is a wedge of angle $\theta(s_2)-\theta(s_1),$ $\theta(s_2)-\theta^-+\frac{\pi}{2},$  $\theta^+-\theta(s_1)+\frac{\pi}{2}$ or $\theta^+-\theta^-+\pi,$ provided that $I=]\theta(s_1),\theta(s_2)[,$  $I=]\theta^-,\theta(s_2)[,$ $I=]\theta(s_1),\theta^+[$ or $I=\Theta_0(\mb),$ respectively, and  $\pi:\overline{\mb^I} \to \overline{\Omega^I}$ is bijective.

\end{lemma}
In the sequel we write $\mb^I=\{(x,u^I(x)) \;:\; x \in \Omega^I\cup \gamma^I\},$ for any good interval $I \subset \Theta_0(\mb).$

\subsection{The blow-down multigraph of $^*$maximal surfaces with connected boundary}
Fix a sequence of positive real numbers $\{\lambda_j\}_{j \in \n}$ satisfying $\lim_{j \to +\infty} \lambda_j=0,$ and consider the associated blow-down sequence of shrunk surfaces $\{\mb_j:=\lambda_j \cdot \mb, \,j \in \n\}.$

From the conformal point of view, $\mb_j=\mb$ and   both Riemann surfaces have the same holomorphic Gauss map $g.$ Therefore, we can choose the same branch $\theta$ of  $\mbox{Im} (\log(g))$ along $\Gamma_j:=\lambda_j \cdot \Gamma.$ We also denote by $\gamma_j:=\pi \circ \Gamma_j$ and observe that $\gamma_j=\lambda_j \cdot  \gamma.$ Lemma \ref{lem:param} applies to $\mb_j$ and the corresponding diffeomorphism $F_j:\r \times [0,+\infty[ \to \mb_j$ is now given by
$$F_j(\lambda_j(s, a))=\lambda_j F_0(s,a), \quad (s,a) \in \r \times [0,+\infty[.$$ Likewise we define the argument function $\Theta_j:\mb_j \to ]\theta^+,\theta^+[$ and the objects $\mb_j^I,$ $\Gamma_j^I$ and $\gamma_j^I,$ for any $I \subset \Theta_j(\mb_j)=\Theta_0(\mb).$  It is obvious that $\Theta_j(\lambda_j p)={\Theta_0} (p)$ for any $p \in \mb,$ and $\mb_j^I=\lambda_j \cdot \mb^I,$ $\Gamma_j^I=\lambda_j \cdot \Gamma^I$ and $\gamma_j^I=\lambda_j \cdot \gamma^I.$ If $I$ is a good interval, Lemma \ref{lem:sheet} gives that $\mb_j^I$ is a maximal graph over the wedge $\Omega_j^I:=\pi(\mb_j^I)=\lambda_j \cdot \Omega^I$ as well.

\begin{lemma} \label{lem:domlim}
If $I \subset \Theta_0(\mb)$ is a good interval, then $\{\Omega_n^I\}_{n \in \n}$ converges in the Hausdorff distance to the domain $\Omega_\infty ^{I^*} \subset \Pi_0$ given by:

\begin{enumerate}[(i)]
\item If ${I}=]\theta(s_1),\theta(s_2)[,$  $s_1,$ $s_2 \in \r,$ then $I=I^*$ and $\Omega_\infty ^I$ is the interior of the wedge of angle $\theta(s_2)-\theta(s_1)$ with boundary $\alpha^0_{s_1} \cup \alpha^0_{s_2},$ where $\alpha^0_{s_j}(a)=a n(s_j)$ for any $a \in [0,+\infty[,$ $j=1,2.$
\item If $I=]\theta^-,\theta(s_2)[,$ then $I^*=]\theta^--\frac{\pi}{2},\theta(s_2)[$ and $\Omega_\infty ^{I^*}$ is the interior of the wedge of angle $\theta(s_2)-\theta^-+\frac{\pi}{2}$ with boundary $\alpha^0_{-\infty}\cup \alpha^0_{s_2},$ where $\alpha^0_{-\infty}(a)=-a\, \gamma'(-\infty)$ for any $a \in [0,+\infty[.$
\item If $I=]\theta(s_1),\theta^+[,$ then $I^*=]\theta(s_1),\theta^++\frac{\pi}{2}[$ and $\Omega_\infty ^{I^*}$ is the interior of the wedge of angle $\theta^+-\theta(s_1)+\frac{\pi}{2}$ with boundary $\alpha^0_{s_1}\cup \alpha^0_{+\infty},$ where $\alpha^0_{+\infty}(a)=a\gamma'(+\infty)$ for any $a \in [0,+\infty[.$
\item  If $I=\Theta_0(\mb),$ then $I^*=]\theta^--\frac{\pi}{2},\theta^++\frac{\pi}{2}[$ and $\Omega_\infty ^{I^*}$ is the interior of the wedge of angle $\theta^+-\theta^-+\pi$ with boundary $\alpha^0_{-\infty} \cup \alpha^0_{+\infty}.$
\end{enumerate}
Moreover, up to  subsequences $\{\mb_j^I-\Gamma_j^I\}_{j \in \n}\to \mb_\infty^{I^*}$  in the ${\cal C}^0$-topology, where $\mb_\infty^{I^*}$ is a PS graph over $\Omega_\infty ^{I^*}$ with boundary given by:
\begin{enumerate}[(i)']
\item If  ${I}=]\theta(s_1),\theta(s_2)[$ then $\partial(\mb_\infty^{I^*})=\widetilde{\alpha}^0_{s_1} \cup \widetilde{\alpha}^0_{s_2},$ where $\widetilde{\alpha}^0_{s_j}$ is a divergent arc with initial point $O$ and projecting in a one to one way onto ${\alpha}^0_{s_j},$ $j=1,2.$
\item If $I=]\theta^-,\theta(s_2)[$ then $\partial(\mb_\infty^{I^*})=\widetilde{\alpha}^0_{-\infty} \cup \widetilde{\alpha}^0_{s_2},$ where
 $\widetilde{\alpha}^0_{-\infty}:=\{-a \Gamma'(-\infty) \,:\, a \in [0,+\infty[\}$ and $\Gamma'(-\infty)=\lim_{s \to -\infty} \Gamma'(s)=\frac{1}{\sqrt{2}}(ie^{i \theta^-},1)$ and $\widetilde{\alpha}^0_{s_2}$ is given as above.
\item If $I=]\theta(s_1),\theta^+[$ then $\partial(\mb_\infty^{I^*})=\widetilde{\alpha}^0_{+\infty} \cup \widetilde{\alpha}^0_{s_1},$ where
 $\widetilde{\alpha}^0_{+\infty}:=\{a \Gamma'(+\infty) \,:\, a \in [0,+\infty[\}$ and $\Gamma'(+\infty)=\lim_{s \to +\infty} \Gamma'(s)=\frac{1}{\sqrt{2}}(ie^{i \theta^+},1)$ and $\widetilde{\alpha}^0_{s_1}$ is given as above.
\item If $I=\Theta_0(\mb)$ then    $\partial(\mb_\infty^{I^*})=\widetilde{\alpha}^0_{-\infty} \cup \widetilde{\alpha}^0_{+\infty}.$
\end{enumerate}
\end{lemma}
\begin{proof}
For any $s \in \r$ write $\alpha_{s,j}:=\lambda_j \alpha_{s}.$

If $I=]\theta(s_1),\theta(s_2)[,$ we have $\partial(\Omega_j^I)=\alpha_{s_1,j} \cup \gamma_j^I \cup \alpha_{s_2,j},$ hence $\lim_{j \to \infty}\mbox{d}_H(\partial(\Omega_j^I),\partial(\Omega_\infty ^I))=0.$ This obviously implies that $\mbox{d}_H(\Omega_j^I,\Omega_\infty ^I)\to 0$ as $j \to \infty$ and proves $(i).$ When $I=]\theta^-,\theta(s_2)[,$ $\partial(\Omega_j^I)=\gamma_j^I \cup \alpha_{s_2,j}.$ Taking into account that the divergent arc $\gamma_j^I$ is sublinear with direction $- \gamma'(-\infty),$ $\{\mbox{d}_H(\Omega_j^I,\Omega_\infty ^I)\}_{j \in \n}\to 0$ as well, proving $(ii).$ The cases $(iii)$ and $(iv)$ are similar.

From Remark \ref{re:convergencia}, and up to taking a subsequence, $\{\mb_n^I-\Gamma_j^I\}_{n \in \n}$ converges in the ${\cal C}^0$-topology  to a  PS graph $\mb_\infty^{I^*}$ over $\Omega_\infty^{I^*}$  which  can be extended continuously to $\overline{\Omega_\infty^{I^*}}.$

Item $(i)'$ is straightforward. If $I=]\theta^-,\theta(s_2)[$ then $\Gamma^I$ is a sublinear arc with direction $\Gamma'(-\infty),$ hence  $\lim_{j \to +\infty} \Gamma_j^I=\widetilde{\alpha}^0_{-\infty}$ and $\widetilde{\alpha}^0_{-\infty}\subset \partial(\mb_\infty^{I^*}),$ proving $(ii)'.$ The cases $(iii)'$ and $(iv)'$ are similar. \end{proof}

Let ${\cal I}$ denote the family of good intervals in $\Theta_0(\mb),$ and set ${\cal I}_0=\{]c,d[ \in {\cal I} \;:\; [c,d] \subset \Theta_0(\mb)\}.$ Note that $I \in {\cal I}_0$ if and only if $I=I^*$ (i.e., $I$ is not a tail interval).

Take a covering ${\cal G}:=\{I_k \;:\; k \in \n\}$ of $\Theta_0(\mb)$ by good intervals containing a tail interval $]\theta^-,b[$ provided $\theta^- \in \r,$ and likewise if $\theta^+ \in \r.$  For each $k \in \n$ and from Lemma \ref{lem:domlim}, we can find a subsequence  of $\{\mb_j^{I_k}-\Gamma_j^{I_k}\}_{j\in \n},$ which depends on $k,$ converging on compact subsets of $\Omega_\infty^{I_k^*}$ to a PS graph. A standard diagonal process leads to a subsequence, namely $\{\mb_{j(k)}\}_{k \in \n},$ such that $\{\mb_{j(k)}^{I_h}-\Gamma_{j(k)}^{I_h}\}_{k \in \n}$
converges in the ${\cal C}^0$-topology to a PS graph $\mb_\infty^{I_h^*}$ on $\Omega_\infty^{I_h^*}$ {\em for all} $h \in \n.$ Up to replacing $\{\mb_j, \;j \in \n\}$  for $\{\mb_{j(k)},\;k \in \n\},$ we can suppose that $\{\mb_j^{I_h}-\Gamma_j^{I_h}\}_{j \in \n}\to \mb_\infty^{I_h^*}$ for any $h \in \n.$ Moreover, since any $I \in {\cal I}$ can be covered by finitely many intervals in ${\cal G},$ we also have $\{\mb_j^{I}-\Gamma_j^I\}_{j \in \n} \to \mb_\infty^{I^*}$ in the ${\cal C}^0$-topology over $\Omega_\infty ^{I^*}.$ In the sequel we write $\mb_\infty^{I^*}=\{(x,u_\infty^{I^*}(x)), \; x \in \Omega_\infty^{I^*}\}$ for any $I \in {\cal I},$ and observe that $u_\infty^{I_1^*}=u_\infty^{I_2^*}$ on  $\Omega_\infty^{I_1^*} \cap \Omega_\infty^{I_2^*},$ for any  $I_1,$ $I_2 \in {\cal I}.$

Let $\hat{\Omega}_\infty^{I^*}$ denote the {\em marked} domain $\{(x,I) \;:\; x \in \Omega_\infty^{I^*}\},$ $I \in {\cal I}.$ Set
${\cal O}:=\cup_{I\in {\cal I}}\hat{\Omega}_\infty^{I^*}$ the direct sum of the topological spaces  $\{\hat{\Omega}_\infty^{I^*}, \;I\in {\cal I} \},$ and likewise define ${\cal O}_0:=\cup_{I\in {\cal I}_0}\hat{\Omega}_\infty^{I^*}=\cup_{I\in {\cal I}_0}\hat{\Omega}_\infty^{I}.$ Consider on both spaces the equivalence relation:  $(x_1,I_1) \sim (x_2,I_2)$ if and only if $J=I_1\cap I_2 \neq \emptyset$ and  $x_1=x_2 \in \Omega_\infty^{J^*}.$
\begin{definition} We set ${\cal W}:={\cal O}/\sim$ and ${\cal W}_0:={\cal O}_0/\sim$ endowed with the quotient topology.
\end{definition}
It is natural to define an argument function  $\Theta_\infty: {\cal W}_0 \to \Theta_0(\mb)$ as the "limit" of
$\{\Theta_j:\mb_j \to \Theta_0(\mb)\}_{j \in \n}.$  Indeed, take $p=[(x,I)] \in {\cal W}_0$ and an arbitrary sequence $\{p_j\}_{j \in \n}$ such that $p_j \in \mb_j^I-\Gamma_j^I$ and $\lim_{j \to +\infty} p_j=(x,u_\infty^I(x))\in \mb_\infty^I.$ Writing $p_j=\lambda_j F_0(s_j,a_j),$ we have $\Theta_j(p_j)=\Theta_0(F_0(s_j,a_j))=\theta(s_j) \in I,$ $j \in \n.$ Taking into account that  $I \in {\cal I}_0,$ we infer that the limit $s:=\lim_{j \to +\infty} s_j$ exists and depends only on $p.$ Furthermore, $\lim_{j \to +\infty} \Theta_j(p_j)= \theta(s)\in I$ and $e^{i \theta(s)}=x/\|x\|_0.$ Thus, it suffices to set $\Theta_\infty([(x,I)]):=\theta(s).$

Since ${\cal W}$ is simply connected,  $\log:{\cal W} \to \c,$ $\log([(x,I)]):=\log(x)$ has a well defined branch and the map $h:{\cal W} \to {\cal R}^*,$  $h([(x,I)])=\left(x,\log([(x,I)])\right)$ is a homeomorphism. Choose the branch of $\log:{\cal W} \to \c$ in such a way that  $\mbox{Im}(\log)|_{{\cal W}_0}=\Theta_\infty,$ that is to say,
$\mbox{Im}(\log)\left((x,\log([(x,I)]))\right)=\Theta_\infty([(x,I)])$ for any $[(x,I)] \in {\cal W}_0.$ Finally, identify ${\cal W}$ with $h({\cal W})$ via $h$ and consider ${\cal W}\subset {\cal R}.$ Up to this identification, $z(p)=x$ provided that $p=h([(x,I)])$ and $\arg|_{{\cal W}_0}=\Theta_\infty.$
Observe that ${\cal W}=\arg^{-1}(]\theta^--\frac{\pi}{2},\theta^++\frac{\pi}{2}[)$ and ${\cal W}_0=\arg^{-1}(\Theta_0(\mb)),$ and its closures $\overline{{\cal W}}$ and $\overline{{\cal W}_0}$ in ${\cal R}$   are wedges of angles $\theta_\mb+\pi$ and $\theta_\mb,$ respectively.

If $\theta^- \in \r,$ we write ${\cal W}_-:=\arg^{-1} \left([\theta^--\frac{\pi}{2},\theta^-]\right) \cup\{[0]\},$ and likewise   ${\cal W}_+:=\arg^{-1} \left([\theta^+,\theta^++\frac{\pi}{2}]\right)\cup\{[0]\}$ provided that $\theta^+ \in \r.$ Obviously, $\overline{{\cal W}}=\overline{{\cal W}_0} \cup{\cal W}_- \cup {\cal W}_+$  (here we are assuming ${\cal W}_{\pm}=\emptyset$ provided that $\theta^\pm=\pm \infty$).

Define $u_\infty:{\cal W} \to \r,$ $u_\infty(p):=u_\infty^{I^*}(z(p)),$ where $I \in {\cal I}$ is any interval satisfying $\arg(p) \in I^*,$ and call with the same name its continuous extension  to $\overline{{\cal W}}.$ It is clear that  $u_\infty([0]))=0.$ Notice that $u_\infty$ is ${\cal C}^1$ on $\overline{{\cal W}}-\{[0]\}$ having  $\|\nabla u_\infty\|_0\leq 1$ (see Proposition \ref{pro:converge}), and so the map $${\cal X}:\overline{{\cal W}} \to \l^3, \quad {\cal X}((p,u_\infty(p))):= (z(p),u_\infty(p))$$ is a PS multigraph of angle $\theta_\mb+\pi.$

\begin{definition}
${\cal X}$ is defined to be the blow-down multigraph of $\mb$ associated to the sequence $\{\lambda_j\}_{j \in \n}.$  We also say that $\mb_\infty:={\cal X}(\overline{{\cal W}})$ is the blow-down surface of $\mb$ associated to  $\{\lambda_j\}_{j \in \n}.$
\end{definition}
Since ${\cal X}([0])=O,$ equation (\ref{eq:acausal}) gives
\begin{equation} \label{eq:acausalinfi}
\mb_\infty:=\{(z(p),u_\infty(p)) \;:\; p \in \overline{{\cal W}}\} \subset \overline{\mbox{Ext}({\cal C}_0)}.
\end{equation}
Taking into account that $\mb_\infty$ is the limit set of a sequence of embedded surfaces and that $\theta(s)$ is increasing, the sheets of the multigraph ${\cal X}:{\cal W} \to \mb_\infty \subset \l^3$ are ordered by height, i.e.,
\begin{equation} \label{eq:monoto}
u_\infty(p) \geq u_\infty(q)\; \mbox{provided that}\; \arg(p)=\arg(q)+2 k \pi,\; k\geq 0.
\end{equation}

If $\theta^+=+\infty,$ (resp., $\theta^-=-\infty$) we label $u^+(x):=\sup \{u_\infty(p)\;:\; p \in z^{-1}(x)\},$ (resp., $u^-(x) :=\inf \{u_\infty(p)\;:\; p \in z^{-1}(x)\}$), for any $x \in \c-\{0\}.$ We make the continuous extension $u^+(0)=0$ (resp., $u^-(0)=0$) and call  $$\mb_\infty^+:=\{(x,u^+(x)) \;:\; x \in \c\} \;\; \mbox{(resp.,}\;\; \mb_\infty^-:=\{(x,u^-(x)) \;:\; x \in \c\})$$ the associated graph. From Remark \ref{re:convergencia} and equation (\ref{eq:acausalinfi}),  $\mb_\infty^+$ and $\mb_\infty^-$ are entire PS graphs containing the origin, and so lying in  $\overline{\mbox{Ext}({\cal C}_0)}.$

Notice that ${\cal X}|_{\arg^{-1}(I^*)}:\arg^{-1}(I^*) \to \mb_\infty^{I^*}$ is a homeomorphism for any good interval $I \in {\cal I},$ and so $\mb_\infty^{I^*}$   can  be identified with $\arg^{-1}(I^*)$ via ${\cal X}.$ Then, it is natural to set ${\cal W}^I:=\arg^{-1}(I) \subset \overline{\cal W}$ and put $\mb_\infty^I={\cal X}({\cal W}^I),$ for {\em any subset} $I \subset \arg(\overline{\cal W}-\{[0]\}).$ When $I$ is connected, the closure of ${\cal W}^I$ is a wedge and ${\cal X}|_{\overline{\cal W}^I}:\overline{\cal W}^I \to \l^3$ is a PS multigraph of angle $|I|.$

Label ${\cal A}_{{\cal X}}=|\nabla u_\infty|^{-1}(1)\subset \overline{{\cal W}}-\{[0]\}$ as the set of singular points of $u_\infty$ (see Theorem \ref{th:converge} and Remark \ref{re:strip}).
A point $\xi \in \arg(\overline{{\cal W}}-\{[0]\})$ is said to be a {\em singular} angle if $\arg^{-1}(\xi) \subset {\cal A}_{{\cal X}}.$ 
We denote by $I_{\cal X} \subset \arg(\overline{{\cal W}}-\{[0]\})$ the subset of singular angles. A singular angle $\xi \in I_{\cal X}$ is said to be {\em conical}  if ${\cal X}(\arg^{-1}(\xi))$ is a lightlike half line (that is to say, if $\arg^{-1}(\xi)$ is a singular segment of $u_\infty$). The set of conical singular angles will be denoted by $I_{\cal X}^c.$ Note that $I_{\cal X}$ and $I_{\cal X}^c$ are closed subsets of $\arg(\overline{{\cal W}}-\{[0]\}).$

\begin{proposition} \label{pro:singu}
If ${\cal A}_{{\cal X}} \neq \emptyset$ then $I_{\cal X}^c\neq \emptyset$  and ${\cal A}_{{\cal X}}={\cal W}^{I_{\cal X}}.$  Moreover:
\begin{enumerate}[(i)]
\item If $I_{\cal X}-I_{\cal X}^c \neq \emptyset$ then $\mb_{\infty}^{I_{\cal X}-I_{\cal X}^c} \subset \Sigma,$ where $\Sigma$ is the lightlike plane passing through $O,$ and $\overline{\mb_\infty^J}$ is  a  half plane in $\Sigma$  bounded by $L:=\Sigma \cap {\cal C}_0,$ for any connected component $J$ of $I_{\cal X}-I_{\cal X}^c.$
\item If $I=]\xi_2,\xi_2[ \subset \arg(\overline{\cal W}-\{0\})-I_{\cal X}$ is a bounded component then $\xi_1,$ $\xi_2\in I_{\cal X}^c$ and  $\xi_2-\xi_1=k \pi,$ $k \in \n,$ $k \geq 2.$  Moreover, ${\mb_\infty^I}$ is an embedded maximal multigraph, the 
limit tangent plane $\Sigma'$ of  ${\mb_\infty^I}$  at infinity is lightlike, and  $\Sigma'$ does 
not depend on $I.$ If in addition $I_{\cal X}-I_{\cal X}^c \neq
\emptyset,$ then $\Sigma'=\Sigma$ and $\partial(\mb_\infty^I)
\subset L.$
\item If $J_0$ is the closure of a connected component of $\mbox{Int}(I_{\cal X}^c)$ then $\arg({\cal W}-[0])-J_0$ is connected.
\end{enumerate}
\end{proposition}
\begin{proof} Assume that ${\cal A}_{\cal X} \neq \emptyset,$ and let us see that $\xi \in I_{\cal X}$ if and only if $\arg^{-1} (\xi) \cap  {\cal A}_{{\cal X}} \neq \emptyset.$ Indeed, suppose $\arg^{-1}(\xi)\cap {\cal A}_{{\cal X}} \neq \emptyset$ and take $p \in \arg^{-1} (\xi)\cap {\cal A}_{{\cal X}} .$ From Remark \ref{re:strip}  there is a divergent arc $l_p\subset \overline{\cal W}-\{[0]\}$ passing through $p$ such that $L_p:=\overline{{\cal X}(l_p)}$ is either a lightlike half line starting from $O$ or a complete lightlike straight line (in particular, $l_p \subset {\cal A}_{{\cal X}}$).

If $O \in L_p$ then $\xi \in I_{\cal X}^c\subset I_{\cal X}$ and we are done. If $O \notin L_p,$ then $p \in {\cal A}_{{\cal X}}-{\cal W}^{I_{\cal X}^c}$  and $L_p$ is a lightlike straight line.  Then, consider the open half plane $H \subset \Pi_0$ satisfying $\pi(L_p) \subset H$ and $O \in \partial(H),$ and label $V_\xi$ as the connected component of $(\pi \circ {\cal X})^{-1}(H)$ containing $L_p.$ From Lemma \ref{lem:basico}, $V_\xi$ lies in the lightlike  plane containing $L_p,$  and so $V_\xi \subset {\cal A}_{{\cal X}}.$ As  $\arg^{-1}(\xi) \subset V_\xi \subset {\cal A}_{{\cal X}}$ then $\xi \in I_{\cal X},$ proving our assertion. Furthermore, note that $\partial(V_\xi)$ is a  lightlike straight line passing through $O.$
As a consequence, ${\cal A}_{{\cal X}}={\cal W}^{I_{\cal X}}$ and  $I_{\cal X}^c \neq \emptyset.$

Assume that $I_{\cal X}-I_{\cal X} ^c\neq \emptyset$ and consider
$\xi \in I_{\cal X}-I_{\cal X} ^c.$ Let $J \subset  I_{\cal
X}-I_{\cal X} ^c$ be the connected component containing $\xi.$ With the previous notation, we
have shown that  $|J|=\pi,$ $L_J:= \partial(V_\xi) \subset
\mb_\infty^{I_{\cal X}^c}$ is a lightlike straight line and the
lightlike half plane $V_\xi=\mb_\infty^J$ is a connected
component of $\mb_\infty^{I_{\cal X}-I_{\cal X}^c}.$  Therefore,
to finish item $(i)$ it suffices to check that the plane
$\Sigma_J$ containing $\mb_\infty^{J}$ (hence $L_J$) does not
depend on $J.$ Indeed, take  two components $J_1,$ $J_2 \subset
I_{\cal X}-I_{\cal X}^c$, and simply observe that the lightlike planes
$\Sigma_{J_1}$ and $\Sigma_{J_2}$ can not meet transversally
because $\mb_\infty$ is the limit set of a sequence of embedded
surfaces. In the sequel we will call $\Sigma:=\Sigma_J$ and
$L:=L_J,$ provided that $I_{\cal X}-I_{\cal X}^c\neq \emptyset$
and $J \subset I_{\cal X}-I_{\cal X}^c$ is any connected
component.

Let us check $(ii).$ to do this, consider $I=]\xi_2,\xi_2[ \subset \arg(\overline{\cal
W}-\{0\})-I_{\cal X}$  a bounded connected component. From item $(i)$,  $\xi_1,$ $\xi_2 \in I_{\cal X}^c.$ Thus, $l_j:=\overline{{\cal X}(\arg^{-1}(\xi_j))}$ is a lightlike half line with initial point at the origin, $j=1,$ $2.$  Since $\mb_\infty$ is the limit of a sequence of embedded maximal surfaces, it is not hard to check that $\mb^I$ is embedded too. From Corollary \ref{co:wedge}, the half lines $l_1$ and $l_2$ must be  parallel, hence $\xi_2=\xi_1+k \pi,$ $k \in \n.$ Furthermore, the limit tangent plane of $\mb_\infty^I$ at infinite, namely $\Sigma'_I,$ is the lightlike plane containing $l_1 \cup l_2.$

If $k=1,$ Lemma \ref{lem:basico} implies that $\mb_\infty^I$ lies in a lightlike plane, contradicting that $]\xi_1,\xi_2[ \cap I_{\cal X}=\emptyset$ and proving that $k \geq 2.$

Finally, let $I_1,$ $I_2  \subset \arg(\overline{\cal W}-\{0\})-I_{\cal X}$ be  two components as in the statement of the claim. Reasoning as in the proof of item $(i)$, the embedded multigraphs $\mb_\infty^{J_1},$ $\mb_\infty^{J_2}$  can not meet transversally, and so $\Sigma':=\Sigma'_{I_1}=\Sigma'_{I_2}.$ Likewise  $\Sigma$ and $\mb_\infty^I$ can not meet transversally,  provided that $I_{\cal X} -I_{\cal X}^c \neq \emptyset$ and the plane $\Sigma$ makes sense, and in this case $\Sigma'=\Sigma.$

Let $J_0$ be closure of a component of $\mbox{Int}(I_{\cal X}^c).$ By a connectedness argument, either $\mb_\infty^{J_0} \subset {\cal C}_0^+$ or  $\mb_\infty^{J_0} \subset {\cal C}_0^-.$ Suppose that $\mb_\infty^{J_0} \subset {\cal C}_0^+$ (the other case is similar), and let us show that  $a:=\sup (J_0)=\theta^++\frac{\pi}{2}.$ 

Reason by contradiction and assume that $a<\theta^++\frac{\pi}{2}.$ Set $b:=\max \{a-2 \pi, \inf (J_0)\}<a,$  label $J_1=[b,a] \subset J_0$ and write 
$J_1(k)=(J_1+2 k \pi) \cap \arg({\cal W}-[0])$ for any $k \in \n.$ From equation (\ref{eq:monoto}) we get that $\mb_\infty^{J_1(k)}$ lies above $\mb_\infty^{J_1},$ and using equation (\ref{eq:acausalinfi}) we infer that $\mb_\infty^{J_1(k) } \subset {\cal C}_0^+,$  for any $k \in \n.$ So,  $J_1(k) \subset  I_{\cal X}^c$ for any $k \in \n.$ As a consequence $b>a-2 \pi,$ because otherwise $[a-2\pi,\theta^++\frac{\pi}{2}[\subset I_{\cal X}^c,$ contradicting that $J_0$ is the closure of a connected component of $I_{\cal X}^c.$  Thus $I:=]a,b+2\pi[$ is a connected component of $\arg({\cal W}-[0])-I_{\cal X}^c,$ and items $(i)$ and $(ii)$ give that $b=a-\pi,$ $I \subset I_{\cal X}-I_{\cal X}^c$ and $\mb_\infty^I$ is an open lightlike half plane bounded
by a lightlike straight line. In particular, $\mb_\infty^I$ can not lie above $\mb_\infty^{J_0},$ contradicting equation (\ref{eq:monoto}) and proving $(iii).$ \end{proof}

\subsection{$\omega^*$-maximal surfaces and the blow-down plane} \label{subsec:w*}
Recall that $\mb$ is asymptotically weakly spacelike, or simply $\omega^*$-maximal, if there is an affinely spacelike arc in $\l^3$ disjoint from $\mb.$ Although the $\omega^*$-condition is a little involved, it is connected with quite natural geometrical properties, as shown in Proposition \ref{pro:w*} below. Moreover it provides us a good control about the geometry of $\mb$ at infinity. Let start with some previous notions and comments.

Given a complete maximal surface $S$ in  $\l^3$ and a divergent curve $\alpha:[0,1[ \to S,$ Theorem \ref{th:calabi} shows that $\lim_{y \to 1} \frac{L(\alpha^y)}{L_0(\pi\circ \alpha^y)}\geq C>0,$ where $\alpha^y=\alpha|_{[0,y]},$ $L$ and $L_0$ are the intrinsic length in $S$ and $\Pi_0,$ respectively,  and $C$ is a constant not depending on $\alpha.$ With this inspiration in mind, $\mb$ is said to be {\em complete far from the boundary} if for any $\alpha:[0,1[ \to \mb$  whose projection $\pi \circ \alpha$ is a divergent arc in $\Pi_0,$ the limit $\lim_{y \to +\infty} \frac{L(\alpha^y)}{L_0(\pi \circ \alpha^y)}$ is positive.

On the other hand, we know that $\mb$ behaves like a "multigraph" in a generalized sense. It is then natural to say that $\mb$ is  {\em asymptotically strongly spacelike} if it admits gradient estimates far from the boundary, or being more precise, if for any good interval $I$ there is $\epsilon(I) \in ]0,1[$ such that $\|\nabla u^I\|_0<1-\epsilon(I)$ in the complement of a neighborhood of $\gamma^I,$ where $u^I:\Omega^I \to \r$ is the function defining the graph $\mb^I.$

\begin{proposition}\label{pro:w*}
Suppose that $\mb$ satisfies any of the following conditions: 
\begin{enumerate}[(a)]
\item  There is  an affinely spacelike arc contained in the surface.
\item There are a good interval $I$ and a real number $\epsilon \in ]0,1[$ such that $\|\nabla u^I\|_0<1-\epsilon$ in the complement of a neighbourhood of $\gamma^I.$ 
\item There is an arc $\alpha:[0,1[ \to \mb$ such that $\pi(\alpha)$ is half line and $\liminf_{y \to 1} \frac{L(\alpha^y)}{L_0(\pi \circ \alpha^y)}>0.$
\end{enumerate}
Then $\mb$ is $\omega^*$-maximal. As a consequence, if $\mb$ is complete far from the boundary or asymptotically weakly spacelike then it is $\omega^*$-maximal.
\end{proposition}
\begin{proof} If $\mb$ contains an affinely spacelike arc $\beta,$ it is easy to take an affinely spacelike  $\alpha \subset \l^3-\mb$ in a  small neighborhood of $\beta$ (recall that $\mb$ is properly embedded), proving $(a).$ As a consequence, if $\widetilde{\alpha}_s$ is affinely spacelike for some $s \in \r$ then $\mb$ is $\omega^*$-maximal. This automatically holds if $u^I$ satisfies that $\|\nabla u^I\|_0<1-\epsilon$ in the complement of a neighbourhood of $\gamma^I$ in $\Omega^I$ for some  $\epsilon \in ]0,1[,$ where $I$ is a good interval and $s \in t(\gamma^I),$ proving $(b).$ 

To see $(c),$ it suffices to check that $\alpha$ is affinely spacelike. Let $a$ be the Euclidean arclength parameter of $\pi \circ \alpha.$ Reason by contradiction, and take a divergent sequence $\{a_n\}_{n \in \n}\subset [0,+\infty[$ such that $\lim_{n \to \infty} \frac{\alpha(a_n)}{a_n}  \in {\cal C}_0.$ Put $\lim_{n \to \infty} \frac{\alpha(a_n)}{a_n} =(e^{i \xi},\pm 1)$ and without loss of generality suppose that $\xi=\lim_{a \to+\infty} \Theta_0(\alpha(a)).$ Let $I$ be a good interval containing $\xi,$ and up to removing a compact subarc  assume that $\alpha \subset \mb^I.$ 
Choose  $\lambda_n=\frac{1}{a_n},$ $n \in \n,$  and consider the blow down surface $\mb_\infty$ associated to $\{\lambda_n\}_{n \in \n}.$ For every $n \in \n,$ $y \in \lambda_n \cdot \Omega^I$ and $a \in [0,+\infty[$ write  $u^I_n(y)=\lambda_n u^I(y/\lambda_n)$ and $\alpha_n(a)=\lambda_n \alpha(a/\lambda_n).$ It is clear that ${\cal X}(\arg^{-1}(\xi))\subset \mb_\infty^I$ is a lightlike half line starting at the origin and $\lim_{n\to\infty} \alpha_n= {\cal X}(\arg^{-1}(\xi))$ in the ${\cal C}^1$ topology over $]0,+\infty[$ (see Proposition \ref{pro:converge}).  Since  $a$ is the arclength parameter of $\pi \circ \alpha_n$ too, we infer that  $\lim_{n \to \infty} \langle (\pi\circ \alpha_n)',\nabla u_n^I \circ \alpha_n\rangle_0=1$ uniformly on $[\delta,1]$ and $\lim_{n \to \infty} L(\alpha_n|_{[\delta,1]})=0,$ for any $\delta \in ]0,1[.$ However, the inequality $L(\alpha_n|_{[0,\delta]})<\delta,$ $n \in \n,$ implies that $\lim_{n \to \infty} L(\alpha_n)=\lim_{n \to \infty} \frac{L(\alpha^{a_n})}{a_n} \leq \delta$ for any $\delta \in ]0,1[.$ We deduce that $\lim_{n \to +\infty}  \frac{L(\alpha^{a_n})}{a_n} =0,$ which is absurd and proves $(c).$ \end{proof}

Next theorem is the main result of this section.

\begin{theorem} \label{th:strong}
The following statements hold:
\begin{enumerate}[(i)]
\item If $\theta_\mb<+\infty$ then $\mb$ is $\omega^*$-maximal, $\mb_\infty$ is a lightlike plane and $\theta_\mb=2k\pi,$ $k \in \n.$
\item If $\Theta_0(\mb)=\r$ and $\mb$ is $\omega^*$-maximal then $\mb_\infty$ is a non timelike plane.
\item If $\Theta_0(\mb) \neq \r,$ $\theta_\mb=+\infty$ and $\mb$ is $\omega^*$-maximal then $\overline{I_{\cal X}-I_{\cal X}^c}$  contains a non compact connected component $J_{\cal X}$ and $\mb_\infty^{J_{\cal X}}$ is a lightlike plane.  
\end{enumerate}
\end{theorem}
\begin{proof} 
Assume that $\theta_\mb<+\infty.$ Since $\mb$ is a multigraph of angle $\theta_\mb+\pi,$  Corollary \ref{co:wedge}  gives that   $\mb$ has sublinear growth over a lightlike plane and $\theta_\mb=m \pi,$ $m \in \n.$ Therefore $\mb_\infty$ is the lightlike plane containing the lightlike straight line $L=\widetilde{\alpha}^0_{-\infty} \cup \widetilde{\alpha}^0_{+\infty}.$ Furthermore,  Lemma \ref{lem:domlim}  gives that $\widetilde{\alpha}_{-\infty}^0 \subset {\cal C}_0^-$ and  $\widetilde{\alpha}_{+\infty}^0 \subset {\cal C}_0^+,$ hence $m$ is even. Finally, observe that any straight line in $\{t=0\}$ not contained in $\mb_\infty$ meet $\mb$ into a compact set, hence $\mb$ is $\omega^*$-maximal and $(i)$ holds. 

In the sequel, and up to a Lorentzian isometry preserving our normalizations, we will suppose that $\theta^+=+\infty.$

Let $\alpha \cong [0,1[ \subset \l^3-\mb$ be an affinely spacelike arc. Since $\pi(\alpha)$ is a proper curve contained in a half space of $\Pi_0,$ any non compact lifting to $\mb$ of $\pi(\alpha)$  lies in $\mb^J$ for a suitable finite interval $J \subset \Theta_0(\mb).$ In particular, $\pi(\alpha)$ has infinitely many   non compact liftings to $\mb,$ namely $\{\beta_k \,:\, k \in F\subset \z\}.$ We are assuming that this family of curves has been  {\em ordered by heights}, that is to say,  $F= \z \cap ]r^-,+\infty[$ where  $r^- \in [-\infty,+\infty[,$ and $\beta_{k_1}$ lies above $\beta_{k_2}$ outside a compact set provided that $k_1 >k_2.$  For any $k \in F,$ let $s_k=\max \{s \in \r \,:\, \beta_k \subset \mb^{[\theta(s),+\infty[}\}$ and define $I_k=[\theta(s_k),+\infty[.$  Obviously,   $k_2 >k_1$ implies that $s_{k_2}>s_{k_1},$ $\beta_{k_2} \subset \mb^{I_{k_1}}$ and $\beta_{k_1} \cap \mb^{I_{k_2}}$ is compact. Furthermore,   $\lim_{k \to +\infty} s_k =+\infty.$

\begin{quote}
{\bf Claim 1:} {\em There is $k_0 \in F$ such that $\beta_k$ lies above $\alpha \cap \pi^{-1}(\pi(\beta_k))$ for any $k > k_0.$ }
\end{quote}
\begin{proof} Up to relabeling assume that $\n \subset F.$ Reason by contradiction and suppose that $\beta_k$ lies below $\alpha \cap \pi^{-1}(\pi(\beta_k))$ for any $k \in F.$ For any $k\in F,$ the initial point $p_k$ of $\beta_k $ lies in either $\Gamma \cap \pi^{-1}(\alpha)$ or $\pi^{-1}(\pi(p)),$ where $p$ is the initial point of $\alpha.$ However, the properness of $\mb$ gives that  $\{p_k \,:\, k \in \n\} \cap \pi^{-1}(\pi(p))$ contains finitely many points below $p,$ hence we can suppose without loss of generality that $p_k \in \Gamma$ for any $k\in \n.$ Therefore $p_k\in \Gamma^{I_k},$ and so $\Gamma^{I_k}$ contains points below $\alpha$ for any $k \in \n.$ By equation (\ref{eq:normal}) we have that $\Gamma^{I_k} \subset \mbox{Int}({\cal C}^+_{\Gamma(s_k)}),$ hence $\alpha \cap  \mbox{Int}({\cal C}^+_{\Gamma(s_k)}) \neq \emptyset$ for any $k \in \n,$ contradicting that $\alpha$ is affinely spacelike. \end{proof}

For the remainder and up to relabeling we will suppose that $k_0=0.$

Let $W$ be an spacelike wedge containing $\alpha,$ and call $Z$ as its axis (i.e., the intersection of its boundary faces). Up to a Lorentzian isometry preserving our previous normalizations, we will suppose that $Z \subset \Pi_0,$ and label $Z_0 \subset \Pi_0$  as  the straight line orthogonal to $Z$ passing through the origin. Furthermore, write $Z_0=\{s e^{i \xi_0}\,:\, s \in \r\},$ $\xi_0 \in [0,2 \pi[$ and suppose that $Z_0^+=\{s e^{i \xi_0}\,:\, s\geq 0\}$ and $W$ meet at a compact set. A compact interval $I \subset \Theta_0(\mb),$ $I$  is said to be {\em centered} if its middle point lies in  $\{\xi_0+2 m \pi\,:\, m \in \z\}$ and $I \subset I_1=[\theta(s_1),+\infty[.$ 

\begin{quote}
{\bf Claim 2:} {\em If $I=[\xi-\rho,\xi+\rho]$  is a centered  and $ \rho>2\pi$ then  $\tau_0(\mb^I)\geq \Xi>0,$ where $\Xi$ is the constant given in Corollary \ref{co:universal}.}
\end{quote}
\begin{proof}
Observe that $\beta_{k} \cup \beta_{k+1}\subset  \mb^I,$ for some $k \in \n,$ and take a maximal graph $G_0 \subset \mb^I$ with boundary $\hat{\beta}_{k} \cup \hat{\beta}_{k+1} \cup \beta,$ where $\hat{\beta}_k \subset \beta_k,$ $\hat{\beta}_{k+1} \subset \beta_{k+1}$ are non compact proper subarcs and  $\beta \subset \mb^I$ is a compact arc. Since $W \subset U_\delta$ for suitable $\delta \in ]0,\frac{\pi}{4}[,$ $\partial(G_0) \subset U_\delta$ up to a compact subset. The claim follows from Corollary \ref{co:universal}. \end{proof}

\begin{quote}
{\bf Claim 3:} {\em  $\mb_\infty^+$  is a non timelike plane.}
\end{quote}
\begin{proof} Any domain in $\mb_\infty^+-\{O\}$  is the limit in the ${\cal C}^1$-topology  of a sequence of maximal graphs. Therefore,
the singular segments of $u^+$ in $\Pi_0-\{O\}$ are either  complete lines or half lines starting from the origin  (see Theorem \ref{th:converge} and Remark \ref{re:strip}). We will call  ${\cal A}^+$  as the singular set of $u^+$ in $\Pi_0-\{O\}.$

Taking into account Remark \ref{re:strip}, we only have four possibilities: $(a)$ ${\cal A}^+=\emptyset,$ $(b)$ ${\cal A}^+$ contains a complete straight line as singular segment, $(c)$ any half line in $\Pi_0-\{O\}$ with endpoint the origin is a singular segment of $u^+,$  or $(d)$ there is a wedge $W \subset \Pi_0-\{O\}$ bounded by two singular half lines of $u^+$ with endpoint $O$ and such that $u^+|_{W-\partial(W)}$ is smooth and defines a  maximal graph $G.$

In  case $(a),$  it is well known that $\mb_\infty^+$ is either a spacelike plane or a half of the Lorentzian catenoid, see \cite{ecker}.
In case $(b),$ Lemma \ref{lem:basico} implies that $\mb_\infty^+$ is a lightlike plane.

Let us see that $(c)$ is impossible. Indeed, in this case $\mb_\infty^+=\cup_{p \in \mb_\infty^+-\{O\}}  L_p,$ where $L_p$ is a lightlike half line containing  $p$ and $O,$ hence  $\mb_\infty^+ \subset {\cal C}_0.$ Furthermore, by a connectedness argument either $\mb_\infty^+ \subset {\cal C}_0^+$ or $\mb_\infty^+ \subset {\cal C}_0^-.$ Let $I$ be a centered interval of length $>4 \pi$ and take a real number $\epsilon\in]0,\Xi[,$ where $\Xi$ is the constant of Corollary \ref{co:universal}. From the definition of $\mb_\infty^+,$  there exists large enough $k_0\in \n$  such that $\mb_\infty^{2k_0\pi+I} \cap \{(x,t) \;:\; 1\leq \|x\|_0\leq 2\}$ lies in the Euclidean neighborhood of radius $\epsilon/2$ of ${\cal C}_0\cap \{(x,t) \;:\; 1\leq \|x\|_0\leq 2\}.$ On the other hand, and from the definition of $\mb_\infty^{2k_0\pi+I},$  we can find $j_0\in \n$ such that $\mb_{j}^{2k_0\pi+I} \cap \{(x,t) \;:\; 1\leq \|x\|_0\leq 2\}$ lies in the Euclidean neighborhood of radius $\epsilon/2$ of $\mb_\infty^{2k_0\pi+I} \cap \{(x,t) \;:\; 1\leq \|x\|_0\leq 2\},$ for any $j \geq j_0.$ Thus $\mb_{j}^{2k_0\pi+I} \cap \{(x,t) \;:\; 1\leq \|x\|_0\leq 2\}$ lies in the Euclidean neighborhood of radius $\epsilon$ of ${\cal C}_0\cap \{(x,t) \;:\; 1\leq \|x\|_0\leq 2\},$ for any $j \geq j_0,$ and proves that $\tau_0(\mb^{2k_0\pi+I})\leq \epsilon<\Xi.$ This contradicts Claim 2 and proves that $(c)$ is impossible. 

Suppose $(d)$ holds and label as $l_j$ as the two lightlike half lines in $\partial(G).$ From Corollary \ref{co:wedge}, $l_1$ and $l_2$ must lie in the same lightlike straight line $l,$ and since $l_1 \cup l_2 \neq l$ (otherwise  from $(b)$ $\mb_\infty^+$ would be a lightlike plane, impossible), we infer that $l_1=l_2$ and $\mb_\infty^+$ is congruent in the Lorentzian sense to the  Enneper graph $E_2$ (see Proposition \ref{pro:unienn}).

Summarizing,  in order to prove the  claim it suffices to check that $\mb_\infty^+$ can not be neither a half of the Lorentzian catenoid  nor an Enneper's graph $E_2.$ Reason by contradiction, and observe that in both cases $\mb_\infty^+$  is asymptotic at the origin to either $ {\cal C}_0^+$ or ${\cal C}_0^-$(see Remark \ref{re:cono}).

As above take a centered interval $I$ of length $>4 \pi$ and a real number $\epsilon \in ]0,\Xi[.$ Up to a dilation assume that  $\mb_\infty^+\cap \{(x,t) \;:\; 1\leq \|x\|_0\leq 2\}$ lies in the Euclidean  neighborhood of radius $\epsilon/3$ of ${\cal C}_0 \cap \{(x,t) \;:\; 1\leq \|x\|_0\leq 2\}.$ Reasoning as above there are $k_0,$ $j_0 \in \n$ such that $\mb_{\infty}^{2 k_0 \pi+I}\cap \{(x,t) \;:\; 1\leq \|x\|_0\leq 2\}$ and $\mb_j^{2 k_0 \pi+I} \cap  \{(x,t) \;:\; 1\leq \|x\|_0\leq 2\}$ 
 lie in the Euclidean neighborhood of radius $\epsilon/3$ of $\mb_\infty^+ \cap \{(x,t) \;:\; 1\leq \|x\|_0\leq 2\}$ and $\mb_{\infty}^{2 k_0 \pi+I}\cap \{(x,t) \;:\; 1\leq \|x\|_0\leq 2\},$ respectively, for any $j \geq j_0.$ 
As above, this shows that  $\tau^+_0(\mb^{2 k_0 \pi +I} )\leq \epsilon<\Xi,$ contradicting Claim 2. \end{proof}

Now we can prove $(ii).$  From Claim 3  $\mb_\infty^+$ is a timelike plane, and by a symmetric argument the same holds for $\mb_\infty^-.$ If the planes $\mb_\infty^+$ and $\mb_\infty^-$ were different they would meet transversally (both planes contain the origin!), which would contradict that $\mb_\infty^+$ and $\mb_\infty^-$ lie in the limit set of a sequence of embedded surfaces.

Finally let us see $(iii).$ We have to deal with the case $\theta^-\in \r,\, \theta^+=+\infty.$ 

Let us show first that $\mb_\infty^+$ is a lightlike plane. Reason by contradiction and assume that it is spacelike.  In particular, $\mb_\infty^+$  intersects transversally  any lightlike  plane or  maximal multigraph of finite angle with lightlike limit tangent plane at infinity. By  Lemma \ref{pro:singu}, we deduce that  $\arg(\overline{\cal W}-\{[0]\})-I_{\cal X}$ contains a non compact component $J=]a,+\infty[,$ where $a \in I_X^c.$ 

Take an  arc $c \subset \mb_\infty^J\cup\{O\}$  projecting onto a divergent arc $\pi(c)$ in $\Pi_0$ with initial point $O$ and satisfying $lim_{x \in c \to \infty}\mbox{d}(x,\mb_\infty^+)=0.$ Let us see  that  $\lim_{x \in c\to \infty} g(x)=\mbox{st}(v),$ where $v\in \h^2_-$ is the Lorentzian normal to $\mb_\infty^+.$ Indeed, take a divergent sequence $\{p_n, \; n \in \n\}\subset c,$
call $J_n=[\arg(p_n)-\frac{\pi}{2},\arg(p_n)+\frac{\pi}{2}]$ and consider the sequence of translated multigraphs $\{G_n:=-p_n+\mb_\infty^{J_n}, \;n \in \n\}.$ By Proposition \ref{pro:converge}, $\{G_n\}_{n \in \n}$  converges in the ${\cal C}^1$-topology to an entire PS graph $G_\infty$ over $\Pi_0$ lying in a closed half space bounded by the spacelike plane $\mb_\infty^+.$  Taking into account Theorem \ref{th:calabi}, we deduce that   $G_\infty$  is a spacelike plane, i.e., $G_\infty=\mb_\infty^+,$ proving our assertion.

On the other hand, set $l$ the lightlike half line $\partial(\mb_\infty^J)={\cal X}(\arg^{-1}(a) \cup \{[0]\})$  and denote by $N_0\subset \overline{{\cal W}^J}$ the proper region  bounded by $c\cup l.$ Consider  a new proper region $N_0' \subset N_0$ having $\partial(\N_0')=c \cup l',$ where $l' \subset \mbox{Int}(N_0)$ is a divergent arc close enough to $l$ in such a way that $\lim_{x \in l' \to \infty} g(x)=\mbox{st}_0(w),$ where $w$ is the lightlike direction of $l$ and $g$ is the holomorphic Gauss map of $N_0'.$ Since $N'_0$ lies in a half space bounded by $\mb_\infty^+$ it is parabolic (see Corollary \ref{co:halfspace}), and Theorem \ref{th:sectorial} implies that $\mbox{st}(v) =\lim_{x \in c \to \infty} g(x)=\lim_{x \in l \to \infty} g(x)=\mbox{st}_0(w).$ This is obviously absurd (note that $\mbox{st}(v) \in \d$  and $\mbox{st}_0(w)\in \partial(\d)$) and shows that the plane $\mb_\infty^+$ must be lightlike.

As a consequence, $J_0=\mbox{Int}(I_{\cal X}^c)$ is either empty or a finite interval with endpoint $\theta^-.$ Indeed, otherwise  Proposition \ref{pro:singu}, $(iii)$ yields that  $J=]b,+\infty[,$ and so $\mb_\infty^{J_0} \subset {\cal C}_0,$ contradicting that $\mb_\infty^+$ is a lightlike plane.

Let $L$ denote the lightlike straight line in $\mb_\infty^+$ passing through $O,$ and call $L^-=L \cap {\cal C}_0^-.$ From  equation (\ref{eq:monoto}), $L^-$ lies above any arc in $\pi^{-1}(\pi(L^-))\cap \mb_\infty,$ and taking into account equation (\ref{eq:acausalinfi}), we deduce that $L^- \subset \mb_\infty^I$ for any closed interval $I \subset \arg({\cal W}-[0])$ of length $2 \pi,$ that is to say, $I \cap I_{\cal X} \neq \emptyset$ for any compact interval of length $2\pi.$ Define $B=\{\xi \in I_1 \,:\,{\cal X}(\arg^{-1}(\xi))=L^-\}$ and take $\theta_0 \in B$ such that $J\cap \mbox{Int}(I_{\cal X}^c)=\emptyset,$ where $J=[\theta_0,+\infty[.$ 

\begin{quote}
{\bf Claim 4:} {\em  $J - I_{\cal X}$ is either empty or bounded.}
\end{quote}
\begin{proof}
 Reason by contradiction, let us see that any connected component of $J - I_{\cal X}$ determines an Enneper's graph with limit tangent plane at infinity parallel to $\mb_\infty^+.$ Consider  a component $I$ of $J - I_{\cal X}$ (hence of $\arg({\cal W}-[0])-I_{\cal X}$), and observe that $J \cap \mbox{Int}(I_{\cal X}^c)=\emptyset.$ Our previous analysis implies that $|I|=2\pi,$ the endpoints of $I$ lie in $B$ and $\mb_\infty^I$ is a maximal graph. By  Proposition \ref{pro:unienn},  $\mb_\infty^I$ is an Enneper graph  with an upward  conelike singularity at the origin (i.e., asymptotic at the origin to  ${\cal C}_0^-$). Furthermore, since $\mb_\infty^I$ lies below $\mb_\infty^+$ then its limit tangent plane at infinity is parallel to $\mb_\infty^+,$ proving our assertion.

Let $T \subset J$ be a compact centered interval of length $>4\pi$ and take a real number $\epsilon\in]0,\Xi[,$ where $\Xi$ is the constant of Corollary \ref{co:universal}. Up to a dilation, assume that $\mb_\infty^T \cap \{(x,t) \;:\; 1\leq \|x\|_0\leq 2\}$ lies in a neighborhood of ${\cal C}_0^-$ of radius $\epsilon/2.$ Then take $j_0\in \n$  large enough such that $\mb_j^{T} \cap \{(x,t) \;:\; 1\leq \|x\|_0\leq 2\}$ lies in the Euclidean neighborhood of radius $\epsilon/2$ of $\mb_\infty^T\cap \{(x,t) \;:\; 1\leq \|x\|_0\leq 2\},$ for any $j \geq j_0.$ We infer that $\mb_{j}^{T} \cap \{(x,t) \;:\; 1\leq \|x\|_0\leq 2\}$ lies in the Euclidean neighborhood of radius $\epsilon$ of ${\cal C}^-_0\cap \{(x,t) \;:\; 1\leq \|x\|_0\leq 2\}$ for any $j \geq j_0,$ proving that $\tau_0(\mb^{T})\leq \epsilon<\Xi.$ This is contrary to Claim 2, proving the claim. \end{proof}

Claim 4 shows that $I_{\cal X}-\mbox{Int}(I_{\cal X}^c)$ contains an interval $J_{\cal X}$ of the form $[b,+\infty[,$ $b \in \r.$ Item $(iii)$ is an  elementary consequence of Proposition \ref{pro:singu}, $(i).$ \end{proof}


\begin{definition} \label{def:limitplane}
In the context of Theorem \ref{th:strong}, if $\Theta_0(\mb)=\r$ or $\theta_\mb<+\infty$ then  $\Sigma_\infty:=\mb_\infty$  is defined to be the blow-down plane of $\mb$ associated to $\{\lambda_n\}_{n \in \n}.$  Likewise, if $\Theta_0(\mb)\neq \r$ and $\theta_\mb=+\infty$ we set  $\Sigma_\infty:=\mb_\infty^{J_{\cal X}}$ and call it the blow-down plane of $\mb$ associated to $\{\lambda_n\}_{n \in \n}$ as well.    
\end{definition}

\subsection{The transversality of $\mb$ and the blow-down plane $\Sigma_\infty$}

This subsection is devoted to prove that  the Lorentzian Gauss map of $\mb$ omits the normal direction to the blow-down plane.  We need some notations and preliminary results.

Let  $c$ be a lightlike ray in $\mb,$ call $l_c$ the lightlike half line to which $c$ is asymptotic  and write $L_c:=\pi (c)=\pi(l_c).$ Putting $L_c=\{x_0+a e^{\xi}\;:\; a\geq 0\},$  there is a unique real number $\xi_c$ congruent to $\xi$ modulo $2 \pi$ such that $c \subset \mb^J,$ where  $J= [\xi_c -\frac{\pi}{2},\xi_c +\frac{\pi}{2}] \cap ]\theta^+,\theta^-[.$

As a consequence, the limit  $\theta_c:=\lim_{x \in c \to \infty} \Theta_0(x) \in [\theta^-,\theta^+]$ exists and is a finite real number.
The arguments $\theta_c$ and $\xi_c$ coincide provided that $\theta_c \in ]\theta^+,\theta^-[.$ If  $\theta_c \in \{\theta^+,\theta^-\},$ then $J$ is a tail interval and either $\xi_c \in [\theta^+,\theta^++\frac{\pi}{2}]$ or  $\xi_c \in [\theta^--\frac{\pi}{2},\theta^-]$ (see Lemma \ref{lem:sheet}).

\begin{lemma} \label{lem:importante}
If $\mb$ admits an upward (resp., downward)   lightlike ray $c,$ then  $\theta^+=\theta_c+\frac{ 3\pi}{2}$  (resp., $\theta^-=\theta_c-\frac{ 3\pi}{2}$) and $\theta_\mb=+\infty.$
\end{lemma}
\begin{proof}
We only deal with the case when $c\subset \mb$ is an upward lightlike ray.

\begin{quote}
{\bf Claim 1:} {\em $\theta_c+\frac{3\pi}{2} \geq \theta^+.$}
\end{quote}
\begin{proof}
Reason by contradiction, and assume that $\theta_c+\frac{3\pi}{2} < \theta^+.$ Let $s_c\in \r$ denote the unique real number such that $\theta(s_c)=\theta_c,$ and let us show that $\widetilde{\alpha}_{s_c}=F_0(s_c,\cdot)$ is an upward lightlike ray too. Indeed, since $\alpha_{s_c}$ and $L_c$  are  parallel half lines, the spacelike condition gives that $\mbox{d}_H(\widetilde{\alpha}_{s_c},c)\leq \sqrt{2} \mbox{d} (\alpha_{s_c},L_c).$ Taking into account that ${\widetilde{\alpha}}_{s_c}$ has slope $<1$   we deduce that the limit $\lim_{x \in {\widetilde{\alpha}}_{s_c} \to \infty} \mbox{d}(x,{\cal C}_{\Gamma(s_c)}^+)$ exists and is finite, proving the assertion.

Write $I_c:=[\theta_c,\theta^+[$ and let  $H\subset \l^3$ be an open half space containing  $\partial(\mb^{I_c})={\widetilde{\alpha}}_{s_c}\cup \Gamma^{I_c}$ and  such that $\partial(H)$ is a lightlike plane parallel to $l_c$ (here we have taken into account  equation (\ref{eq:normal})).

Let us see that $\mb^{I_c}\subset H.$ To do this, label $L_0$ as the complete straight line in $\Pi_0$ containing  $\alpha_{s_c}$ and call  $l_0=\partial(H)\cap \pi^{-1}(L_0).$ First, we observe that $\mb^{I_c}$ is disjoint from $l_0.$ Indeed, if $c'$ a connected component of $\pi^{-1}(L_0)\cap \mb^{I_c},$ then either $c'$ has en endpoint  in $\partial(\mb^{I_c}) \subset H$ or $\pi(c') \cap \alpha_{s_c} \neq \emptyset,$ and in the second case  $c'\cap \pi^{-1}(\alpha_{s_c})$ lies above $\widetilde{\alpha}_{s_c}\cap \pi^{-1}(c').$ By property  (\ref{eq:acausal}) we deduce that $c'\subset H,$ hence $\mb^{I_c}\cap l_0=\emptyset.$

Reason by contradiction and assume that $\mb^{I_c} -H\neq \emptyset.$ Then take a connected
component $S$ of $\mb^{I_c} -H.$ By the spacelike condition, there are no compact arcs in $S$ with endpoints in $\partial(S) \cap \partial(H)$ projecting onto a segment parallel to $L_0.$ Since  $S \cap
\partial(\mb^{I_c})=\emptyset$ then we deduce that  $S$ is simply connected. On the other hand,  $\pi|_{S}:S \to
\Pi_0$ is a proper local embedding, hence $S$ is a graph over
$\Pi_0.$ Taking into account that $S \cap l_0=\emptyset$ and
$\partial(S) \subset \partial(H),$  $G:=S \cup
(\partial (H)-\pi^{-1}(\pi(S)))$ is a PS entire graph over
$\Pi_0$ containing $l_0.$ Lemma \ref{lem:basico} gives that
$G=\partial(H),$ which is absurd and  proves that
$\mb^{I_c}\subset H.$

Let $s_0$ be the unique real number satisfying $\Theta_0(\Gamma(s_0))=\theta_c+\frac{3 \pi}{2},$ and call $T_{\Gamma(s_0)} \mb$ as the tangent plane to $\mb$ at $\Gamma(s_0)$ (that is to say, the lightlike plane parallel to the vector $\Gamma'(s_0)$). Observe that $T_{\Gamma(s_0)} \mb$ and $\partial (H)$ are parallel, and call $H'$ as the open half space containing $\partial(H)$ and with boundary $T_{\Gamma(s_0)} \mb.$  Let $y:\mb \to \r$ denote the harmonic coordinate function  vanishing on $\partial(H')\cap \mb.$ From equation (\ref{eq:wei}), the holomorphic 1-form $dy$ (that can be reflected holomorphically to the mirror surface $\mb^*$), has a zero or order $\geq 2$ at $\Gamma(s_0).$ Therefore, $\mb^{I_c}-\partial(H')$ has at least a connected component $S$ lying in the slab $H'\cap H$  and with boundary $\partial(S) \subset \partial(H').$ As before, $S$ is  a graph over $\Pi_0$ and  $S_0:=S \cup(\pi|_{\partial(H')})^{-1}(\Pi_0-\pi(S))$ is an entire PS  graph over $\Pi_0.$ Take $p_0 \in S-\partial(S)$ and a neighborhood $D_0 \subset  S-\partial(S)$ of $p_0$ projecting via $\pi$ onto a closed disc. From equation (\ref{eq:acausal}), $D_0-\{p_0\} \subset \mbox{Ext}({\cal C}_{p_0})$ and  $\delta:= \mbox{d}(\partial(D_0),{\cal C}_{p_0})$ is  positive. As a consequence, the PS graph $S_0-D_0$ lies in $\mbox{Ext}({\cal C}_{p_0})$  and
$\mbox{d}(S_0-D_0,{\cal C}_{p_0})\geq \delta>0.$ This shows that $S_0$ lies in $\{x \in \l^3 \;:\; \|x\| \geq \delta\}$ up to a compact subset, and the same holds for $S.$  From Theorem \ref{th:parabo}, $S$ is a parabolic. But $\partial(S) \subset \partial(H')$ gives $S \subset \partial(H'),$ which is absurd and proves the claim. \end{proof}

Claim 1 gives that  $\theta^+ \in \r,$ and so  $\mb^{I_c}$ is a multigraph of finite angle bounded by  two sublinear arcs, namely $\Gamma^{I_ c}$ and $\widetilde{\alpha}_{s_c},$  with lightlike direction.  By Corollary \ref{co:wedge}, these arcs have the same direction, that is to say, $\theta^+=\theta_c+\frac{3\pi}{2},$ concluding the first part of the lemma.

It remains to check that $\theta^-=-\infty.$ Reason by contradiction and assume $\theta^- \in \r.$ As above, Corollary \ref{co:wedge} gives that $\mb$ is parabolic, hence $\mb$ is biholomorphic to $\overline{\u}=\{z \in \c \;:\; \mbox{Im}(z) \geq 0\}.$ Let $X:\overline{\u} \to \l^3$ be a conformal maximal embedding satisfying $X(\overline{\u})=\mb.$ Set $(\phi_3,g)$  the Weierstrass representation of $X,$ see equation (\ref{eq:wei}). The holomorphic map $g$ extends by Schwarz reflection to a meromorphic map on $\c$ of finite degree $n,$ and so  we can put $g(z)=\frac{P(z)}{Q(z)},$ where $P(z)=\sum_{j=0}^n a_j z^j,$ $Q(z)=\sum_{j=0}^n \overline{a_j} z^{j}$ and $a_n \neq 0.$ Since the 1-forms $\phi_1,$ $\phi_2$ and $\phi_3$ have no common zeroes in $\overline{\u},$ we get $\phi_3=-i B P(z) Q(z)dz, \quad B\in ]0,+\infty[.$ Up to a Lorentzian isometry, we can suppose that $g(\infty)=1,$ $a_n=1$ and $\theta^+=\lim_{r \to +\infty} g(r)$ (note that $X^{-1}(\Gamma)=\r$).  Therefore we also have $\frac{\theta^+}{2 \pi} \in \z$ and  $\theta^-=\theta^+-2 n \pi.$

Observe that $f_2(z):=\int_0^z \phi_2= \frac{B}{2}\int_0^z (P(w)^2+Q(w)^2) dw,$ $f_1(z):=\int_0^z \phi_1= \frac{i B}{2}\int_0^z(P^2(w)-Q^2(w))dw$ and
$f(z):=\int_0^z (\phi_2-i\phi_3)= \frac{B}{2}\int_0^z(P(z)-Q(z))^2$ are polynomial functions of degrees $2n+1,$ $n+n_0+1$ and $2n_0+1,$ respectively, where $0\leq n_0=\mbox{Deg}(P(z)-Q(z)) <n.$ Since $\widetilde{\alpha}_{s_c}$ is a lightlike ray with direction $(0,1,1),$ then the limits $\lim_{z \in X^{-1}(\widetilde{\alpha}_{s_c})\to \infty} \mbox{Re}(f_1(z))$ and $\lim_{z \in X^{-1}(\widetilde{\alpha}_{s_c})\to \infty} \mbox{Re}(f(z))$ are finite,  and so there are positive odd integers $m_1$ and $m$ such that $\lim_{z \in X^{-1}(\widetilde{\alpha}_{s_c}) \to \infty}\arg(z)=\frac{m_1 \pi}{2(n+ n_0+1)}=\frac{m \pi}{2(2n_0+1)}.$
Taking into account that $X^{-1}(\partial(\mb^{I_c}))=[r,+\infty[ \cup \widetilde{\alpha}_{s_c},$ where $r=X^{-1}(\Gamma(s_c)),$ and that $\mb^{I_c}\subset \{(x,y,t) \;:\; t-y \geq R\}$ for a suitable $R \in \r,$ we infer that $m=1,$ hence $n-n_0=(2n_0+1) (m_1-1) \geq 2(2 n_0+1)$ and $n \geq 5 n_0+2.$

On the other hand, note that for any $k \in \r,$ the set $\{z \in \overline{\u} \,:\, X(z) \in \mb^{I_c} \,\mbox{and}\, \mbox{Re}(f_2(z))=k\}$  consists of either a proper arc $\cong ]0,1[$ or two proper arcs homeomorphic to $[0,1[.$ Indeed, just take into account that $\mb^{I_c}$ is (up to removing a compact set) a multigraph of angle $2\pi$ with sublinear boundary arcs $\Gamma^{I_c}$ and $\widetilde{\alpha}_{s_c}$ of direction $(0,1,1).$ As any divergent nodal arc in $\overline{\u}$ of the harmonic function $\mbox{Re}(f_2)-k$  is asymptotic to $\{s e^{\frac{j \pi}{2 (2n+1)}},\;s \geq 0\}$ for a suitable odd integer $j\leq n,$ and $\mb^{I_c}$ contains only two such arcs for any $k \in \r,$  then  $\frac{5 \pi}{2 (2n+1)}>\lim_{z \in X^{-1}(\widetilde{\alpha}_{s_c}) \to \infty}\mbox{Im}(\log(z))=\frac{\pi}{2(2n_0+1)},$ or equivalently $n<5 n_0+2,$ which is absurd and concludes the proof. \end{proof}

Set $t_{\mb}:=\pi^{-1}(O)\cap \mb$ the intersection of $\mb$ and the $t$-axis, and for any $q=F_0(s,a) \in t_{\mb}$ write $r_q=F_0(\{s\} \times [0,a]).$ Consider the simply connected surface ${\cal S}:=\mb-\cup_{q \in t_{\mb}} r_q$ and fix a branch $\textsf{f}$ of $\log \circ \pi$ along ${\cal S}.$ It is clear that $\kappa:{\cal S} \to \mb \times \overline{\c},$ $\kappa(q)=\left(q,\textsf{f}(q)\right)$ is an embedding, and that $\hat{\mb}:=\overline{\kappa({\cal S})}\subset \mb \times \overline{\c}$ is a surface with piecewise analytical boundary  homeomorphic to $\overline{\d}-\{1\}.$ 
Set $Y:\hat{\mb} \to \mb$ the projection map $Y(q,\textsf{f}(q))=q,$  and  note that for any $q \in t_{\mb}$ we have $Y^{-1}(r_q)=r_q^+\cup r_q^-,$ where $r_q^+\cap r_q^- =Y^{-1}(q)=(q,\infty)$ and $Y|_{r_q^+}:r_q^+ \to r_q,$  $Y|_{r_q^-}:r_q^- \to r_q$ are homeomorphisms.
Furthermore, the restriction of $Y$  to $ \hat{\mb}-\cup_{q \in t_{\mb}} Y^{-1}(r_q-\{q\})$ is one to one, and $\partial(\hat{\mb})=Y^{-1}\left(\Gamma \cup \left(\cup_{q \in t_{\mb}} r_q\right)\right).$

\begin{figure}[htpb]
\begin{center}\includegraphics[width=.7\textwidth]{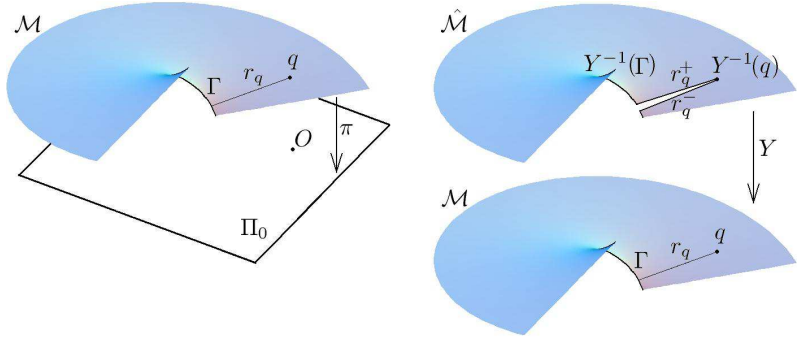}
\caption{The surface $\hat{\mb}$ and the projection $Y.$}\label{fig:segunda}
\end{center}
\end{figure}

As $\hat{\mb}-Y^{-1}(t_{\mb})$ is simply connected and $\textsf{f}|_{\hat{\mb}-Y^{-1}(t_{\mb})}$ never vanishes,  $\textsf{Arg}:=\mbox{Im}(\textsf{f})$ is well defined on $\hat{\mb}-Y^{-1}(t_{\mb}).$

\begin{remark} \label{re:dosargu} Up to a suitable choice of the branch,  $\lim_{x \in \mb^J \to \infty} \textsf{Arg}(Y^{-1}(x))-\Theta_0(x)=0,$ provided that $J \subset \Theta_0(\mb)$ is compact. 

Therefore,  $\theta^--\frac{\pi}{2} \leq \liminf_{x \in \mb \to \infty} \textsf{Arg}(Y^{-1}(x)),$ $\limsup_{x \in \mb \to \infty} \textsf{Arg}(Y^{-1}(x)) \leq \theta^++\frac{\pi}{2},$ and
 if $\{x_n\}_{n \in \n}\subset \mb$ is a divergent sequence satisfying that $ \liminf_{n \to \infty} \textsf{Arg}(Y^{-1}(x_n)) \subset [\theta^+,\theta^++\frac{\pi}{2}]$ (resp., $ \limsup_{n \to \infty} \textsf{Arg}(Y^{-1}(x_n)) \subset[\theta^--\frac{\pi}{2},\theta^-]$), then $\lim_{n \to \infty} \Theta_0(x_n)=\theta^+$ (resp., $\theta^-$).

Moreover, standard  monodromy  arguments also show that:

\begin{enumerate}[(i)]
\item $]\theta^--\frac{\pi}{2},\theta^++\frac{\pi}{2}[ \subset \textsf{Arg}(\partial(\hat{\mb}))$ and $\textsf{Arg}^{-1} (J) \cap \partial(\hat{\mb})$  is compact, for any compact interval $J \subset \r-\{\theta^--\frac{\pi}{2},\theta^++\frac{\pi}{2}\}.$
\item If $\{q_n\}_{n \in \n} \subset \partial(\hat{\mb})$ and $\{\Theta_0\left(Y(q_n)\right)\}_{n \in \n} \to \theta^+$ (resp., $\theta^-$), then $\{\textsf{Arg}(q_n)\}_{n \in \n} \to \theta^++\frac{\pi}{2}$ (resp., $\theta^--\frac{\pi}{2}$).
\end{enumerate}
\end{remark}

Denote by ${\cal V}_r:=\overline{\mbox{Int}({\cal C}^+_{(0,0,r)})}$ and $\Gamma_r:=\Gamma \cap {\cal V}_r.$ Since $(1,0,0) \in \Gamma,$ equation (\ref{eq:normal}) gives that $\Gamma_r\neq \emptyset$ provided that $r \leq -1.$ In this case, $\Gamma \cap \partial({\cal V}_r)$ consists of a unique point $\Gamma(s_r)$  and $\Gamma_r=\Gamma^{I_r},$ where $I_r=[\theta(s_r),\theta^+[.$
In the sequel we will  suppose $r \leq -1.$

We label ${\mb}(r)$ as  the connected component of ${\cal V}_r \cap \mb$ containing  $\Gamma_r,$ and write ${\mb}'(r)=({\cal V}_r \cap \mb)-\mb (r).$  Likewise we put  $\hat{\mb} (r):=Y^{-1}({\mb} (r)).$

It is interesting to observe that $\mb=\cup_{n \in \n} \mb(r_n),$ provided that $\{r_n\}_{n \in \n} \subset ]-\infty,-1]$ is divergent. Indeed, 
fix an arbitrary point $q = F_0(s,a) \in \mb$  and take $n\in \n$ large enough  in such a way that $\Gamma(s)=F_0(s,0)$ and $q$ lie in ${\cal V}_{r_n}.$ By equation (\ref{eq:acausal}), $F_0(\{s\}\times [0,a])$ is contained in ${\cal V}_{r_n},$ and so $q \in F_0(\{s\} \times [0,a])\subset \mb(r_n).$ Therefore $q \in \cup_{n \in \n} \mb(r_n)$ and we are done.

\begin{lemma} \label{lem:argmon}
If $c\subset \partial({\cal V}_r)-\{(0,0,r)\}$ is a spacelike arc and $\rho$ is a branch of $\mbox{Im}\left(\log \circ \pi \right)|_c,$ then $\rho$ is monotone and $\|q_2-q_1\|>0$ provided that $q_1,$ $q_2 \in c$ and $0<|\rho(q_2)-\rho(q_1)|<2\pi.$
\end{lemma}
\begin{proof} The spacelike property gives that $\rho$  has neither local maxima  nor minima, hence it is monotone.
For any $q \in \partial({\cal V}_r)-\{(0,0,r)\},$ label $l_q$ as the closed lightlike half line in $\partial({\cal V}_r)$ containing $q.$ It is obvious that $\partial({\cal V}_r)-l_q \subset \mbox{Ext}({\cal C}_q),$ which simply means that $\|q'-q\|>0$ for any $q' \in \partial({\cal V}_r)-l_q.$ If $q \in c,$ the  monotonicity of  $\rho$ yield that $\{q' \in c-\{q\} \;:\; |\rho(q')-\rho(q)|<2 \pi\} \subset  c-l_q \subset \mbox{Ext}({\cal C}_q).$ This concludes the proof. \end{proof}

\begin{corollary} \label{co:importante}
Assume that $\mb$  contains no upward   lightlike rays,  and fix  $r \leq -1,$ $(0,0,r) \notin \mb.$
Then the following statements hold:
\begin{enumerate}[(a)]
\item If $\hat{c}$ is an arc in $Y^{-1}(\partial({\cal V}_r)\cap \mb)$   then $\textsf{Arg}|_{\hat{c}}$ is monotone.  Furthermore, $\|Y(q_2)-Y(q_1)\|>0$ for any $q_1,$ $q_2 \in \hat{c}$ satisfying that $0<|\textsf{Arg}(q_2)-\textsf{Arg}(q_1)|<2\pi.$
\item $\partial(\mb(r))=\Gamma_r \cup \beta_r,$ where $\beta_r\cong [0,1[$ is a proper divergent arc in $\mb$ with  initial point  $\Gamma(s_r)$ and meeting $\Gamma$ only at this point.
\item If $D \subset \mb'(r)$ is a connected component, then $D $ is a closed disc, $\partial(D) \subset \partial ({\cal V}_r)$ and $t_{\mb}\cap D$ consists of a single point.
\end{enumerate}
\end{corollary}
\begin{proof} Since $\mb$ is spacelike and $\partial({\cal V}_r)$  is lightlike, they meet transversally and $\left(\mb \cap\partial ({\cal V}_r)\right)-\{(0,0,r)\}$  consists of a family of pairwise disjoint properly embedded analytical regular curves. Item $(a)$ is an elementary consequence of Lemma \ref{lem:argmon}.

From our hypothesis and Lemma \ref{lem:importante} we get $\theta^+=+\infty,$ hence from Remark \ref{re:dosargu} we have $\textsf{Arg}(\hat{\Gamma}_r)\subseteq [a,+\infty[,$ $a \in \r,$ where $\hat{\Gamma}_r$ is the smallest arc in $\partial(\hat{\mb})$ containing $Y^{-1}(\Gamma_r)$ (obvioulsy contained in $\partial(\hat{\mb}(r))$).

Let us show that $\textsf{Arg}\left(Y^{-1}({\cal V}_r \cap \mb-t_\mb)\right)= [b,+\infty[,$ $b \in \r.$
Reason by contradiction, and suppose there exists a divergent arc $\hat{c} \in Y^{-1}(\partial({\cal V}_r) \cap \mb)$ homeomorphic to $[0,1[$ such that $\textsf{Arg}(\hat{c})=]-\infty,d],$ $d \in \r,$ $d<a-\pi.$ For any $q \in \hat{c},$ let $L_q\subset \Pi_0$ denote the straight line passing through $\pi(Y(q))$ and the origin, and label $l_q$ as the connected component of  $\pi^{-1}(L_q)\cap \mb$ containing $Y(q).$ The choice of $\hat{c}$ yields that $\textsf{Arg}(Y^{-1}(l_q-t_\mb)) \cap [a,+\infty[=\emptyset,$ and so $l_q \cap {\cal V}_r$ is disjoint from $\Gamma.$ Since $\mb$ has no upward   lightlike rays, we infer that $c_q:=l_q \cap {\cal V}_r$ is a compact arc with endpoints in $\partial({\cal V}_r)-\{(0,0,r)\}$ and passing through a point of $t_{\mb}.$ However, $t_{\mb}$ is a closed discrete set, and therefore the point $c_q \cap t_{\mb}$ does not depend on $q \in \hat{c}.$ This obviously contradicts that the family of compact curves $\{c_q, \;q \in \hat{c}\}$ diverge in $\l^3$ as $q$ diverges in $\hat{c}.$

Now we can prove $(b).$ 

Indeed, first note that  $\partial(\hat{\mb} (r))$ contains no (closed) Jordan curves. To see this, recall that $\hat{\mb}$ is simply connected, and so  any such curve must bound a compact disc $V\subset \hat{\mb}.$ By the convex hull property, $Y(V) \subset {\cal V}_r,$ and since $\hat{\mb}(r)$ is a connected component of $Y^{-1}({\cal V}_r),$ then we get that $V=\hat{\mb}(r),$ which is absurd. 

Suppose there are two different divergent arcs $\hat{c}_1,$ $\hat{c}_2$ in $\partial(\hat{\mb} (r))$ homeomorphic to $[0,1[$ and disjoint from $\hat{\Gamma}_r.$  From the previous arguments, $\textsf{Arg}(\hat{c}_j)=[a_j,+\infty[,$ $j=1,2,$ hence there are points  $q_1 \in \hat{c}_1$ and $q_2 \in \hat{c}_2$ satisfying $\textsf{Arg}(q_1)=\textsf{Arg}(q_2).$  As above, set $L_j \subset \Pi_0$  and $l_j$ the straight line passing through $O$ and $\pi(Y(q_j))$ and its lifting to $\mb$ with initial condition $Y(q_j),$ respectively, $j=1,2.$
Let us check that  $l_1 \cap l_2=\emptyset.$ Indeed, the fact $L_1= L_2$ and the uniqueness of the lifting give that either $l_1=l_2$ or $l_1 \cap l_2=\emptyset,$ and the first option  leads to $Y(q_1),$ $Y(q_2) \in l_1=l_2,$ contradicting that $Y(q_2)-Y(q_1)$ is a lightlike vector. As a consequence, $\lim_{x \in l_1 \to \infty} \Theta_0(x)\neq \lim_{x \in l_2 \to \infty} \Theta_0(x).$  On the other hand, Remark \ref{re:dosargu} gives that $\textsf{Arg}(q_j)=\lim_{x \in Y^{-1}(l_j) \to \infty} \textsf{Arg}(x)=  \lim_{x \in l_j \to \infty} \Theta_0(x),$ $j=1,2,$ which contradicts that $\textsf{Arg}(q_2)=\textsf{Arg}(q_1)$ and proves $(b).$

To finish, consider a connected component $D \subset \mb'(r).$ If $\partial(D)$  contains a proper divergent arc $\alpha \cong ]0,1[,$ we can split $\alpha$ into  two divergent arcs homeomorphic to $[0,1[,$ getting a contradiction as above. Therefore, any connected component of  $\partial(D)$ is a Jordan curve. Since $D$ is simply connected, $D$ is a closed disc. Furthermore,  $\textsf{Arg}$ is monotone along $Y^{-1} (\partial(D)-\cup_{q \in t_{\mb}} r_q),$ hence  $Y(D)$  is a closed disc meeting $t_\mb$ at a unique point, proving  $(c).$ \end{proof}

\begin{lemma}  \label{lem:foliacion}

Assume that $\mb$  contains no  upward   lightlike rays,  and fix $r \leq -1,$ $(0,0,r) \notin \mb.$ Then, there exists a smooth foliation ${\cal D}(r):=\{D_s(r)\;:\; s \in [r,+\infty[\}$ of ${\cal V}_r$ satisfying:

\begin{enumerate} [(a)]
\item For any $s>r,$   $D_s(r)$ is a maximal disc, $\partial(D_s(r)) \subset \partial({\cal V}_r)$ and $(0,0,s)\in D_s(r).$ Moreover, $D_r(r)=\{(0,0,r)\}.$ 
\item $D_s(r) \cap \mb'(r) \neq \emptyset$ if and only if $D_s(r)$ is a connected component of $\mb'(r).$
\item $D_s(r) \cap \mb(r) \neq \emptyset$ if and only if $D_s(r) \cap \partial(\mb(r)) \neq \emptyset,$ and in this case $D_{s}(r) \cap \mb$ is an embedded compact arc lying in $\mb(r)$ with initial point at $\Gamma_r$ and final point at $\beta_r.$
\end{enumerate}
\end{lemma}
\begin{proof} Consider a spacelike smooth divergent embedded arc  $\delta_r \subset \partial ({\cal V}_r)$ containing $\beta_r$ and with initial point $(0,0,r).$ 

\begin{quote}
{\bf Claim:} {\em There exists a smooth foliation ${\cal F}_r:[0,+\infty[ \times \s^1 \to \partial({\cal V}_r)$ of $\partial({\cal V}_r)$ satisfying
\begin{enumerate}[(i)]
\item  $c_y:={\cal F}_r(y,\cdot):\s^1 \to \partial({\cal V}_r)$ is a Jordan curve for any $y>0,$  and $c_0$ is the constant curve $c_0(\xi)=(0,0,r),$ $\xi \in \s^1.$
\item $c_y$ and $\delta_r$ meet at a unique point in a transversal way, $y>0.$
\item  $\|c_y(\xi)-c_y(\xi')\|>0$ for any $\xi,$ $\xi' \in \s^1,$ $\xi \neq \xi',$ and any $y>0.$
\item For any connected component $D$ of $\mb'(r),$ there is an unique $y \in ]0,+\infty[$ such that  $c_y=\partial (D).$
\end{enumerate}}
\end{quote}
\begin{proof}
From Lemma \ref{lem:argmon}, any branch of   $\mbox{Im}\left(\log\circ \pi\right)|_{\delta_r}$ is monotone. Since $\pi|_{\partial ({\cal V}_r)}$ is injective, we deduce that $\pi(\delta_r)$ is an embedded divergent arc of spiral type with  initial point at $O.$ Hence, we can take a smooth foliation ${\cal F}_r^*:[0,+\infty[ \times \s^1 \to \Pi_0$ of $\Pi_0$ by Jordan curves $d_y:={\cal F}^*_r(y,\cdot):\s^1 \to \Pi_0$ (where $F_r^*(0,\cdot)$ is constant and equal to $O$) in such a way that $d_y$ bounds a starshaped domain centered at the origin (i.e., $d_y-\{O\}$ can be parameterized by the principal argument) and  $d_y\cap \pi(\delta_r)$ consists of an unique point where both curves meet transversally for any $y>0.$ Furthermore, constructing ${\cal F}_r^*$ with a little care, we can ensure that for every connected component $D$ of $\mb'(r),$ there exists an unique $y_D \in ]0,+\infty[$ such that  $\partial(\pi(D))=F_r^*(y_D,\cdot).$

It suffices to define ${\cal F}_r:=(\pi|_{\partial({\cal V}_r)})^{-1} \circ {\cal F}_r^*.$ Items $(i),$  $(ii)$ and $(iv)$ are clear, and item $(iii)$ follows from  Lemma \ref{lem:argmon}. \end{proof}

From Theorem \ref{th:plateau} and item $(iii)$, there is a unique maximal disc $K_y\subset {\cal V}_r$ with boundary $c_y,$ $y\geq 0$ (we have made the convention $K_0=\{(0,0,r)\}$). Furthermore, since $\pi|_{K_y}$ is a local homeomorphism, $K_y$ is  a graph over the planar domain bounded by $d_y,$ $y>0.$

The convex hull property for maximal surfaces gives $K_y \subset {\cal V}_r$ (even more, $K_y -c_y \subset {\cal V}_r-\partial({\cal V}_r)$).

If $y_1>y_2>0,$ then $c_{y_1}>c_{y_2}$ (that is to say, $t(c_{y_1}(\xi))>t(c_{y_2}(\xi))$ for any $\xi \in \s^1$). A standard application of the maximum principle gives that $K_{y_1}$ lies above $K_{y_2},$ and so $K_{y_1}\cap K_{y_2}=\emptyset.$ The smooth dependence of Plateau's problem solutions with respect to the boundary data implies that there is a unique $D_s(r) \in \{K_y \;:\; y \in [0,+\infty[\}$ such that $(0,0,s) \in D_s(r),$ $s \geq r.$ Furthermore, ${\cal D}(r)=\{D_s(r) \;:\; s \in [r,+\infty[\}$  defines a smooth foliation of ${\cal V}_r$ satisfying $(a)$ and $(b).$

In order to prove $(c),$ let us see that $D_s(r) \cap
\partial(\mb(r)) \neq \emptyset$ if and only if $D_s(r) \cap
\mb(r) \neq \emptyset.$ Suppose $D_s(r) \cap \mb(r) \neq
\emptyset,$ and reasoning by contradiction, assume that $D_s(r)
\cap \partial(\mb(r)) = \emptyset.$ As $\partial (D_s(r))\cap
\mb(r)=\emptyset,$ then $D_s(r) \cap \mb(r)$ is a family of
piecewise analytical Jordan curves lying in the interior of both simply connected 
surfaces.\footnote{If $M_1$ and $M_2$ are maximal surfaces and $\partial(M_j) \cap M_i=\emptyset,$ $\{i,j\}=\{1,2\},$ then
$M_1\cap M_2$ consists of a family of analytical proper
analytical arcs in $M_j-(\partial(M_1) \cup \partial(M_2)),$
$j=1,2,$ meeting equiangularly at points with the same normal.}
Hence we can find compact discs $S_1 \subset \mbox{Int}(D_s(r))$
and $S_2 \in \mbox{Int}(\mb(r))$ with common boundary in $D_s(r)
\cap \mb(r)$ and common projection on the plane $\Pi_0.$ Since
both discs are graphs over $\Pi_0,$ the maximum principle gives
$S_1=S_2,$ and by  an analytic continuation argument  $D_s(r)
\subset \mb(r).$ This is absurd and shows that $D_s(r) \cap
\partial(\mb(r)) \neq \emptyset.$

Finally, assume that $D_s(r) \cap \mb(r) \neq \emptyset.$ From
equations (\ref{eq:acausal}) and (\ref{eq:normal}),  $q_1:=D_s(r)
\cap \Gamma=D_s(r) \cap \Gamma_r$ consists of at most one point
where $D_s(r)$ and $\Gamma$  meet transversally (in case $D_s(r)
\cap \Gamma=\emptyset$ we make the convention $q_1=\emptyset$).
Likewise, from $(ii)$  in the preceding claim,
$q_2:=\partial(D_s(r)) \cap \mb(r)=c_s(r) \cap \beta_r$ is a
point. If $\alpha$ is an inextendible arc in $D_s(r) \cap \mb(r)
\subset {\cal V}_r$ then $\alpha$  is compact and with endpoints
lying in $(\partial(D_s(r)) \cap \mb(r)) \cup (\partial(\mb(r))
\cap D_s(r))=\{q_1,q_2\}.$ Taking into account that $D_s(r) \cap
\mb(r)$ contains no Jordan curves (reason as above), we get that
$q_1$ is a point and   $\alpha$ joins $q_1 $ with
$q_2,$  concluding $(c)$ and the lemma. \end{proof}

%
%

The following theorem has been inspired by  Meeks and Rosenberg ideas in \cite{meeks-ros}.

\begin{theorem} \label{th:omit}
If $\mb$ is $\omega^*$-maximal then any plane parallel to $\Sigma_\infty$ is transverse to $\mb.$

As a consequence, either $\Sigma_\infty$ is spacelike  and $\Theta_0(\mb)=\r$ or $\Sigma_\infty$ is lightlike and $\theta_\mb=2 \pi.$
\end{theorem}
\begin{proof} Up to a Lorentzian isometry, we assume that either $\theta^+=+\infty$ or $\theta_\mb<+\infty.$ In any case, Lemma \ref{lem:importante} guarantees that $\mb$  contains no upward   lightlike rays, and consequently, the foliation ${\cal D}(r)$ in Lemma \ref{lem:foliacion} makes sense, for any $r \in ]-\infty,-1],$ $(0,0,r) \notin \mb.$

From Theorem \ref{th:strong}, $\theta_\mb =2k \pi,$ $k \in \n,$ provided that  $\theta_\mb<+\infty,$ and in this case $\mb_\infty$ is a lightlike plane. 
Furthermore, by Corollary \ref{co:wedge} $\lim_{s \to +\infty} \Gamma'(s)$ and  $\lim_{s \to -\infty} \Gamma'(s)$ are lightlike vectors parallel to $\mb\infty,$ hence the theorem holds provided that $k=1.$

Therefore, it suffices to deal with the case    $\theta_\mb \in [4\pi,+\infty].$

Take $\theta_0 \in \Theta_0(\mb)$ such that  $I_0=]\theta_0-2\pi,\theta_0+2\pi[ \subset \Theta_0(\mb).$ In case $\theta^+=+\infty$ and $\theta^+>-\infty,$ and with the notation of Theorem \ref{th:strong},  we also impose that $I_0 \subset J_{\cal X}.$ Let $s_0 \in \r$ be the unique real number such that $\Theta_0(\Gamma(s_0))=\theta_0.$

In the sequel we only consider $r \in]-\infty,-1]$  such that $(0,0,r) \notin \mb$ and $\Gamma^{I_0} \subset \mbox{Int}({\cal C}^+_{\{(0,0,r)\}}).$

For any $s \in \r,$ set $\Sigma_\infty(s)$ the plane parallel to $\Sigma_\infty$ and passing through $(0,0,s).$

\begin{quote}
{\bf Claim 1:} {\em There exists a divergent sequence $\{R_k\}_{k \in \n} \subset ]-\infty,-1]$ such that $\{{\cal D}(R_k)\}_{k \in \n}$ converges in the ${\cal C}^1$-topology to the foliation of $\l^3$ by planes parallel to $\Sigma_\infty.$\footnote{This means that for any compact interval $I \subset \r,$ $\{D_{s}(R_k)\}_{k \in \n} \to \Sigma_\infty(s)$ in the ${\cal C}^1$-topology  uniformly on $s\in I.$}}
\end{quote}
\begin{proof}
For $r\leq -1$ and  $n \in \n,$  label $r(n):=\frac{r}{\lambda_n}.$ Since $\mb$ has no upward lightlike rays, then $F_0(s_0,\cdot)$ and $\partial ({\cal V}_{r(n)})$ meet at a unique point $q_{r(n)}\in \beta_{r(n)}.$ Call $E(r(n))$  as the unique maximal disc in ${\cal D}({r(n)})$ containing  $q_{r(n)},$ and let us show that $\{\lambda_n E(r(n))\}_{n \in \n}$ converges in the ${\cal C}^0$-topology to $\Sigma_\infty \cap {\cal V}_{r}$ as graphs over $\Pi_0.$

Since  $\{\mb_n^{I_0}\}_{n \in \n}$ converges uniformly on compact subsets to the twice-covered once punctured plane $\Sigma_\infty-\{O\}$ (see Theorem \ref{th:strong}) and $\mb^{I_0} \cap {\cal V}_{r(n)} \subset \mb(r(n)),$ then $c_{r(n)}:=\lambda_n \left( \mb^{I_0}\cap \partial({\cal V}_{r(n)})\right) \subset \lambda_n \beta_{r(n)}$ converge as $n \to \infty$ to the twice-covered Lorentzian circle $c:=\Sigma_\infty \cap \partial({\cal V}_{r})$ (a parabola when  $\Sigma_\infty$ is lightlike). Furthermore, $\lambda_n \partial(E(r(n)))$ and $c_{r(n)}$ meet only at $\lambda_n q_{r(n)}$ in a transversal way, and any of the two components of $c_{r(n)}-\lambda_n q_{r(n)}$  converges uniformly  as $n \to \infty$ to $c.$ Taking into account that $\lambda_n \partial(E(r(n)))$ lies in between these components, we deduce that $\{\lambda_n \partial(E(r(n)))\}_{n \in \n}\to c$ too.
If $\Sigma_\infty$ is spacelike,  $c$ is a closed curve and the  continuous dependence of Plateau's problem solutions with respect to the boundary data gives that $\{\lambda_n E(r(n))\}_{n \in \n}\to \Sigma_\infty \cap {\cal V}_{r}$ in the ${\cal C}^0$-topology.

Assume now that $\Sigma_\infty$ is lightlike (and $c$ is a parabola), and call  $E_{r(\infty)}:=\lim_{n \to \infty} \lambda_n E(r(n)).$ Note that
$\lambda_n E(r(n)) \subset {\cal V}_{r}$ for every $n \in \n,$ hence $E_{r(\infty)} \subset {\cal V}_r.$ Since $E_{r(\infty)}$ is a PS graph  and $\partial(E_{r(\infty)})=c,$  equation (\ref{eq:acausal})  yields that $E_{r(\infty)} \subset \left(\cap_{p \in c} \overline{\mbox{Ext} ({\cal C}_p)}\right) \cap {\cal V}_r.$ This proves that $E_{r(\infty)}$ lies in the slab bounded by $\Sigma_\infty$ and $\Sigma_\infty(r),$ and as a consequence, $E_{r(\infty)}$ must contain lightlike segments (otherwise, $E_{r(\infty)}$ would be a parabolic maximal graph, hence a planar domain in $\Sigma_\infty$  by Corollary \ref{co:halfspace}, which is  absurd). From Remark \ref{re:strip},  $E_{r(\infty)} \cap \Sigma_\infty$ must contain a lightlike half line $L$ with initial point in $c,$ and therefore $S:=E_{r(\infty)}-\left(\Sigma_\infty \cap \pi^{-1}(\Pi_0-\pi(E_{r(\infty)}))\right)$ is an entire PS graph over $\Pi_0$ containing the complete straight line determined by $L.$ Lemma \ref{lem:basico} shows that $S=\Sigma_\infty$ and $E_{r(\infty)}=\Sigma_\infty \cap {\cal V}_{r},$  proving that $\{\lambda_n E(r(n))\}_{n \in \n}\to \Sigma_\infty \cap {\cal V}_{r}$ too.

Take a divergent sequence $\{r_k\}_{k \in \n}$ in $]-\infty,-1]$ and observe that $\{\Sigma_\infty \cap {\cal V}_{r_k}\}_{k \in \n} \to \Sigma_\infty$ in the ${\cal C}^0$-topology. By a standard diagonal process, we can find a divergent sequence $\{n_k\}_{k \in \n}\subset \n$  such that $\{\lambda_{n_k} E(r_k(n_k))\}_{k \in \n} \to \Sigma_\infty$ in the ${\cal C}^0$-topology as graphs over $\Pi_0.$
Define $R_k:=r_k(n_k),$ $k \in \n,$ and let us show that $\{R_k\}_{k \in \n}$ solves the claim.

To do this, let $I \subset \r$ be a compact interval, and take a sequence $\{s_k, \;k \in \n\} \subset I$ converging to $s \in I.$ It suffices to check that $\{D_{s_k}(R_k)\}_{k \in \n} \to \Sigma_\infty(s)$ in the ${\cal C}^1$-topology as graphs over $\Pi_0.$

Let us see first that  $\{\lambda_{n_{k}} D_{s_k}(R_k)\}_{k \in
\n}\to \Sigma_\infty$ in the ${\cal C}^1$-topology.   Indeed,
note that $\lambda_{n_{k}} (0,0,s_k) \in \lambda_{n_{k}}
D_{s_k}(R_k)$ and $\{\lambda_{n_{k}} (0,0,s_k)\}_{k \in \n} \to
O.$ Thus,  $\{\lambda_{n_{k}} D_{s_k}(R_k)\}_{k \in \n}$
converges in the ${\cal C}^0$-topology to an entire PS graph
$\Sigma'_\infty$ passing through $O.$ As  either
$D_{s_k}(R_k)=E(R_k)$  or $D_{s_k}(R_k)\cap E(R_k)=\emptyset$ for
any $k,$ then $\Sigma'_\infty$ lies in one of the closed half
spaces bounded by $\Sigma_\infty.$ Using Calabi's theorem, Remark
\ref{re:strip} and Lemma \ref{lem:basico} we infer that
$\Sigma'_\infty$ is a non timelike plane passing through $O,$
hence $\Sigma'_\infty=\Sigma_\infty.$ From Proposition
\ref{pro:converge}, $\{\lambda_{n_{k}} D_{s_k}(R_k)\}_{k \in \n}
\to \Sigma_\infty$ in the ${\cal C}^1$ topology and we are done.

Finally, let  $u_k:\pi(D_{s_k}(R_k))\to\r$ and $v_k:\lambda_{n_k}
\pi(D_{s_k}(R_k)) \to \r$ be the functions determining the graphs
$D_{s_k}(R_k)$ and $\lambda_{n_k} D_{s_k}(R_k)$ respectively,  $k
\in \n.$  Proposition \ref{pro:converge} gives that $\{\nabla
v_k\}_{k \in \n}\to \sigma$ in the ${\cal C}^0$ topology, where
$\sigma$ is the gradient of the linear function defining
$\Sigma_\infty.$ Since $v_k(\lambda_{n_k} x)=\lambda_{n_k}
u_k(x),$ we infer that $\nabla v_k(\lambda_{n_k}x)=\nabla u_k(x)$
for any $x \in \pi(D_{s_k}(R_k)).$ This gives $\{\nabla u_k\}_{k
\in \n}\to \sigma$ in the ${\cal C}^0$-topology, and so $\lim_{k
\to +\infty} D_{s_k}(R_k)=\Sigma_\infty(s)$ in the ${\cal
C}^1$-topology. \end{proof}

To finish the theorem, reason by contradiction and suppose there is $q \in \l^3$ such that $T_q \mb=\Sigma_\infty.$ Recall that the conformal parameterization of $\mb$   extends to the conformal mirror $\mb^*$ of $\mb$ by folding back at $\Gamma.$ In particular, the holomorphic Gauss map $g$ extends by Schwarz reflection to $\mb^*$  as well. Take a closed disc $U \in \mb \cup \mb^*$ such that  $q\in \stackrel{\circ}{U}$ and $U \cap g^{-1}(g(q))=q.$ Let $m\geq 1$ denote the multiplicity of $g$ at $q,$ that is to say, the winding number of $g({\partial(U)})$ around $g(q).$  For any $k \in \n$ such that $U\cap \mb \subset {\cal V}_{R_k}$ and for any $p \in U\cap \mb,$  let $s_k(p) \in \r$ denote the unique real number such that $p \in D_{s_k(p)}(R_k)$ and call
$U_k=\cup_{p \in U\cap \mb} D_{s_k(p)}(R_k).$ From equation (\ref{eq:acausal}), $U_k \subset \cup_{p \in U\cap \mb} \overline{\mbox{Ext}({\cal C}_{p})},$ and so $I_k:=\{s_k(p) \;:\; p \in U\cap \mb\}$ lies in the compact interval $I=\{s \in \r \;:\; (0,0,s) \in   \cup_{p \in U\cap \mb} \overline{\mbox{Ext}({\cal C}_{p})}\}.$ In other words, $U_k \subset \cup_{s \in I} D_{s}(R_k)$ for any $k \in \n$ satisfying $U\cap \mb \subset {\cal V}_{R_k}.$

For any $k \in \n$ such that $U\cap \mb \subset {\cal V}_{R_k}$ and $p \in U \cap \mb,$ let $g_{k,p}:D_{s_k(p)}(R_k) \to \d$ be the holomorphic Gauss map of $D_{s_k(p)}(R_k),$ and set $h_k:U \to \c,$  $(h_k|_{U \cap \mb})(p):=g(p)-g_{k,p}(p)$  and $(h_k|_{U \cap \mb^*})(p):=g(p)-g_{k,p^*}(p^*).$ 
Labeling $s(p)\in I$ as  the unique real number such that $p \in \Sigma_\infty(s(p)),$ $p \in U \cap \mb,$  Claim 1 gives that $\lim_{k \to \infty} s_k(p)=s(p)$ and
$\lim_{k \to +\infty}D_{s_k(p)}(R_k) =\Sigma_\infty(s(p))$ in the ${\cal C}^1$-topology  uniformly on $p \in U \cap \mb.$ Therefore,  $\lim_{k \to +\infty} h_k = g-g(q)$ uniformly on $U,$ hence for large enough $k$ the winding number of $h_k({\partial(U)})$ around the origin is equal to $m.$ However $h_k|_{U \cap \mb^*}$ never vanishes, and so  we can find for $k$ large enough a point $q_k \in (U\cap \mb)-\partial(U)$  such that $h_k(q_k)=0$   (that is to say, $T_{q_k}  D_{s_k(q_k)}(R_k)=T_{q_k} \mb$),  contradicting that ${\cal D}(R_k)$ is transverse to $\mb$ for any $k$ and proving that $\mb$ is transverse to $\Sigma(s),$ $ s \in \r.$ 

If $\Sigma_\infty$ is lightlike, we can find $q \in \Gamma$ such that $T_q \mb$ is parallel to $\Sigma_\infty$ (recall that $\theta_\mb\geq 4\pi$),  a contradiction. Thus $\Sigma_\infty$ is spacelike, and by  Theorem \ref{th:strong}, $]\theta^+,\theta^-[=\r.$  This concludes the proof. \end{proof}

\begin{corollary} \label{co:fujimoto}
If $\mb$ is $\omega^*$-maximal then the blow-down plane $\Sigma_\infty$ does not depend on the blow-down sequence.
\end{corollary}

\begin{proof} If $\theta_\mb=2 \pi,$ $\Sigma_\infty$ is the limit plane of $\mb$ at infinity and the corollary holds.

Assume that $\theta_\mb=+\infty,$ take a new blow-down sequence $\{\lambda_n'\}_{n \in \n}$ and construct the corresponding blow-down PS multigraph ${\cal X}':{\mb}'_\infty \to \l^3$ and blow-down spacelike plane  ${\Sigma}'_\infty:={{\cal X}'}({\mb}'_\infty).$

Reason by contradiction and  suppose ${\Sigma}'_\infty \neq \Sigma_\infty.$ Consider an interval $I \subset \Theta_0(\mb)=\r$ of length $2\pi,$ fix $p \in \Pi_0-O$ and take $p_n\in \pi^{-1}(p)\cap  (\lambda_n \cdot \mb^I),$  $p'_n\in \pi^{-1}(p) \cap (\lambda_n' \cdot \mb^I),$ $n \in \n$ (well defined provided that $n$ is large enough). We know that $\lim_{n \to \infty} {\cal N}(\frac{1}{\lambda_n}p_n)=\zeta$ and  $\lim_{n \to \infty} {\cal N}(\frac{1}{\lambda'_n} p'_n)=\zeta',$ where ${\cal N}$ is the Lorentzian Gauss map of $\mb$ and $\zeta$ and $\zeta' \subset \h^2_-$ are the unitary normal vectors to $\Sigma_\infty$ and $\Sigma'_\infty,$ respectively.
By a connectedness argument, we can find $\zeta_0 \in \h^2_--\{\zeta,\zeta'\}$ and divergent sequence $\{q_n\}_{n \in \n} \subset \mb^I$ such that $\lim_{n \to \infty} {\cal N}(q_n)=\zeta_0.$ As above, $\zeta_0$ is the unitary normal to the blow-down plane $\Sigma_\infty''$ associated to  blow-down sequence $\{\frac{1}{\|q_n\|_0}\}_{n \in \n}.$

On the other hand, let $(g,\phi_3)$ denote the Weierstrass data of $\mb$ (see equation (\ref{eq:wei})) and  extend $(g,\phi_3)$ by Schwarz reflection  to the double $\hat{\mb}$ of $\mb.$  Then, consider the conjugate minimal immersion $X^*:\hat{\mb} \to \r^3$ associated to the same Weierstrass data $(g,\phi_3),$  see equation (\ref{eq:weiconju}), and recall that the metrics on $\hat{\mb}$  induced by the maximal and minimal immersions are given by  $ds^2=\frac{1}{4}|\phi_3|^2(\frac{1}{|g|}-|g|)^2$ and  $ds_0^2=\frac{1}{4}|\phi_3|^2(\frac{1}{|g|}+|g|)^2,$ respectively.

Let us show that $ds_0^2$ is complete. Since  the mirror involution is an isometry of $(\hat{\mb},ds_0^2),$ it suffices to check that divergent curves in $\mb$ have infinite length with respect to $ds_0^2.$ Indeed, set $\alpha \subset \mb$ a divergent curve. Since the first and second coordinate functions of a maximal surface and its conjugate minimal one are the same, then the Euclidean length $L_0(\alpha)$ of  $\alpha$  is greater than or equal to the one  $L_0(\pi(\alpha))$ of $\pi(\alpha).$  The spacelike property of $\mb$ and the divergence of $\alpha$ give that $L_0(\pi(\alpha))=+\infty,$ and so $L_0(\alpha)=+\infty.$

By Theorem \ref{th:omit}, ${\cal  N}:\mb \to \h^2_-$ omits the values $\zeta,$ $\zeta'$ and $\zeta_0 \in \h^2_-,$ hence the Gauss map of $X^*$ omits six complex values. This contradicts Fujimoto's theorem \cite{fujimoto1,fujimoto2} and proves the corollary. \end{proof}

\section{The Uniqueness Theorems} \label{sec:uniq}
In this section we prove the main results of this paper. We start with the following:
\begin{theorem}[Uniqueness of the Enneper surface]
The only properly embedded  $^*$maximal surface  with connected boundary and finite rotation number is, up to Lorentzian congruence, the Enneper surface $E_1.$
\end{theorem}
\begin{proof} Let $\mb$ be a properly embedded $^*$maximal surface  with connected boundary and $\theta_\mb<+\infty.$

Since $\mb$ is a multigraph of finite angle, Corollary \ref{co:wedge} gives that $\mb$ is conformally equivalent to $\overline{\d}-\{1\}.$ Let $(g,\phi_3)$ denote the Weierstrass data of $\mb.$ From Theorem \ref{th:omit}, the holomorphic map $g:\overline{\d}-\{1\} \to \overline{\d}$ is one to one on $\partial (\d)-\{1\},$ and so, up to a Lorentzian isometry, we can suppose that $g(z)=z.$
On the other hand, equation (\ref{def:simetria}) leads to $\phi_3=h(z)\frac{z}{(z-1)^4} dz$ (note that the mirror involution is given by  $J(z)=1/\overline{z}$), where $h:\c \to \c$ is a meromorphic function satisfying $h \circ J=\overline{h}.$ Since the 1-forms $\phi_j$ given in (\ref{eq:wei}) have no common zeroes, $h$ never vanishes on ${\d}.$ Furthermore, as the unique end of $\mb$ corresponds to $z=1,$ $h$ has no poles in ${\d}$ as well. The symmetry condition $h \circ J=\overline{h}$ gives that the zeroes and poles of $h,$ if they occur, lie in $\partial{\d}.$ However, Lemma \ref{lem:prime} implies that $\phi_3$ never vanishes on $\partial{\d}-\{1\},$ and so $h$ must be a real number different from zero. Up to scaling and a conformal reparameterization, $(g,\phi_3)$
are the Weierstrass data of $E_1$ (see Section \ref{sec:ejemplos}). \end{proof}

In the sequel we will deal with the uniqueness of  properly embedded  $\omega^*$-maximal surfaces with connected boundary and infinite rotation number. This part of the paper has been mainly inspired by Meeks-Rosenberg work \cite{meeks-ros}.

Let $\mb$ denote a properly embedded $\omega^*$-maximal surface  with connected boundary and $\theta_\mb=+\infty.$ From Theorem \ref{th:omit}, $\Theta_0(\mb)=\r,$ $\Sigma_\infty$ is a spacelike plane, any plane parallel to $\Sigma_\infty$ meets $\mb$ transversally into a family of pairwise disjoint proper analytical arcs.

Let $\Sigma$ be a plane parallel to $\Sigma_\infty,$ and label $\Sigma^+$ and $\Sigma^-$ as the two closed half spaces in $\l^3$ bounded by $\Sigma.$ As $\Sigma$ is spacelike then $q:=\Sigma \cap \Gamma$ is a single point. We set $\mb(q)=\mb \cap \overline{\mbox{Int}({\cal C}_q)},$ $\mb^+(q):=\mb(q) \cap \Sigma^+$ and $\mb^-(q):=\mb(q) \cap \Sigma^-.$ Note that  equation (\ref{eq:normal}) gives $\Gamma \subset \mb (q).$
Since the arcs $F_0(s,\cdot)$ have slope $\leq 1$ and $\mb$ has no lightlike asymptotic rays (see Lemma \ref{lem:importante}), $F_0(s,\cdot)\cap \mb (q)$ is compact and connected for any $s \in \r.$ We deduce that $\mb^+(q)$ and $\mb^-(q)$ are simply connected regions in $\mb$ with connected boundary. Moreover, $\mb^+(q)\cap \mb^-(q)=\{ q\},$ and so $\mb(q)$ is connected too.

Consider a region $S \subset \mb$ (in most cases we will deal with $S=\mb$). The closure of a connected component of $\mbox{Int}(S)-\Sigma$ is defined to be a $\Sigma$-{\em region} of $S.$

 A $\Sigma$-region of $\mb$ is said to be a {\em finite} (resp., {\em infinite})  if its boundary has finitely (resp., infinitely) many pairwise disjoint proper arcs. A finite $\Sigma$-region of $\mb$ is said to be {\em simple} if it has connected boundary. Any $\Sigma$-region ${W}$ of $\mb$ is parabolic (see Corollary \ref{co:halfspace}) and simply connected, hence  conformally equivalent to $\overline {\d}-E,$ where $E \subset \partial (\d)$ is a totally disconnected compact zero measure subset. For convenience, we will identify ${W}$ and $\overline {\d}-E$ and call $E$ as the set of ends of $W.$

\begin{lemma} \label{lem:primero}
If $W \equiv \overline {\d}-E$ is an infinite $\Sigma$-region of $\mb$  then either $\mb^+(q) \subset W$ or $\mb^-(q) \subset W.$

Moreover, the endpoint of  $\Gamma  \cap W$ is the unique limit end $*$ of $E.$ 
\end{lemma}
\begin{proof}

Up to a Lorentzian isometry,  we will suppose that $\Sigma=\Sigma_\infty=\Pi_0,$ $q=\Gamma \cap \Sigma=O$ and $\Sigma^+=\{t \geq 0\}$ (hence $\Sigma^-=\{t \leq 0\}$). For simplicity we write $\Gamma_0$ instead of  $\Gamma \cap W.$

First of all recall that $\mb^+(O)$ and $\mb^-(O)$ are connected, hence they lie in the $\Sigma$-regions of $\mb$  containing $\Gamma\cap \{t>0\}$ and $\Gamma\cap \{t<0\},$ respectively.

Take a limit end $*$ of $E$ and an auxiliary point $q_0 \in \partial (\d)-\{*\}.$ Label $c_1$ and $c_2$ as the two open arcs in $\partial (\d)-\{*,q_0\},$ and consider sequences $\{e_n\;:\; n \in \n\} \subset \overline{c_1} \cap E$ and  $\{e'_n\;:\; n \in \n\} \subset \overline{c_2}\cap E$ converging to $*.$ Without loss of generality, suppose that $\{e_n \;:\; n \in \n\}$ is not a finite set of ends.

Reason by contradiction, and assume that either $\Gamma_0=\emptyset$ or $\Gamma_0 \neq \emptyset$ and $\Gamma_0$ does not diverge to $*.$ Thus there exists a compact arc $c \subset W-\Gamma$  connecting two points of $\partial (W),$ and such that $W-c$ has a connected component $W'$ with infinitely many boundary components, disjoint from $\Gamma$ and containing $*$ among its  limit ends. Without loss of generality, we can also suppose that $\{e_n\;:\; n \in \n\}\cup \{e'_n\;:\; n \in \n\}$ are ends of $W'.$

Since $W' \cap \mb(O)$ is compact (just observe that $\partial(W') \cap \mb(O)=c \cap \mb(O)$ is compact and  $W'\cap \Gamma=\emptyset$), then  $\pi|_{W'}:W' \to \Pi_0$ is proper. Take $R>0$ such that $c \subset \{(x,t) \in \l^3 \;:\; \|x\|_0<R\}$ and consider the connected component $W_R$ of $W' \cap \{(x,t) \in \l^3 \;:\; \|x\|_0\geq R\} $ with infinitely many boundary arcs. It is clear that  $W_R$ is biholomorphic to $D_R-E_R,$ where $D_R \subset \overline {\d}$ is a closed topological  disc  and $E_R=\partial(D_R) \cap E.$ Furthermore and as above we can suppose that  $\{e_n\;:\; n \in \n\}\cup \{e'_n\;:\; n \in \n\} \cup \{*\} \subset E_R.$

Put $\alpha_R:=\partial(W_R)\cap \{(x,t) \in \l^3 \;:\; \|x\|_0=R\},$ and  set $\sigma$ the connected component of $\partial(W_R)$  containing $\alpha_R.$
Let $\hat{W}_R$ denote the region in $\mb$ bounded by $\sigma$ and disjoint from $\Gamma$ (obviously $W_R\subset \hat{W}_R$).

Let us see that $\pi|_{\hat{W}_R}:\hat{W}_R \to \pi(\hat{W}_R)$ is a diffeomorphism. Let us see that $\pi|_{\partial(W_R)}$ is injective. Indeed, since $\pi|_{\partial (W_R)-\alpha_R}$ is the identity map, it suffices to prove that $\pi|_{\alpha_R}$ is injective. Assume without loss that $W_R \subset \{t \geq 0\}$ and note that ${W}_R$ separates the region $\{(x,t) \in \l^3 \;:\; \|x\|_0\geq R,\; t \geq 0 \}.$ Therefore, for an arc $\alpha \subset \alpha_R= W_R \cap \{(x,t) \in \l^3 \;:\; \|x\|_0= R\}$ there can not be another arc $\beta \subset \alpha_R$ immediately above of below $\alpha.$ Otherwise, the Euclidean normal vectors to $W_R$ along $\alpha$ and $\beta$ would lie in different hemispheres, which contradicts that the projection $\pi$ orients $W_R.$  
Therefore, $\pi|_{\hat{W}_R}:\hat{W}_R \to \pi(\hat{W}_R)$ is  a proper local diffeomorphism, hence a global diffeomorphism from the simply connectedness of $\pi(\hat{W}_R).$

Set $\{\Omega_n \;:\; n \in \n\}$  the countable family of connected components in $\pi(\hat{W}_R-W_R),$  and let  $G_n=\{(x,u_n(x)) \;:\; x \in \overline{\Omega}_n\}$ denote the maximal graph in  $\hat{W}_R-\mbox{Int}(W_R)$ satisfying $\pi(G_n)=\overline{\Omega}_n,$ $n \in \n.$ It is clear that   $\overline{\Omega_n} \cap \overline{\Omega_m} =\emptyset,$ $m \neq n,$ and $u_n|_{\partial(\Omega_n)}=0.$ The desired contradiction will comes from Theorem \ref{th:li-wang2}, provided that  $|\nabla u_n| <1-\epsilon$ for any $n$ for  a suitable $\epsilon>0.$

To check the last inequality,  reason by contradiction and  suppose there exists a sequence $\{p_n\}_{n \in \n},$ where $p_n \in G_n,$  such that  $\{\nabla u_n(p_n)\}_{n \in \n} \to 1$ (or in other words, $\{|g(p_n)|\}_{n \in \n} \to 1,$ where $g$ is the holomorphic Gauss map of $\mb$). Since $G_n$ is parabolic (see Corollary \ref{co:halfspace}) and $|g|$ subharmonic, then $|g|(p_n) \leq \int_{\partial(G_n)} |g| \, d\mu_{p'_n},$ where $d\mu_{p'_n}$ is the harmonic measure respect to a given point $p'_n\in G_n.$ As $d\mu_{p'_n}$ is a probabilistic measure (i.e.,  $\int_{\partial(G_n)}  d\mu_{p'_n}=1$), then we can find $q_n \in \partial (G_n) \subset \Pi_0$ satisfying $|g(q_n)| \geq |g(p_n)|.$
Taking into account that $\mb$ is proper and the  $G_n$'s are pairwise disjoint, we deduce that $\{\lambda'_n:=\frac{1}{\|q_n\|_0}\}_{n \in \n} \to 0.$ Like in Section \ref{sec:main}, consider the sequence $\{{\mb}'_n:=\lambda'_n \mb\}_{n \in \n},$ and  the corresponding blow-down  PS multigraph ${{\cal X}'}:{\mb}'_\infty \to \l^3.$
Assume that $\{p_n:=\frac{q_n}{\|q_n\|_0}\}_{n \in \n} \to p \in \Pi_0,$ take an open disc $D \subset \Pi_0$ centered at $p$ of radius $<\frac{1}{2}$ and call $N_n$ as the closure of the connected component of ${\mb}'_n\cap \pi^{-1}(D)$ containing $p_n,$ $n \in \n.$
From equation (\ref{eq:acausal}), $N_n \subset \overline{\mbox{Ext}({\cal C}_{p_n})} \cap \pi^{-1}(D),$ and therefore $N_n \cap \overline{\mbox{Int}({\cal C}_0)}=\emptyset,$ for any  $n \in \n.$ On the other hand, equation (\ref{eq:normal}) and the fact $O \in \Gamma$ give that  $\lambda_n' \Gamma \subset \overline{\mbox{Int}({\cal C}_{0})}$ and  prove that $N_n\cap (\lambda_n' \cdot \Gamma)=\emptyset,$ for any $n \in \n.$ As a consequence $\pi(N_n)=D$ and $\pi|_{N_n}:N_n \to D$ is a diffeomorphism,  $n \in \n.$

The hypothesis $\{|g(p_n)|\}_{n \in \n} \to 1$ implies that
$N_\infty:=\lim_{n \to \infty} N_n$ contains a lightlike half
line  passing through $p$ (see Theorem \ref{th:converge}).
Therefore $N_\infty$ and the plane ${{\cal
X}'}({\mb}'_\infty)=\Pi_0$ (see Corollary \ref{co:fujimoto}) meet
transversally at $p.$ This contradicts that $\Pi_0\cup N_\infty$ lies in the
limit set of the sequence of embedded  surfaces  $\{{\mb}'_n\}_{n
\in \n}$ and  concludes the proof. \end{proof}

\begin{proposition} \label{pro:porfin}
If $\Sigma$ si a plane parallel to $\Sigma_\infty$ then $\mb\cap \Sigma$ consists of a proper regular arc.
\end{proposition}
\begin{proof} Like in the preceeding lemma, assume that $\Sigma=\Sigma_\infty=\Pi_0$ and $\Gamma \cap \Sigma=O.$

Let $c_0$ denote the unique proper divergent arc in $\mb \cap \Pi_0$ meeting $\Gamma$ (that is to say, the one with initial point $O$), and  set $\mb^+$ (resp., $\mb^-$) the region in $\mb$ bounded by $c_0 \cup (\Gamma \cap \{t \geq 0\})$ (resp., $c_0 \cup (\Gamma \cap \{t \leq 0\})$). It is clear that $W^+ \subset \mb^+$ and $W^- \subset \mb^-,$ where $W^+$ and $W^-$ are the $\Sigma$-regions of $\mb$ containing $\mb^+(O)$ and $\mb^-(O),$ respectively.

Lemma \ref{lem:primero} implies that any region $U \subset \mb$  disjoint from $\Gamma$ with $\partial(U) \subset \Pi_0$   contains  finitely many $\Sigma$-regions of $\mb.$ Thus, we can find a divergent
arc $\beta^+ \subset \mb^+$ disjoint from $\Gamma,$  meeting
$c_0$ just at the initial point of $\beta^+,$  and meeting twice
any connected component of $(\mb^+ \cap \Pi_0)-c_0.$ In a similar
way we define $\beta^-,$ and without loss of generality  suppose
$c_0 \cap \beta^+=c_0 \cap \beta^-.$  The proper arc
$\beta=\beta^+ \cup \beta^-$ is disjoint from $\Gamma$ and splits $\mb$ into two connected components. We set $\mb_\beta$  the closure
of the connected component of $\mb-\beta$ disjoint from $\Gamma.$

Let  $V$ be a $\Sigma$-region of  $\mb_\beta,$ obviously non compact. $V$ is said to be a {\em middle} $\Sigma$-region of $\mb_\beta$ if $\partial (V) \cap \beta$ is compact. Otherwise, $V$ is said to be a {\em tail} $\Sigma$-region of $\mb_\beta.$  Two different $\Sigma$-regions of $\mb_\beta$ are said to be contiguous if they share a non compact boundary arc in $\Pi_0.$
Obviously $\mb_\beta \cap \mb^+$ (resp., $\mb_\beta \cap \mb^-$) contains at most one tail $\Sigma$-region, and it contains a tail $\Sigma$-region if an only if $\partial(W^+)$ (resp., $\partial(W^-)$) contains finitely many components. Moreover, $\mb_\beta$ contains no middle $\Sigma$-regions  if and only if $\mb_\beta \cap\Pi_0=c_0.$

Let $t^*:\mb \to \r$ denote the harmonic conjugate of the third coordinate function $t:\mb \to \r,$ and  consider the holomorphic function $h:=t+i t^*:\mb \to \c.$

\begin{quote}
{\bf Claim 1:} {\em If $V_0$ is a middle $\Sigma$-region of $\mb_\beta$ then  $t|_{V_0}$ is unbounded.}
\end{quote}
\begin{proof}
Suppose that $t|_{V_0}$ is bounded, and without loss of generality  assume that $t(V_0)\subset ]-\infty,0].$
Since $V_0$ is parabolic (see Corollary \ref{co:halfspace}), there is a biholomorphism $T:V_0 \to A=\{z \in \overline{\d}-\{0\}\;:\; |\arg(z)|\leq \frac{\pi}{2}\}$  such that $T(\partial(V_0) \cap \Pi_0)=\{z \in A \;:\; \mbox{Re}(z)=0\}.$  Up to the identification $V_0 \equiv A$ via $T,$  $f:A \to \c,$ $f:=e^h,$  extends by Schwarz reflection to a bounded function, that we keep calling $f,$ on $\overline{\d}-\{0\}.$ by Riemann's removable singularity theorem, $f$ extends holomorphically to $\overline{\d}.$ Furthermore, $f$ has no zeroes in $\overline{\d}$ because $t=\log(|f|)\leq 0$ is bounded, and thus $h=\log(f)$ has well defined limit at $0.$ Thus, $h|_{V_0}$ is bounded and has well defined limit at its unique end.

Let $V_1$ be a middle $\Sigma$-region of $\mb_\beta$ contiguous
to $V_0.$ Since $h(V_1) \subset \{z \in \c \;:\; \mbox{Re}(z)
\geq 0\}$ and $h(V_0)$ is bounded,  $h(V_0 \cup V_1)$ omits
infinitely many complex values, and consequently it is a normal
function. From the conformal point of view, $D_1:=V_0 \cup V_1$
is biholomorphic to the  $\overline{\d}-\{1\}.$ Identifying $D_1
\equiv \overline{\d}-\{1\},$  $h$ has a well defined finite
limit $w_0$ along arcs $\alpha \subset V_0 \subset
\overline{\d}-\{1\}$ diverging to $1.$ Basic sectorial theorems
for normal functions imply that $h|_{D_1}$ has well defined {\em
finite} angular limit $w_0$  at the end $1.$ In particular,
$t|_{V_1}$ can not have asymptotic curves  with asymptotic value
$\infty,$ which proves that $t|_{V_1}$ is bounded. Reasoning as
at the beginning of the claim, $h|_{V_1}$ is bounded and  has
well defined finite limit $w_0$ at its unique end. Repeating this
argument for successive contiguous middle $\Sigma$-regions, we
conclude that $h$  has  limit $w_0$ at the end of any middle
$\Sigma$-region.

Now we can finish the claim. As we are assuming that $\mb_\beta$
contains middle $\Sigma$-regions, then there is a
$\Sigma$-region $U$ of $\mb$ with $\partial (U) \subset \Pi_0.$
The parabolicity of $U$ and Claim 1 show that $U$ is
biholomorphic to $\overline{\d}-\{w_1,\ldots,w_k\},$  where
$\{w_1,\ldots,w_k\}\subset \partial(\d).$ Since
$t|_{\overline{U}}$ is bounded  and $t|_{\partial (U)}=0,$ we get
$t|_U=0,$ which is absurd. \end{proof}

From Theorem \ref{th:omit} (see also the comments below Definition \ref{def:simetria}),  ${\cal Y}:=\frac{1}{(t \circ {\cal N})^2}\nabla t^*$ is well defined and never vanishes  on $\mb.$ Thus, any integral curve of ${\cal Y}$ is a proper arc in $\mb$  contained in a horizontal plane\footnote{Thw symbol $\nabla$ means  gradient  with respect to the metric $ds^2$ induced by $\langle,\rangle$.}. Furthermore, since ${\cal Y}$ is a spacelike field and $\Gamma$ is lightlike, the integral curves of ${\cal Y}$  are transverse to $\Gamma.$
Consider the flow
${\cal F}:\Gamma \times [0,+\infty[$ of ${\cal Y}$ and define ${\cal D}\subset \mb$ as the open subdomain that is the image ${\cal F}(\Gamma \times [0,+\infty[).$
For the sake of simplicity, we denote by $c_s$ integral curve ${\cal F}(\Gamma(s),[0,+\infty[),$  $s \in \r.$ In other words, $c_s$ is the connected component of $\mb \cap \{t=s\}$ meeting $\Gamma.$

To finish the theorem, it suffices to check that ${\cal D}=\mb.$  Reason by contradiction, and assume that $\mb-{\cal D} \neq \emptyset.$ Therefore, $\partial ({\cal D})=\Gamma \cup C,$ where $C \neq \emptyset$ is a collection of pairwise disjoint proper integral curves of ${\cal Y}$ disjoint from $\Gamma.$

Since  arcs in $C$ lie in horizontal planes,  we can suppose up to a translations that $C \cap \Pi_0\neq \emptyset.$ Let $c$ be a proper arc in $C \cap \Pi_0.$ Fix $p_0 \in c$ and let $\delta:[-\epsilon,\epsilon] \to \mb$ be the integral curve of $\nabla t$ with initial condition $\delta(0)=p_0.$ Since $p_0 \in C \subset \partial ({\cal D}),$ we infer that $\delta(]0,\epsilon]) \subset {\cal D},$ provided that $\epsilon$ is small enough. Write $t(\delta(\epsilon))=a>0$ and note that $t(\delta(]0,\epsilon]))=]0,a].$ For any $s \in ]0,a],$ let $\hat{c}_s\subset c_s$ denote the compact arc joining $\Gamma(s)={\cal F}(\Gamma(s),0)=\Gamma \cap c_s$ and $(t \circ\delta)^{-1}(s) .$
From the  choice of $c,$ the curves $\{\hat{c}_s \;:\; s \in ]0,a]\}$ converge as $s \to 0$ uniformly on compact subsets of $\mb$ to $c_0 \cup \hat{c},$ where $\hat{c}\subset C\cap \Pi_0$ is a  collection of proper subarcs in $C$ and  $p_0 \in \hat{c}.$  By Lemma \ref{lem:primero}, $\hat{c}$ has finitely many connected components, one of then being a divergent subarc of $c$ with initial point $p_0.$

Set $V=(\cup_{s \in ]0,a]}\hat{c}_s) \cup c_0 \cup \hat{c},$ and note that $V$ is a region in $\mb$ homeomorphic to  a closed disc minus a finite set of boundary points and with boundary $\partial(V)= \Gamma([0,a]) \cup \delta(]0,\epsilon]) \cup \hat{c}_a\cup c_0 \cup \hat{c}.$ Since $\Gamma([0,a])\cup \delta(]0,\epsilon])\cup \hat{c}_a $ is compact and $c_0 \cup \hat{c} \subset \Sigma,$  $V$ contains, up to a compact set, at least one $\Sigma$-region of $\mb_\beta.$  However, $t|_V$ is bounded, contradicting Claim 1 and concluding the proof. \end{proof}

\begin{theorem}[Uniqueness of the Lorentzian Helicoid]

The unique properly embedded $\omega^*$-maximal surface  with connected boundary and infinite rotation number is, up to Lorentzian congruence, the Lorentzian helicoid.

\end{theorem}
\begin{proof}
Let $\mb$ be a properly embedded  $^*$maximal surface  with connected boundary and $\theta_\mb=+\infty.$ Up to isometries, suppose $\Sigma_\infty=\Pi_0.$ From Proposition \ref{pro:porfin},  $h:=t+i t^*:\mb \to \c$ is a injective holomorphic map. Furthermore, since $\Gamma$ is a lightlike arc of mirror symmetry, $t^*|_\Gamma$ is constant ( without loss of generality suppose $t^*|_\Gamma=0$).

Let us see that $\lim_{x \in c_s  \to \infty} t^*(x)=+\infty$ for any $s \in \r,$ where as in the proof of Proposition \ref{pro:porfin} $c_s$ is  the integral curve of ${\cal Y}$ with initial condition $\Gamma(s).$ Indeed, as  $t^*|_{c_s}$ is monotone then the limit $r_s:=\lim_{x \in c_s  \to \infty} t^*(x)$ exists, and without loss of generality, belongs to $]0,+\infty],$ for any $s \in \r.$ In particular
$$\lim_{x \in c_s  \to \infty} h(x)=s+i r_s, \quad s \in \r.$$
Let $V_s$ denote the parabolic region in $\mb$ bounded by $\Gamma([0,s]) \cup c_0 \cup c_s.$ The holomorphic function $h|_{V_s}$ omits infinitely many complex values, hence from Theorem \ref{th:sectorial} the limits along $c_0$ and $c_s$ must coincide for any $s \in \r.$ Therefore,  $r_s=+\infty$ for any $s \in \r,$ proving our assertion.

As a consequence, $h(\mb)=\overline{\u}$ and  $h:\mb \to \overline{\u}$ is a biholomorphism. Furthermore, identifying $\mb$ and $\overline{\u}$ via $h,$ we get $\phi_3=-iB dz,$ $B>0.$

On the other hand, Theorem \ref{th:omit} gives that $g(\overline{\u}) \subset \overline{\d}-\{0\},$ and so $\log(g):\overline{\u} \to \c$ is  well defined. As $|g|^{-1}(1)=\partial(U),$ then $\mbox{Re}(\log(g))$ only vanish on the real axis and $\log(g)|_{\partial(\u)}$ is one to one. Therefore,   $g(z)=e^{a i z+ i b},$ where $a,$ $b \in \r.$ Up to Lorentzian congruence, $\mb$ is  the Lorentzian helicoid, which concludes the proof. \end{proof}


{\bf ISABEL FERNANDEZ,} \newline Departamento de Matemática
Aplicada I \newline Facultad de Informática, Universidad de
Sevilla
\newline 41012 - SEVILLA (SPAIN)\newline
e-mail: isafer@us.es.\\

 {\bf FRANCISCO J. LOPEZ,} \newline
Departamento de Geometr\'{\i}a y Topolog\'{\i}a \newline Facultad
de Ciencias, Universidad de Granada \newline 18071 - GRANADA
(SPAIN) \newline e-mail: fjlopez@ugr.es

\end{document}